\documentclass[oneside, a4paper,reqno]{amsart}
\usepackage{pdfsync}
\usepackage{mathabx}
\usepackage{stmaryrd}
\usepackage{mathrsfs}
\usepackage{amsmath, amsthm, amscd, amssymb, latexsym, eucal}
\usepackage[all]{xy}
 \calclayout \makeatletter
 \calclayout \makeatletter
\def\serieslogo@{} \def\@setcopyright{} \makeatother
\makeatletter
\renewcommand*\env@matrix[1][c]{\hskip -\arraycolsep
  \let\@ifnextchar\new@ifnextchar
  \array{*\c@MaxMatrixCols #1}}
\makeatother
\numberwithin{equation}{section}
\newtheorem{thm}{Theorem}[section]
\newtheorem{cor}[thm]{Corollary}
\newtheorem{lem}[thm]{Lemma}
\newtheorem{prop}[thm]{Proposition}

\theoremstyle{definition}
\newtheorem{defn}[thm]{Definition}
\newtheorem{rem}[thm]{Remark}
\newtheorem{exam}[thm]{Example}

\newcommand{\xr}{\xrightarrow}
\newcommand{\lxr}{\longrightarrow}

\newcommand{\epic}{\twoheadrightarrow}
\newcommand{\monic}{\rightarrowtail}
\newcommand{\ab}{\mathfrak{Ab}}
\newcommand{\A}{\mathscr A}
\newcommand{\B}{\mathscr B}
\newcommand{\C}{\mathscr C}

\newcommand{\E}{\mathscr E}
\newcommand{\F}{\mathcal F}

\newcommand{\R}{\mathcal R}
\newcommand{\T}{\mathcal T}
\newcommand{\U}{\mathcal U}
\newcommand{\V}{\mathcal V}
\newcommand{\W}{\mathcal W}
\newcommand{\X}{\mathcal X}
\newcommand{\Y}{\mathcal Y}
\newcommand{\Z}{\mathcal Z}

\newcommand{\cpt}{\mathsf{cpt}}
\newcommand{\bbot}{\pmb{\bot}}

\DeclareMathOperator*{\Ker}{\mathsf{Ker}}
 \DeclareMathOperator*{\Image}{\mathsf{Im}}
\DeclareMathOperator*{\Coker}{\mathsf{Coker}}

 \DeclareMathOperator{\pd}{\mathsf{pd}}
\DeclareMathOperator*{\id}{\mathsf{id}}
\DeclareMathOperator*{\fd}{\mathsf{fd}}

 \DeclareMathOperator*{\gd}{\mathsf{gl.dim}}

\DeclareMathOperator*{\Mod}{\mathsf{Mod}-\!}

\DeclareMathOperator*{\End}{\mathsf{End}}
 \DeclareMathOperator*{\smod}{\mathsf{mod}-\!}
 \DeclareMathOperator*{\lsmod}{\!-\mathsf{mod}}
\DeclareMathOperator*{\lMod}{\!-\mathsf{Mod}}
\DeclareMathOperator*{\Gor}{\mathsf{G}-\!}
 \DeclareMathOperator*{\umod}{\underline{\mathsf{mod}}-\!}
\DeclareMathOperator*{\uMod}{\underline{\mathsf{Mod}}-\!}
\DeclareMathOperator*{\omod}{\overline{\mathsf{mod}}-\!}
\DeclareMathOperator*{\oMod}{\overline{\mathsf{Mod}}-\!}
\DeclareMathOperator*{\inj}{\mathsf{inj}}
\DeclareMathOperator*{\proj}{\mathsf{proj}}
\DeclareMathOperator*{\Inj}{\mathsf{Inj}}
\DeclareMathOperator*{\Proj}{\mathsf{Proj}}
\DeclareMathOperator*{\Flat}{\mathsf{Flat}}

\newcommand{\GProj}{\operatorname{\mathsf{GProj}}\nolimits}
 \newcommand{\Gproj}{\operatorname{\mathsf{Gproj}}\nolimits}
 \newcommand{\GInj}{\operatorname{\mathsf{GInj}}\nolimits}
 \newcommand{\Ginj}{\operatorname{\mathsf{Ginj}}\nolimits}
 \newcommand{\uGProj}{\operatorname{\underline{\mathsf{GProj}}}\nolimits}
 \newcommand{\oGInj}{\operatorname{\overline{\mathsf{GInj}}}\nolimits}
  \newcommand{\uGproj}{\operatorname{\underline{\mathsf{Gproj}}}\nolimits}
 \newcommand{\oGinj}{\operatorname{\overline{\mathsf{Ginj}}}\nolimits}

\newcommand{\GFlat}{\operatorname{\mathsf{GFlat}}\nolimits}

\DeclareMathOperator*{\ddom}{\mathsf{dom.dim}}
\DeclareMathOperator*{\Sub}{\mathsf{Sub}}

\DeclareMathOperator*{\Add}{\mathsf{Add}}
 
\DeclareMathOperator*{\add}{\mathsf{add}}

\DeclareMathOperator{\Hom}{\mathsf{Hom}}
 
\DeclareMathOperator*{\rd}{\mathsf{res.dim}}
\DeclareMathOperator*{\cd}{\mathsf{cores.dim}}
\DeclareMathOperator*{\Ext}{\mathsf{Ext}}
 
 \DeclareMathOperator*{\colim}{\varinjlim}
 \DeclareMathOperator*{\plim}{\varprojlim}
  \DeclareMathOperator*{\op}{\mathsf{op}}

  \DeclareMathOperator*{\fin}{\,\mathsf{fin}}

\DeclareMathOperator*{\fp}{\mathsf{fp}}

   \DeclareMathOperator*{\filt}{\mathsf{Filt}}
    
     \DeclareMathOperator*{\loc}{\mathsf{Loc}}
     \DeclareMathOperator*{\rp}{\mathsf{rep.dim}}
      \DeclareMathOperator*{\Supp}{\mathsf{Supp}}
       
       \DeclareMathOperator*{\dd}{\mathsf{dom.dim}}
        \DeclareMathOperator*{\uHom}{\underline{\mathsf{Hom}}}
         \DeclareMathOperator*{\uEnd}{\underline{\mathsf{End}}}
          \DeclareMathOperator*{\bmax}{\mathsf{max}}
          \DeclareMathOperator*{\bdim}{\mathsf{dim}}
\newcommand{\La}{\Lambda}
\newcommand{\unA}{\underline A}
\newcommand{\unB}{\underline B}

\newcommand{\unX}{\underline X}

\newcommand{\unf}{\underline f}

\newcommand{\un}{\underline}

\newsavebox{\proofbox}
\savebox{\proofbox}{\begin{picture}(7,7)%
  \put(0,0){\framebox(7,7){}}\end{picture}}

\begin{document}

\title[On Algebras of Finite CM-Type]{On Algebras of Finite Cohen-Macaulay Type}

\author[A. Beligiannis]{Apostolos Beligiannis}
\address{Department of Mathematics, University of Ioannina, GR-45110
Ioannina, Greece}

\email{abeligia@cc.uoi.gr}



\begin{abstract}
We study Artin algebras $\Lambda$ and commutative Noetherian
complete local rings $R$ in connection with the following
decomposition  property of Gorenstein-projective modules:
\[(\dagger) \,\, \,\,\,\,\,\ \text{\em any Gorenstein-projective module is a direct sum of
finitely generated modules.}\] We show that the class of algebras
$\Lambda$ enjoying $\mathbf{(\dagger)}$ coincides with the class of
virtually Gorenstein algebras of finite Cohen-Macaulay type,
introduced in \cite{BR, B:cm}. Thus we solve the  problem stated in
\cite{Chen}. This is proved by characterizing when a resolving
 subcategory is of finite representation type
in terms of decomposition properties of its closure under filtered
colimits, thus generalizing a classical result of Auslander
\cite{Auslander:large} and Ringel-Tachikawa \cite{RT}. In the
commutative case, if  $R$ admits a non-free finitely generated
Gorenstein-projective module, then we show that $R$ is of finite
Cohen-Macaulay type iff $R$ is Gorenstein and satisfies $(\dagger)$.
We also generalize
 a result of Yoshino \cite{Yoshino} by
characterizing when finitely generated modules without extensions
with the ring are Gorenstein-projective. Finally we study the
(stable) relative Auslander algebra of a virtually Gorenstein
algebra of finite Cohen-Macaulay type and, under the presence of a
cluster-tilting object, we give descriptions of the stable category
of Gorenstein-projective modules in terms of the cluster category
associated to the quiver of the stable relative Auslander algebra.
In this setting we show that the cluster category is invariant under
derived equivalences.
\end{abstract}

\maketitle

\setcounter{tocdepth}{1} \tableofcontents

\section{Introduction}
Finitely generated modules of Gorenstein dimension zero over a
Noetherian ring $R$ were introduced by Maurice Auslander in the
mid-sixties, see \cite{Auslander:ENS, ABr}, as a natural
generalization of finite projective modules in order to provide
finer homological invariants of Noetherian rings. Since then they
found important applications in commutative algebra, Algebraic
Geometry, singularity theory and relative homological algebra. Later
on the influential works of Auslander and Buchweitz \cite{ABu,
Buchweitz:tate} gave a more categorical approach to the study of
modules of Gorenstein-dimension zero.

 In general the size and the homological
complexity of the category $\Gproj R$ of finitely generated modules
of Gorenstein dimension zero, also called, depending on the setting,
Gorenstein-projective, totally reflexive, or maximal Cohen-Macaulay,
modules, measure how far the ring $R$ is from being Gorenstein. A
particularly  nice feature of $\Gproj R$ is that its stable version
$\uGproj R$ modulo projectives is triangulated, so stable
homological phenomena concerning the ring take place in $\Gproj R$. For instance $\uGproj R$
admits a full triangulated embedding into the triangulated category of singularities of $R$, in Orlov's sense \cite{Orlov}, which is an equivalence if and only if $R$ is Gorenstein.
On the other hand it was a discovery of many people that also
representation-theoretic properties of the category $\Gproj R$ have
important consequences for the structural shape of the ring, see
Yoshino's book \cite{Yoshino:book} for more details.  For instance
if the Noetherian ring $R$ is commutative, complete, local and
Gorenstein, then $R$ is a simple singularity provided that the
category $\Gproj R$ is of finite representation type, i.e. the set
of isomorphism classes of its indecomposable objects is finite. In
this case  $R$ is said to be of {\em finite Cohen-Macaulay type},
finite CM-type for short.  In many cases also the converse holds.
For instance in the recent paper \cite{CPST}, Christensen,
Piepmeyer, Striuli and Takahashi proved that if a commutative
Noetherian local ring $R$ is of finite CM-type and admits a non-free
finitely generated module of Gorenstein dimension zero, then $R$ is
Gorenstein and an isolated singularity, which is  a simple hypersurface singularity if in addition $R$ is
complete.

In the non-commutative setting, Enochs-Jenda in \cite{EJ}
generalized the notion of finitely generated modules of Gorenstein
dimension zero over a Noetherian ring in order to cover, not
necessarily finitely generated, modules over any ring $R$, under the
name Gorenstein-projective modules; the corresponding full
subcategory is henceforth denoted by $\GProj R$. Enochs and his
collaborators studied extensively relative homological algebra based
on Gorenstein-projective modules, see the book \cite{EJ:book} for
more details and Christensen's book \cite{Christensen} for an
hyperhomological approach. In the context of Artin algebras,
Auslander-Reiten \cite{AR:applications, AR:cm} and Happel
\cite{Happel}, see also \cite{Buchweitz:tate,BR}, studied finitely
generated Gorenstein-projective modules and they established strong
connections with cotilting theory and the shape of the singularity
category.  The results in the
commutative setting suggest that there is a strong link between
representation-theoretic properties of the category $\Gproj R$ and
global structural properties of the ring $R$, especially when $R$ is
Gorenstein. In this connection the class of virtually Gorenstein
Artin algebras was introduced in \cite{BR}, and studied extensively
in \cite{B:cm, BK}, as a convenient common generalization of
Gorenstein  algebras and  algebras of finite representation type.
Note that $\Lambda$ is virtually Gorenstein if and only if any
Gorenstein-projective module is a filtered colimit of finitely
generated Gorenstein-projective modules, so in this case one has a
nice control of the objects of $\GProj R$ in terms of the objects of
$\Gproj R$.

Our aim in this paper is to study representation theoretic aspects
of Gorenstein-projective modules.  In particular we are mainly
interested in the class of Artin algebras and/or commutative
Noetherian complete local rings of finite CM-type, and the structure
of the category of Gorenstein-projective modules and how they are
built from finitely generated ones.

The organization and the main results of the paper are as follows.

In section 2, for completeness and reader's convenience, we collect
several preliminary notions and results that will be useful
throughout the paper and we fix notation. In section 3, generalizing
and extending  classical results of Auslander
\cite{Auslander:large}, Ringel-Tachikawa \cite{RT}, and others, we
give several characterizations of when a resolving subcategory of
the category $\smod\Lambda$ of finitely generated modules over an
Artin algebra $\Lambda$ is of finite representation type in terms of
decomposition properties of its closure under filtered colimits, see
Theorem $3.1$ and Corollary $3.5$. We also point-out some
consequences concerning submodule categories, Weyl modules over
quasi-hereditary algebras and torsion pairs.

In section 4, by applying the results of the previous section, we
present a host of characterizations, more precisely 13 equivalent
conditions, of when an Artin algebra enjoys the decomposition
property $\pmb{(\dagger)}$: any Gorenstein-projective module is a
direct sum of finitely generated modules; see Theorem $4.10$ and
Corollary $4.11$.  In particular we prove that this class of
algebras coincides with the class of virtually Gorenstein algebras
of finite CM-type. Since virtually Gorenstein algebras contain
properly the Gorenstein algebras, we obtain as a direct consequence
an extension of  the main result of Chen \cite{Chen}, and at the
same time we give a complete solution to the problem/conjecture
stated in \cite{Chen}, see Corollary $4.13$. We also characterize
arbitrary Artin algebras of finite CM-type, see Theorem $4.18$.
In the commutative case, we characterize the commutative Noetherian
complete local rings of finite CM-type admitting a non-free finitely
generated Gorenstein-projective module as the Gorenstein rings
satisfying $(\dagger)$; see Theorem $4.20$. As a consequence any
commutative Noetherian
 complete local Gorenstein ring satisfying $(\dagger)$ is a simple singularity.

 Due to their  importance in the homological or
representation-theoretic structure of an Artin algebra, it is useful
to have convenient descriptions of the category $\Gproj\Lambda$. In
this connection Jorgensen-Sega \cite{JS} and Yoshino \cite{Yoshino}
studied the problem of when a finitely generated module without
non-zero extensions with the ring, i.e. a {\em stable} module, is
Gorenstein-projective. By the results of \cite{JS} this is not
always true for certain commutative Noetherian local rings
$(R,\mathfrak{m},k)$. On the other hand Yoshino proved that all
stable $R$-modules are Gorenstein-projective provided that $R$ is
Henselian and the full subcategory ${^{\bot}}R$ of all stable
modules is of finite representation type, see \cite{Yoshino}. Recall
that for a full subcategory $\X \subseteq \smod R$, ${^{\bot}}\X$
denotes the left $\Ext$-orthogonal subcategory of $\X$ consisting of
all modules $A \in \smod R$ such that $\Ext^{n}_{R}(A,X) = 0$,
$\forall X \in \X$, $\forall n \geq 1$. The right $\Ext$-orthogonal
subcategory $\X^{\bot}$ of $\X$ is defined in a dual way. Recently
Takahashi \cite{Takahashi} proved a far reaching generalization: if
${^{\bot}}R$ contains a non-projective module and $k$ admits a right
${^{\bot}}R$-approximation, then $R$ is Gorenstein and then
${^{\bot}}R = \Gproj R$. In section 5 we extend and generalize
Yoshino's results in the non-commutative setting, in fact to a large
class of abelian categories.  In particular  we give a complete
answer to the question of when all stable modules over a virtually
Gorenstein algebra $\Lambda$, or more generally over an algebra for
which $\Gproj\Lambda$ is contravariantly finite, are
Gorenstein-projective, i.e. when ${^{\bot}}\Lambda = \Gproj\Lambda$.
In fact the main result of section 5 shows that this happens if and
only if the category ${^{\bot}}\Lambda \cap (\Gproj\Lambda)^{\bot}$
is of finite representation type (in that case it coincides with the
full subcategory of finite projective modules), see Theorem $5.2$.
This has some consequences for the, still open, Gorenstein Symmetry
Conjecture, \textsf{(GSC)} for short, see \cite{BR, ARS, B:cm},
which predicts that any Artin algebra with finite one-sided
self-injective dimension is Gorenstein. We show that \textsf{(GSC)}
admits the following representation-theoretic formulation:  an Artin
algebra  $\Lambda$ with finite right self-injective dimension is
Gorenstein if and only if $\Gproj\Lambda$ is contravariantly finite
and the full subcategory ${^{\bot}}\Lambda \cap
(\Gproj\Lambda)^{\bot}$ of $\smod\Lambda$  is of finite
representation type. On the other hand we show that if $R$ is a
commutative Noetherian complete local Gorenstein ring, then all
infinitely generated stable modules are Gorenstein-projective if and
only if $R$ is Artinian, so the infinitely generated versions of the
above results rarely hold in the commutative setting.

In section 6, motivated by Auslander's homological theory of
representation-finite Artin algebras $\Lambda$, see
\cite{Auslander:queen}, we consider relative Auslander algebras
$\mathsf{A}(\X)$, defined as endomorphism rings of representation
generators, i.e. modules containing as direct summands all types of
indecomposable modules, of contravariantly finite resolving
subcategories $\X \subseteq \smod\Lambda$. Recall that Auslander
proved that $\mathsf{A}(\X)$ has
global dimension at most $2$, when $\X = \smod\Lambda$ for an Artin algebra $\Lambda$
of finite representation type. In section $6$ we compute the global
dimension of $\mathsf{A}(\X)$  in terms of the resolution dimension
of $\smod\Lambda$ with respect to $\X$. In particular any natural
number or infinity can occur. As a consequence we characterize Gorenstein abelian categories of finite Cohen-Macaulay type, thus generalizing and,  at the same time, giving a simple proof to the main result of \cite{LiZhang}, see Corollary 6.13. Then we concentrate on relative, resp.
stable, Auslander algebras $\mathsf{A}(\Gproj\Lambda)$,  resp.
$\mathsf{A}(\uGproj\Lambda)$, of Artin algebras $\Lambda$ of finite
CM-type which we call (stable) Cohen-Macaulay Auslander algebras. We
show that the global dimension of $\mathsf{A}(\Gproj\Lambda)$ is
bounded by the Gorenstein dimension $\Gor\bdim\Lambda$ of $\Lambda$,
so it is finite if and only if $\Lambda$ is Gorenstein, whereas
$\mathsf{A}(\uGproj\Lambda)$ is always self-injective. As a
consequence $\bmax\{2,\Gor\bdim\Lambda\}$ is an upper bound for the
dimension, in the sense of Rouquier \cite{Rouquier:Kth}, of the
derived category ${\bf D}^{b}(\smod\Lambda)$ of $\Lambda$.  In the
commutative case we show that the global dimension of the
Cohen-Macaulay Auslander algebra $\mathsf{A}(\Gproj R)$ of a
commutative Noetherian complete local ring $R$ of finite CM-type is
at most $\bmax\{2,\bdim R\}$, and exactly $\bdim R$ if $\bdim R \geq
2$, so this number is an upper bound of the dimension of ${\bf
D}^{b}(\smod R)$. Finally we give an application to the periodicity
of Hochschild (co)homology of the stable Cohen-Macaulay Auslander
algebra $\mathsf{A}(\uGproj\Lambda)$
 which is based on
results of Buchweitz \cite{Buchweitz}.

The final section 7 of the paper is devoted to the study of finitely
generated Gorenstein-projective rigid modules, i.e. modules without
self-extensions, or cluster tilting Gorenstein-projective modules,
in the sense of Keller-Reiten \cite{KR} or Iyama \cite{Iyama}, over
an Artin algebra $\Lambda$ in connection with the property that
$\Lambda$  is of finite CM-type. We show that Cohen-Macaulay
finiteness of $\Lambda$ is intimately  related to the representation
finiteness of the stable endomorphism ring $\uEnd_{\Lambda}(T)$ of a
Gorenstein-projective rigid module $T$, and in case $T$ is cluster
tilting the two finiteness conditions are equivalent. As a
consequence we show that if $\uEnd_{\Lambda}(T)$ is
of finite CM-type, then the representation dimension
$\rp\uEnd_{\Lambda}(T)$ of $\uEnd_{\Lambda}(T)$, in the sense of
Auslander \cite{Auslander:queen}, is at most $3$, and is exactly $3$ if and only if $\Lambda$
is of infinite CM-type. Finally by using that the stable
triangulated category $\uGproj\Lambda$ of finitely generated
Gorenstein projective modules over a virtually Gorenstein algebra
$\Lambda$ admits Serre duality \cite{B:cm} and recent results of
Amiot \cite{Amiot} and Keller-Reiten \cite{KR,KR:acyclic}, we give
descriptions of $\uGproj\Lambda$, in case the latter is
$2$-Calabi-Yau in the sense of \cite{KR} and $\uEnd_{\Lambda}(T)$
has finite global dimension, in terms of the cluster category
$\C_{Q}$ associated to the quiver $Q$ of $\uEnd_{\Lambda}(T)$. In
particular we show that, under the above assumptions, the cluster
category $\C_{Q}$ is invariant under derived equivalences.

\section{Preliminaries}

In this section we fix notation and, for completeness and reader's
convenience, recall some well-known notions and collect several
preliminary results which will be used throughout the paper.

\subsection{Locally Finitely Presented Categories} Let $\A$ be an
additive category with filtered colimits. An object $X$ in $\A$ is
called {\bf finitely presented} if the functor $\A(X,-) : \A \lxr
\ab$ commutes with filtered colimits. We denote by $\fp\A$ the full
subcategory of $\A$ consisting of the finitely presented objects. If
$\U$ is a full subcategory of $\A$ we denote by $\colim\U$ the
closure of $\U$ in $\A$ under filtered colimits, that is $\colim\U$
is the full subcategory of $\A$ consisting of all filtered colimits
of objects from $\U$.  We recall that $\A$ is called {\bf locally
finitely presented} if $\fp\A$ is skeletally small and $\A =
\colim\fp\A$.

For instance let  $\C$ be a skeletally small additive category with
split idempotents.  Then the category $\Mod\C$ of contravariant
additive functors $\C^{\op} \lxr \ab$ is locally finitely presented
and $\fp\Mod\C = \smod\C$ is the full subcategory of coherent
functors over $\C$. Recall that an additive functor $F \colon \C \to
\ab$ is called {\bf coherent} \cite{Auslander:coherent}, if there
exists an exact sequence $\C(-,X) \lxr \C(-,Y) \lxr F \lxr 0$. The
categories $\C\lMod$ and $\C\lsmod$ are defined by $\C\lMod = \Mod\C^{\op}$,
and $\C\lsmod = \smod\C^{\op}$ respectively. Note that $\smod\C$ is abelian if
and only if $\C$ has weak kernels. Recall that $\C$ has weak kernels
if any map $X_{2} \lxr X_{3}$ in $\C$ can be extended to a diagram
$X_{1} \lxr X_{2} \lxr X_{3}$ in $\C$ such that the induced sequence
of functors $\C(-,X_{1}) \lxr \C(-,X_{2}) \lxr \C(-,X_{3})$ is
exact; in this case the map $X_{1} \lxr X_{2}$ is called a {\bf weak
kernel} of $X_{2} \lxr X_{3}$. {\em Weak cokernels} are defined
dually. Note that a locally finitely presented category $\A$ has
products iff $\fp\A$ has weak cokernels and  $\A$ is abelian iff
$\A$ is a Grothendieck category, see \cite{WCB}. In this case $\A$
is called {\bf locally finite} if  any finitely presented object has
finite Jordan-H\"{o}lder length, equivalently $\fp \A$ is Artinian
and Noetherian.

\subsubsection{Representation Categories, Purity and the Ziegler Spectrum} Let $\A = \colim \fp\A$ be a
locally finitely presented additive category. The {\bf
representation category} $\mathcal L(\A)$ of $\A$ is defined to be
the functor category $\mathcal L(\A) := \Mod\fp\A$. It is well-known
that $\mathcal L(\A)$ reflects important representation theoretic
properties of $\A$ via the {\bf representation functor}
\[
\mathsf{H} \, : \, \A \lxr \mathcal L(\A), \,\,\,\,\ A \,
\longmapsto \, \mathsf{H}(A) = \A(-,A)|_{\fp\A}
\]
which induces an equivalence $\mathsf{H} \ : \ \A \,
\stackrel{\approx}{\lxr}\, \Flat(\fp\A)$, where for an additive
category $\C$, $\Flat\C$ denotes the full subcategory of $\Mod\C$
consisting of all flat functors; recall that $F \colon
\C^{\op} \lxr \ab$ is called {\bf flat} if $F$ is a filtered colimit of
representable functors. Assume now that $\A$ has products. A
sequence $(*): 0\lxr A \lxr B \lxr C \lxr 0$ in $\A$ is called {\bf
pure-exact} if the induced sequence $0\lxr \mathsf{H}(A) \lxr
\mathsf{H}(B) \lxr \mathsf{H}(C) \lxr 0$ is exact in $\mathcal
L(\A)$. An object $E$ in $\A$ is called {\bf pure-projective}, resp.
{\bf pure-injective}, if for any pure-exact sequence $(*)$ as above,
the induced sequence $0 \lxr \A(E,A)\lxr \A(E,B) \lxr \A(E,C) \lxr
0$, resp. $0 \lxr \A(C,E)\lxr \A(B,E) \lxr \A(A,E) \lxr 0$ is exact.
It is easy to see that the full subcategory of $\A$ consisting of
the pure-projective objects coincides with $\Add(\fp\A)$. Recall
that if $\X$ is a full subcategory of $\A$, then $\Add \X$, resp.
$\add\X$, denotes the full subcategory of $\A$ consisting of the
direct summands of all small, resp. finite, coproducts of objects of
$\X$. The locally finitely presented category $\A$ is called {\bf
pure-semisimple} if any pure-exact sequence in $\A$ splits.
Equivalently  $\A = \Add(\fp\A)$, or the representation category
$\mathcal L(\A)$ is {\bf perfect}, i.e. any flat functor is
projective.

We denote by ${\bf Zg}(\A)$ the collection of the isomorphism
classes of indecomposable pure-injective objects of $\A$. It is
well-known that ${\bf Zg}(\A)$ is a small set and admits a natural
topology, the {\bf Ziegler topology}, defined as follows. Let $\Phi$
be a collection of maps between finitely presented objects. An
object $E$ is called $\Phi$-{\em injective} if for any map $\phi
\colon X\lxr Y$ in $\Phi$, the induced map $\A(\phi,E)$ is
surjective. We denote by $\Inj(\Phi)$ the collection of
$\Phi$-injective objects of $\A$. Then the subsets of ${\bf Zg}(\A)$
of the form ${\bf U}_{\Phi} := {\bf Zg}(\A) \bigcap \Inj(\Phi)$
constitute the closed sets of the Ziegler topology of ${\bf
Zg}(\A)$, i.e. a subset ${\bf U}$ of ${\bf Zg}(\A)$ is
Ziegler-closed if ${\bf U} = {\bf U}_{\Phi}$ for some collection of
maps $\Phi$ in $\fp(\A)$. A subset ${\bf C}$ of ${\bf Zg}(\A)$ is
Ziegler-open if its complement ${\bf Zg}(\A)\setminus {\bf C}$ is
Ziegler-closed. We refer to \cite{Krause:memoirs} for more details.

\subsection{Contravariantly Finite Subcategories, Filtrations, and Cotorsion Pairs}
Let $\A$ be an additive category. If $\V$ is a full subcategory of
$\A$, then a map $f : A \xr{} B$ in $\A$ is called $\V$-{\bf epic}
if the map $\A(\V,f) : \A(\V,A) \xr{} \A(\V,B)$ is surjective.
 $\V$ is called {\bf contravariantly finite} if there exists a $\V$-epic $f_{A} : V_{A} \xr{} A$ with $V_{A}$
in $\V$. In this case $f_{A}$, or $V_{A}$,  is called a {\bf right} $\V$-{\bf
approximation} of $A$.  {\em Covariantly finite} subcategories,
$\V$-{\em monics} and {\em left} $\V$-{\em approximations} are
defined dually. $\V$ is called {\em functorially finite} if it is
both contravariantly and covariantly finite. A map $f : V \lxr A$ is
called {\em right minimal} if any endomorphism $\alpha$ of $V$ such
that $\alpha \circ f = f$ is an automorphism. A right minimal right
$\V$-approximation is called a {\em minimal right
$\V$-approximation}. {\em Left minimal} maps and {\em minimal left
approximations} are defined dually. The objects of the {\bf stable
category} $\A/\V$ of $\A$ modulo $\V$ are the objects of $\A$; the
morphism spaces are defined by $\A/\V(A,B) := \A(A,B)/\A_{\V}(A,B)$
where $\A_{\V}(A,B)$ is the subgroup of $\A(A,B)$ consisting of all
maps factorizing through an object from $\V$.  We have the natural
projection functor $\pi : \A \lxr \A/\V$, $\pi(A) = \unA$ and
$\pi(f) := \un{f}$, which is universal for additive functors out of
$\A$ killing the objects of $\V$.

Now let  $\A$ be  abelian. We denote by $\Proj\A$, resp.
$\Inj\A$, the full subcategory of the projective, resp. injective,
objects of $\A$. Let $\U$ be a full subcategory of $\A$. The stable
categories $\U/\Proj\A$ and $\U/\Inj\A$ are denoted by $\un\U$ and
$\overline{\U}$ respectively.  $\U$ is called {\bf resolving} if
$\U$ contains the projectives and is closed under extensions and
kernels of epimorphisms. {\em Coresolving} subcategories are defined
in a dual way. Typical examples of (co)resolving subcategories are
the {\bf left} and {\bf right} $\Ext$-{\bf orthogonal}
subcategories ${^{\bot}}\U$ and $\U^{\bot}$ of $\U$ defined as follows:
\[{^{\bot}}\U := \big\{A \in \A \,\ | \,\,  {\Ext}^{n}_{\A}(A,U) = 0, \, \forall n\geq 1, \, \forall U \in
\U \big\}\]
\[\U^{\bot} :=
\big\{A \in \A \,\ | \,\ {\Ext}^{n}_{\A}(U,A) = 0, \, \forall n\geq
1, \, \forall U \in \U \big\}\]

\begin{defn} A pair $(\X,\Y)$ of full subcategories of $\A$ is called a
{\bf cotorsion pair} if: $\Ext^{n}(\X,\Y) = 0$, $\forall n \geq 1$,
and for any object $A\in \A$, there exist short exact sequences
\begin{equation}
  0\,\ \lxr Y_{A}  \,\ \stackrel{g_{A}}{\lxr} \,\  X_{A} \,\ \stackrel{f_{A}}{\lxr} \,\  A \,\ \lxr \,\ 0
\,\,\,\,\,\ \text{and}\,\,\,\,\,\
  0 \,\ \lxr\,\  A \,\ \stackrel{g^{A}}{\lxr} \,\ Y^{A} \,\, \stackrel{f^{A}}{\lxr} \,\
  X^{A} \,\ \lxr \,\, 0
\end{equation}
where $X^{A}, X_{A} \in \X$ and $Y_{A}, Y^{A} \in \Y$. The {\bf
heart} of the cotorsion pair $(\X,\Y)$ is the full subcategory
$\X\cap \Y$.
\end{defn}

Cotorsion pairs in the above sense are also known in the literature
as complete hereditary cotorsion pairs. It follows easily that if
$(\X,\Y)$ is a cotorsion pair in $\A$, then $\X$ is contravariantly
finite resolving and $\X^{\bot} = \Y$, and $\Y$ is covariantly
coresolving and $\X = {^{\bot}}\Y$. Moreover the heart $\omega =
\X\cap \Y$ is an Ext-injective cogenerator of $\X$ and an
Ext-projective generator of $\Y$. Recall that a full subcategory
$\omega$ of $\X$ is called an {\bf Ext-injective cogenerator} of
$\X$ provided that for any object $X \in \X$, there exists a short
exact sequence $0 \lxr X \lxr T \lxr X^{\prime} \lxr 0$, where $T$
lies in $\omega$, $X^{\prime}$ lies in $\X$, and $\Ext^{1}_{\A}(\X,T) =
0$. {\em Ext-projective generators} are defined
dually. A cotorsion pair $(\X,\Y)$ in $\A$ is called {\bf
projective} if  $\X \cap \Y = \Proj\A$. Note that if $\A$ has enough
projectives and $(\X,\Y)$ is a projective cotorsion pair in $\A$,
then the stable category $\un\X$ is triangulated  and the
subcategory $\Y$ is {\em thick} in the sense that $\Y$ is closed
under extensions, direct summands, kernels of epimorphisms and
cokernels of monomorphisms. Moreover the canonical map
$\Ext^{n}_{\A}(X,A) \lxr \un\A(\Omega^{n}X,A)$ is invertible,
$\forall n \geq 1$, $\forall X \in \X$, $\forall A \in \A$. Here
$\Omega$  is the Heller's loop space functor on the stable category
$\un\A$. {\em Injective} cotorsion pairs are defined dually; we
refer to \cite{BR} for more details.
 If the abelian category $\A$ is cocomplete, then  a cotorsion pair
$(\X,\Y)$ in $\A$ is called {\bf smashing} if $\Y$ is closed under
all small coproducts. If $\A$ is in addition locally finitely
presented, then the cotorsion pair $(\X,\Y)$ is called {\bf of
finite type} if there is a subset $\mathcal S$ of $\fp\A$ such that
$\Y = \mathcal S^{\bot}$.

If $\X$ is a full subcategory of $\A$, we denote by $\rd_{\X}A$ the
$\X$-{\em resolution dimension} of $A \in\A$ which is defined as the
minimal number $n$ such that there exists an exact sequence $0 \lxr
X_{n} \lxr X_{n-1} \lxr \cdots \lxr X_{0} \lxr A \lxr 0$, where the
$X_{i}$ lie in $\X$. If no such number $n$ exists, then we set
$\rd_{\X}A = \infty$. The $\X$-{\bf resolution dimension}
$\rd_{\X}\A$ of $\A$ is defined by $\rd_{\X}\A = \sup\{\rd_{\X}A\, |
\, A\in \A\}$. The {\em coresolution dimensions} $\cd_{\X}A$ and
$\cd_{\X}\A$ are defined dually. Finally we denote by $\filt\X$ the
full subcategory of $\A$ consisting of all direct summands of
objects $A$ admitting  a finite filtration
\begin{equation}
 0 = A_{n+1}\, \subseteq\, A_{n} \, \subseteq\, A_{n-1}\, \subseteq \,\
 \cdots \,
 \cdots \,\
 \subseteq \, A_{2}\, \subseteq A_{1} \, \subseteq\, A_{0} = A
\end{equation}
such that the subquotients $A_{i}/A_{i+1}$, for $i = 0,1,\cdots,n$,
lie in $\X$.

\subsection{Compactly Generated Triangulated Categories} Let $\T$ be
a triangulated category which admits all small coproducts. An object
$T$ in $\T$ is called {\bf compact} if the functor $\T(T,-) \colon
\T \lxr \ab$ preserves all small coproducts. We denote by
$\T^{\cpt}$ the full subcategory of $\T$ consisting of all compact
objects. It is easy to see that $\T^{\cpt}$ is a {\em thick}
subcategory of $\T$, that is, $\T^{\cpt}$ is closed under direct
summands and if two of the objects in a triangle $X \lxr Y \lxr Z
\lxr X[1]$ in $\T$ lie in $\T^{\cpt}$, then so does the third. The
triangulated category $\T$ is called {\bf compactly generated} if
$\T^{\cpt}$ is skeletally small and {\em generates} $\T$ in the
following sense: an object $X$ in $\T$ is zero provided that
$\T(C,X) = 0$, for any compact object $C \in \T^{\cpt}$.

\subsubsection{Representation Categories} Let $\T$ be a
compactly generated triangulated category. Then clearly the category
$\Flat(\T^{\cpt})$ of flat contravariant functors
$\{\T^{\cpt}\}^{\op} \lxr \ab$, which coincides with the category of
cohomological functors over $\T^{\cpt}$, has products. We call
$\Mod\T^{\cpt} = \mathcal L(\Flat(\T^{\cpt}))$ the {\bf
representation category} of $\T$ and we denote it by $\mathcal L(\T)
:= \Mod\T^{\cpt}$. In contrast to the case of locally finitely presented categories,
the representation functor
\[
\mathsf{H} \, \colon \, \T \, \lxr \, \mathcal L(\T) =
\Mod\T^{\cpt}, \,\,\,\,\, \mathsf{H}(A) = \T(-,A)|_{\T^{\cpt}}
\]
is in general not fully faithful, due to the presence of {\bf phantom maps},
i.e. maps $A \lxr B$ in $\T$ such that the induced map $\T(X,A) \lxr
\T(X,B)$ is zero for any compact object $X$.  In fact $\mathsf{H}$
is faithful if and only if $\mathsf{H}$ induces an equivalence
$\mathsf{H} : \T \stackrel{\approx}{\lxr} \Flat(\T^{\cpt})$ if and
only if $\T$ is {\bf phantomless}, i.e. there are no non-zero phantom
maps in $\T$; equivalently $\T$ is a (pure-semisimple) locally
finitely presented category. However $\mathsf{H}$ is homological,
has image in $\Flat(\T^{\cpt})$, and detects zero objects:
$\mathsf{H}(A) = 0$ if and only if $A = 0$. In general for the
representation functor $\mathsf{H} : \T \lxr \Flat(\T^{\cpt})$ we
have, in general non-reversible, implications: $\mathsf{H}$ is
faithful $\Rightarrow$ $\mathsf{H}$ is full $\Rightarrow$
$\mathsf{H}$ is surjective on objects, see \cite{B:3cats} for
details.

\subsubsection{Purity and the Ziegler Spectrum}
A triangle $(T):  A \lxr B \lxr C \lxr A[1]$ in  $\T$ is called {\bf
pure} if the induced sequence $0\lxr\mathsf{H}(A) \lxr \mathsf{H}(B)
\lxr \mathsf{H}(C) \lxr 0$ is exact in $\mathcal L(\T)$. An object
$E$ in $\T$ is called {\bf pure-projective}, resp. {\bf
pure-injective}, if for any pure-triangle $(T)$ as above, the
induced sequence $0 \lxr \T(E,A)\lxr \T(E,B) \lxr \T(E,C) \lxr 0$,
resp. $0 \lxr \T(C,E)\lxr \T(B,E) \lxr \T(A,E) \lxr 0$ is exact.
Then $\T$ is called {\bf pure-semisimple} if any pure-triangle
splits. This happens if and only if any object of $\T$ is
pure-projective, i.e. $\T = \Add\T^{\cpt}$, if and only if any
object of $\T$ is pure-injective if and only if $\T$ is phantomless.
Note that the collection ${\bf Zg}(\T)$ of isoclasses of
indecomposable pure-injective objects of $\T$ form a small set,
called the {\em Ziegler spectrum}, and admits a natural topology,
the  Ziegler topology, as in the case of locally finitely presented
categories. We refer to \cite{B:art} for more details.

\medskip

{\bf Conventions and Notations.} If $\Lambda$ is a ring, we denote
by $\Mod\Lambda$ the category of right $\Lambda$-modules. Left
$\Lambda$-modules are treated as right $\Lambda^{\op}$-modules. The
category of finitely presented right $\Lambda$-modules is denoted by
$\smod\Lambda$. The category of projective, resp. injective, modules
is denoted by $\Proj\Lambda$, resp. $\Inj\Lambda$.   The stable
category of $\Mod\Lambda$, resp. $\smod\Lambda$, modulo projectives
is denoted by $\uMod\Lambda$, resp. $\umod\Lambda$, and the stable
category of $\Mod\Lambda$, resp. $\smod\Lambda$, modulo injectives
is denoted by $\oMod\Lambda$, resp. $\omod\Lambda$. If $\X$ is a
subcategory of $\Mod\Lambda$, we use the notations $\un\X$ and
$\overline{\X}$ for the stable categories of $\X$ modulo projectives
and injectives respectively. Moreover we set $\X^{\fin} := \X \cap
\smod\Lambda$. In particular we set $(\Proj\Lambda)^{\fin} =
\proj\Lambda$ and $(\Inj\Lambda)^{\fin} = \inj\Lambda$. If $\X
\subseteq \smod\Lambda$, then  we denote by ${^{\bbot}}\X$ the full
subcategory $\{A \in \Mod\Lambda \, | \, \Ext^{n}_{\Lambda}(A,X) =
0, \,\, \forall X \in \X\}$, and dually for $\X{^{\bbot}}$.  Then
${^{\bot}}\X = ({^{\bbot}}\X)^{\fin} = {^{\bbot}}\X \cap
\smod\Lambda$.  If $\Lambda$ is an Artin $R$-algebra over a
commutative Artin ring $R$, we denote by  $\mathsf{D}$ the usual
duality of Artin algebras which is given by $\Hom_{R}(-,E) \colon
\Mod\Lambda \lxr \Mod\Lambda^{\op}$ where $E$ is the injective
envelope of $R/\mathsf{rad}(R)$. The Jacobson radical of $\Lambda$
is denoted by $\mathfrak{r}$. For all unexplained notions and
results  concerning the representation theory of Artin algebras we
refer to the book \cite{ARS}.

The additive categories considered in this paper are assumed to
admit finite direct sums and all full subcategories are closed under
isomorphisms.  The composition of morphisms in a given category is
meant in the diagrammatic order: the composition of  $f \colon A
\lxr B$ and $g \colon B \lxr C$ is denoted by $f \circ g$. However functors are composed in the usual anti-diagrammatic way.

 \section{Contravariantly Finite Resolving Subcategories}
Recall that an additive category $\A$ is called a {\em
Krull-Schmidt} category if any object of $\A$ is a finite coproduct
of indecomposable objects and any indecomposable object has local
endomorphism ring.
 Assume that $\A$ is a skeletally small
Krull-Schmidt category. We denote by $\mathsf{Ind}\A$ the set of
isoclasses of indecomposable objects of $\A$. We say that $\A$ is of
{\em finite representation type} if the set $\mathsf{Ind}\A$ is
finite. In this case $\A$ admits a {\em representation generator},
i.e. an object $T$ such that $\A = \add T$. If $F : \A \lxr \ab$ is
a covariant or contravariant additive functor,  then the {\em
support} of $F$ is defined by $\Supp F = \{A \in \mathsf{Ind}\A \,
| \, F(A) \neq 0\}$.

Let $\Lambda$ be an Artin algebra. Our aim in this section is to
prove the following result which gives
 characterizations for a contravariantly finite resolving
 subcategory $\X$ of $\smod\Lambda$
 to be of finite representation type in terms of decomposition properties of
  the closure of $\X$ under filtered colimits.

\begin{thm} Let $\X$ be a contravariantly finite resolving
subcategory of $\smod\Lambda$, and let $\A = \colim\X$. Then the
following statements are equivalent.
\begin{enumerate}
\item $\mathcal X$ is of finite representation type.
\item Any module in $\A$ is a direct sum of  finitely generated
modules.
\item The representation category $\mathcal L(\A) = \Mod\X$ is locally finite.
\item The functor $(-,X_{\Lambda/\mathfrak{r}}) \in \Mod\X$ has finite length.
\item Any indecomposable module in $\A$ is finitely generated.
\item Any module in $\A$ is a direct sum of indecomposable
modules.
\item $\A$ is equivalent to the category of projective modules over
an Artin algebra.
\item $\X$ is equivalent to the category of finitely generated projective modules over
an Artin algebra.
\end{enumerate}
\begin{proof} Consider  the representation functor $\mathsf{H} \colon \A \lxr
\mathcal L(\A) = \Mod\X$, $\mathsf{H}(A) =
\Hom_{\Lambda}(-,A)|_{\X}$, which induces an equivalence $\A \approx
\Flat(\X)$. By a result of Krause-Solberg \cite{KS:appl}, the
category $\X$ is also covariantly finite. It follows that $\X$ is a
functorially finite extension closed subcategory of $\smod\Lambda$.
Therefore $\X$ is a dualizing $R$-subvariety of $\smod\Lambda$ in
the sense of \cite{AS:subcategories}. In particular $\X$ has left
and right almost split maps, in the sense of
\cite{AS:subcategories}, and there exists an equivalence
$(\smod\X^{\op})^{\op} \approx \smod\X$, see \cite{AR:dual}.

(i) $\Rightarrow$ (ii) Assuming (i), let $T$ be a representation
generator of $\X$, i.e. $\X = \add T$. If $\Gamma :=
\End_{\Lambda}(T)^{\op}$, then $\mathcal L(\A) = \Mod\X \approx
\Mod\Gamma$. Since $\Gamma$ is an Artin algebra we have $\A =
\colim\X \approx \Flat(\proj\Gamma) \approx \Proj\Gamma$ and
therefore any module in $\A$ is a direct sum of (indecomposable)
finitely generated modules.

(ii) $\Rightarrow$ (iii) $\Rightarrow$ (iv)  By hypothesis the
locally finitely presented category $\A = \colim\X$ is
pure-semisimple. By the several object version of the well-known
result of Bass, see \cite{WCB}, this implies that for any chain
\begin{equation}
X_{0} \,\stackrel{f_{0}}{\lxr} \, X_{1} \,\stackrel{f_{1}}{\lxr} \,
X_{2} \, \lxr \,\, \cdots \cdots \,\ \lxr \, X_{n}
\,\stackrel{f_{n}}{\lxr} \, X_{n+1} \,\lxr \, \, \cdots \cdots
\end{equation}
of non-isomorphisms between indecomposable modules in $\X$, the
composition $f_{0} \circ f_{1} \circ \cdots \circ f_{m}$ is zero for
large $m$. We claim that any $0 \neq F \in \X\lMod$ has a simple
subfunctor. Assuming  that this is not the case, we construct a
chain of maps as in (3.1) as follows. Since $F \neq 0$, there exists
an indecomposable module $X_{0}$ in $\X$ such that $F(X_{0}) \neq
0$. This means that there is a non-zero map $\alpha_{0} \colon
(X_{0},-) \lxr F$; in particular the functor $F_{0} =
\Image\alpha_{0}$ is non-zero. Let $\alpha_{0} = \varepsilon_{0}
\circ i_{0} : (X_{0},-) \epic F_{0} \monic F$ be the canonical
factorization of $\alpha_{0}$. By hypothesis $F_{0}$ is not simple
and therefore $F_{0}$ contains a proper non-zero subfunctor $\mu_{0}
\colon F_{1} \monic F_{0}$. Since $F_{1}$ is non-zero,
there exists an indecomposable module $X_{1}$ in $\X$ and a non-zero
map $\alpha_{1} \colon (X_{1},-) \lxr F_{1}$.
 Since $(X_{1},-)$ is a projective
functor, the composition $\alpha_{1} \circ \mu_{0} \colon (X_{1},-)
\lxr F_{1} \lxr F_{0}$ factors through the epimorphism
$\varepsilon_{0} \colon (X_{0},-) \epic F_{0}$, i.e. $(f_{0},-) \circ
\varepsilon_{0} = \alpha_{1} \circ \mu_{0}$ for some map $f_{0}
\colon X_{0} \lxr X_{1}$. Observe that $f_{0}$ is non-zero since
$\alpha_{1} \neq 0$. Let $F_{2} = \Image\alpha_{1}$ and let
$\alpha_{1} = \varepsilon_{1} \circ i_{1} : (X_{1},-) \epic F_{2}
\monic F_{1}$ be the canonical factorization of $\alpha_{1}$. Since
the non-zero functor $F_{2}$ is a subfunctor of $F$, by hypothesis
it is not simple. Hence $F_{2}$ contains a non-zero proper
subfunctor $\mu_{1} \colon F_{3} \monic F_{2}$. Then as
above there exists an indecomposable module $X_{2}$ in $\X$ and a
non-zero map $\alpha_{2} \colon (X_{2},-) \lxr F_{3}$. Since the
functor $(X_{2},-)$ is projective, the map $\alpha_{2}\circ \mu_{1}$
factors through $\varepsilon_{1}$: $(f_{1},-) \circ \varepsilon_{1}
= \alpha_{2}\circ \mu_{1}$ for some map $f_{1} \colon X_{1} \lxr
X_{2}$. Observe that $(f_{0}\circ f_{1},-) \circ \alpha_{0} =
(f_{1},-) \circ (f_{0},-) \circ \alpha_{0} = \alpha_{2} \circ \mu_{1}\circ i_{1}
\circ \mu_{0} \circ i_{0}$. Then $f_{0} \circ f_{1} \neq 0$ since
$\alpha_{2}$ is non-zero and $\mu_{1}\circ i_{1} \circ \mu_{0} \circ
i_{0}$ is a composition of monics. Pictorially then we have the
following diagram:
 {\small\[\xymatrix{
   \cdot \, \,\  \cdot \,\,\  \cdot        & (X_3,-) \ar[dd]^{\alpha_3} \ar@{-->}[rr]^{(f_2,-)}  &
   &  (X_2,-) \ar[dd]^{\alpha_2} \ar@{-->}[rr]^{(f_1,-)} \ar@{->>}[dl]^{\varepsilon_2} &   & (X_1,-) \ar[dd]^{\alpha_1}
   \ar@{-->}[rr]^{(f_0,-)} \ar@{->>}[dl]^{\varepsilon_1}
  &  & (X_0,-) \ar@{->>}[dl]^{\varepsilon_0}  \ar[dd]^{\alpha_0}  \\
    \cdot \, \,\  \cdot \,\,\  \cdot &   &  F_4 \ar@{^{(}->}[dr]_{i_2} &  &  F_2  \ar@{^{(}->}[dr]_{i_1} &
     &  F_0  \ar@{^{(}->}[dr]_{i_0} &  \ar[d]^{} \\
   \cdot \, \,\  \cdot \,\,\  \cdot &  F_5 \ar@{^{(}->}[ur]_{\mu_2} &   & F_3 \ar@{^{(}->}[ur]_{\mu_1}  &
    &  F_1 \ar@{^{(}->}[ur]_{\mu_0} &   & F   }\]}
Continuing in this way we produce a chain of maps between
indecomposable modules in $\X$ as in (3.1) with non-zero composition
$f_{0} \circ f_{1} \circ \cdots \circ f_{n} \colon X_{0} \lxr
X_{n+1}$, $\forall n \geq 0$. This contradiction shows that any
non-zero functor $\X \lxr \ab$ contains a simple subfunctor. On the
other hand since $\X$ has left and right almost split maps, by a
well-known result of Auslander \cite{Auslander:II} it follows that
the functor category $\X\lMod$ is locally finite, i.e. the coherent
functors $\X\lsmod$ form a length category. Then the existence of a
duality $(\X\lsmod)^{\op} \approx \smod\X$ implies that $\smod\X$ is
a length category, i.e. $\Mod\X$ is also locally finite. If this
holds, then any representable functor, in particular
$(-,X_{\Lambda/\mathfrak{r}})$, is of finite length.

(iv) $\Rightarrow$ (i) By a result of Auslander-Reiten
\cite{AR:applications} we have $\X =
\filt(X_{\Lambda/\mathfrak{r}})$, that is, $\X$  consists of the
direct summands of modules $A$ admitting  a finite filtration
\begin{equation}
 0 = A_{t+1}\, \subseteq\, A_{t} \, \subseteq\, A_{t-1}\, \subseteq \,\
 \cdots \,
 \cdots \,\
 \subseteq \, A_{2}\, \subseteq A_{1} \, \subseteq\, A_{0} = A
\end{equation}
such that the subquotients $A_{i}/A_{i+1}$, for $0\leq i \leq t$,
lie in the set $\X(\mathcal S) =
\{X_{S_{1}},X_{S_{2}},\cdots,X_{S_{m}}\}$, where $\mathcal S =
\{S_{1},S_{2},\cdots,S_{m}\}$ is the set of isoclasses of the simple
$\Lambda$-modules and $X_{S_{i}} \lxr S_{i}$ is the minimal right
$\X$-approximation of $S_{i}$. Note that $\add\X(\mathcal S) =
\add(X_{\Lambda/\mathfrak{r}})$.  By induction on the length $t$ of
the filtration $(3.2)$ it follows directly that $A \in
\mathsf{Supp}(-,X_{\Lambda/\mathfrak{r}})$ and consequently $\mathsf{Ind}(\X)
\subseteq \mathsf{Supp}(-,X_{\Lambda/\mathfrak{r}})$. Hence
\begin{equation}
\mathsf{Ind}(\X) = \mathsf{Supp}(-,X_{\Lambda/\mathfrak{r}})
\end{equation}
Since the functor $(-,X_{\Lambda/\mathfrak{r}})$ has finite length,
there exists a composition series in $\Mod\X$:
\begin{equation}
0 = F_{n+1} \, \subseteq\, F_{n} \, \subseteq\, F_{n-1}\, \subseteq
\,\ \cdots\,
 \cdots \,\
 \subseteq \, F_{2} \, \subseteq F_{1}\, \subseteq \, F_{0} = (-,X_{\Lambda/\mathfrak{r}})
\end{equation}
where each subquotient $\mathbb S_{i} := F_{i}/F_{i+1}$ is simple.
It follows that $|\mathsf{Supp}\mathbb S_{i}| < \infty$, $0\leq i
\leq n$, since if $X_{1}, X_{2}$ are indecomposable objects in $\X$
with $\mathbb S_{i}(X_{1}) \neq 0 \neq \mathbb S_{i}(X_{2})$,
for some $i$, then $X_{1} \cong X_{2}$. Since clearly
\begin{equation}
\mathsf{Supp}\,(-,X_{\Lambda/\mathfrak{r}}) =
\bigcup^{n}_{i=0}\mathsf{Supp}\,\mathbb S_{i}
\end{equation}
we have $\mathsf{Supp}\,(-,X_{\Lambda/\mathfrak{r}}) < \infty$. Then
(3.3), shows that $|\mathsf{Ind}\X| < \infty$, so $\X$ is of finite
representation type.

(ii)  $\Leftrightarrow$ (v) Clearly (ii) $\Rightarrow$ (v). Assume
that any indecomposable module in $\A$ is finitely generated.
Consider the Ziegler spectra $\mathbf{Zg}(\A)$ and ${\bf
Zg}(\Mod\Lambda)$ of the locally finitely presented categories $\A =
\colim\X$ and $\Mod\Lambda = \colim\smod\Lambda$. Then
$\mathbf{Zg}(\A) = \mathsf{Ind}(\X)$ since any indecomposable module
in $\A$ is finitely generated and any finitely generated module over
an Artin algebra is pure-injective in $\Mod\Lambda$. On the other
hand since $\X$ has left almost split morphisms, there exists a left
almost split morphism $\phi \colon X \lxr Y$ in $\X$, for any point
$X$ in $\mathsf{Ind}(\X)$. If $F$ lies in $\mathsf{Ind}(\X)\setminus
\{X\}$ and $\alpha \colon X\lxr F$ is a map, then $\alpha$ factors
through $\phi$ since otherwise $\alpha$ would be split monic
implying that $X \cong F$ which is not the case. Hence
$\mathsf{Ind}(\X)\setminus \{X\}$ consists of $\phi$-injectives,
i.e. $\mathsf{Ind}(\X)\setminus \{X\}$ is closed and therefore
$\{X\}$ is Ziegler-open. This means that any point $X$ in ${\bf
Zg}(\A)$ is isolated in the Ziegler topology of ${\bf Zg}(\A)$.
We
claim that ${\bf Zg}(\A)$  is a quasi-compact space, i.e. any open
cover of ${\bf Zg}(\A)$ admits a finite subcover.
So assume that $
{\bf Zg}(\A) = \bigcup_{i\in I}{\bf V}_{i}$ where each ${\bf V}_{i}$
is Ziegler-open subset. Then ${\bf V}_{i}$ is of the form ${\bf
V}_{i} = \{G \in {\bf Zg}(\A) \,\ | \,\ \text{the map}\,\
\Hom_{\Lambda}(Y,G) \lxr \Hom_{\Lambda}(X,G) \,\,\, \text{is not
surjective} \big\}$ for some collection of maps $X \lxr Y$ between
modules in $\X$. However each ${\bf V}_{i}$ is also a Ziegler-open
in ${\bf Zg}(\Mod\Lambda)$ since pure-injective modules in $\A$
remain pure-injective in $\Mod\Lambda$ and $\X$ consists of finitely
generated modules. Since the space ${\bf Zg}(\Mod\Lambda)$ is
quasi-compact, see \cite{Krause:memoirs}, we infer that there exists
a finite subcover: $ {\bf Zg}(\A) = \bigcup_{i\in J}{\bf V}_{i}$ for
some finite subset $J$ of $I$. It follows that ${\bf Zg}(\A)$ is
quasi-compact. Then the open cover ${\bf Zg}(\A) = \bigcup_{X\in
\mathsf{Ind}(\X)}\{X\}$ admits a finite subcover ${\bf Zg}(\A) =
\bigcup^{n}_{i = 1}\{X_{i}\}$. This means that $\mathsf{Ind}(\X) =
\{X_{1},X_{2},\cdots,X_{n}\}$, i.e. $\X$ is of finite representation
type.

(i) $\Leftrightarrow$ (vi) Taking into account that  (i) is
equivalent to (ii) and the fact that $\X$ is a Krull-Schmidt
category, we infer that (i) implies (vi). If (vi) holds then, by
\cite[Corollary 2.7]{Krause:memoirs}, $\A$ is pure-semisimple and
then as in the proof of the direction (ii) $\Rightarrow$ (iii) we
deduce that $\X$ is of finite representation type.

Finally the equivalences  (i) $\Leftrightarrow$ (vii)
$\Leftrightarrow$ (viii) are straightforward.
\end{proof}
\end{thm}

\begin{rem} (i) There is a dual version of Theorem $3.1$ concerning
coresolving covariantly finite subcategories of $\smod\Lambda$. We
leave its formulation to the reader.

(ii) Let $\Lambda$ be connected and let $\X$ be a contravariantly
finite resolving subcategory of $\smod\Lambda$. Using Harada-Sai's
Lemma as in \cite{Ringel:quasi} it follows that $\X$ is of finite
representation type if and only if the (relative) Auslander-Reiten
quiver $\Gamma(\X)$, defined by using irreducible maps in $\X$,
contains a component whose indecomposable modules are of bounded
length.
\end{rem}

Subcategories of $\Mod\Lambda$ of the form $\colim\X$, where $\X$ is
a covariantly finite subcategory of $\smod\Lambda$, are important
examples of {\em definable} subcategories of $\Mod\Lambda$, that is,
subcategories closed under filtered colimits, products and pure
submodules.  On the other hand if $\A$ is a coresolving definable
subcategory of $\Mod\Lambda$ and any module in $\A$ is a direct sum
of finitely generated modules, then clearly $\A = \colim\X$ where
$\X$ is a coresolving subcategory of $\smod\Lambda$ and therefore
$\A$ is locally finitely presented. Thus as a consequence of Theorem
$3.1$ we have the following.

\begin{cor} For a coresolving definable subcategory $\A$ of
$\Mod\Lambda$ the following are equivalent.
\begin{enumerate}
\item Any module in  $\A$ is a direct sum of finitely generated
modules.
\item $\A = \colim\Y$, where $\Y \subseteq \smod\Lambda$ is a coresolving subcategory of
finite representation type.
\end{enumerate}
\end{cor}

\begin{cor} Let $(\A,\B)$ be a cotorsion pair in
$\Mod \Lambda$, Then the following are equivalent.
\begin{enumerate}
\item $(\A,\B)$ is of finite type and the subcategory $\A^{\mathsf{fin}} := \A \cap \smod\Lambda$ is
of finite representation type.
\item $\A^{\fin}$ is contravariantly finite and any module in $\A$ is a direct sum of finitely generated
modules.
\end{enumerate}
\begin{proof} If (i) holds, then since $\A^{\fin}$ is contravariantly finite and $(\A,\B)$ is of finite type,
we have $\A = \colim\A^{\mathsf{fin}}$ by \cite[Theorem 5.3]{AST}.
Since $\A^{\fin}$ is resolving, the assertion in (ii) follows from
Theorem $3.1$. If (ii) holds, then $\A = \colim\A^{\mathsf{fin}} =
\Add\A^{\fin}$. This implies easily that $\B = (\A^{\fin})^{\bot}$
so the cotorsion pair $(\A,\B)$ is of finite type. Since $\A^{\fin}$
is contravariantly finite, $\A^{\fin}$ is representation finite by
Theorem $3.1$.
\end{proof}
\end{cor}

Recall that if $\C$ is a class of modules over a ring $R$, then a
module $X$ is called $\C$-{\bf filtered} if there exists an ordinal
$\tau$ and an increasing chain $\{X_{\alpha} \, | \, \alpha <
\tau\}$ of submodules of $X$ such that $X_{0} = X$, $X_{\alpha} =
\bigcup_{\beta <\alpha}X_{\beta}$ for each limit ordinal $\alpha <
\tau$, $X = \bigcup_{\alpha <\tau}X_{\alpha}$, and each subquotient
$X_{\alpha+1}/X_{\alpha}$ is isomorphic to a module in $\C$, for any
$\alpha + 1  < \tau$. Note that if $\X \subseteq \smod R$, then the
category $\mathfrak{F}(\X)$ of $\X$-filtered modules is contained in
$\colim\X$; this inclusion is an equality if $\mathfrak{F}(\X)$ is
closed under filtered colimits.

\begin{cor} For a resolving subcategory  $\X$ of $\smod\Lambda$, the
following are equivalent.
\begin{enumerate}
\item $\X$ is of finite representation type.
\item $\X$ is contravariantly finite and any module in $\colim\X$ is
a  direct sum of finitely generated modules.
\item $\X$ is contravariantly finite and any $\X$-filtered module is
a  direct sum of finitely generated modules.
\end{enumerate}
\begin{proof} (i) $\Leftrightarrow$ (ii) follows from Theorem $3.1$.
By \cite{AT}, the subcategory $\X$ cogenerates a cotorsion pair of
finite type $(\A,\B)$  in $\Mod\Lambda$, where $\B =
\X^{\pmb{\bot}}$ and $\A = {^{\bot}}\B$, and moreover we have
$\A^{\fin} = \X$ and $\A$ consists of all direct summands of
$\X$-filtered modules. Then (i) $\Leftrightarrow$ (iii) follows from
Corollary $3.4$.
\end{proof}
\end{cor}

\begin{exam} Let $T$ be a (possibly infinitely generated) tilting module over
 $\Lambda$, see for instance \cite{AST}.  Then
$(\X,T^{\bot})$ is a cotorsion pair of finite type in $\Mod\Lambda$,
where $\X = {^{\bot}}(T^{\bot})$, see \cite{BS}. Hence $\X^{\fin}$
is of finite representation type if and only if $\X^{\fin}$ is
contravariantly finite and any module in $\X$ is a direct sum of
finitely generated modules. If $T$ is finitely generated, then any
module in $T^{\pmb \bot}$ is a direct sum of indecomposable modules
if and only if $T^{\bot}$ is of finite representation type and
$T^{\pmb \bot} = \colim T^{\bot}$. Similar remarks apply for
finitely generated cotilting modules.
\end{exam}

Applying Theorem $3.1$ to $\X = \smod\Lambda$ we have the following
version of the classical result, due to Auslander and
Ringel-Tachikawa, characterizing Artin algebras of finite
representation type.

\begin{cor} \cite{Auslander:large, Auslander:II, RT} For an Artin algebra $\Lambda$ the
following are equivalent.
\begin{enumerate}
\item $\Lambda$ is of finite representation type.
\item Any $\Lambda$-module is a direct sum of finitely generated
modules.
\item Any $\Lambda$-module is a direct sum of indecomposable
modules.
\item Any indecomposable $\Lambda$-module is finitely generated.
\item The category $\Mod\,(\smod\Lambda)$ is locally finite.
\item The functor $(-,\Lambda/\mathfrak{r}) \in \Mod\,(\smod\Lambda)$ has finite length.
\item $\Mod\Lambda$ is equivalent to the category of projective modules over
an Artin algebra.
\item $\smod\Lambda$ is equivalent to the category of finitely generated projective modules over
an Artin algebra.
\end{enumerate}
\end{cor}

If $\E$ is an exact category, let $\mathsf{K}_{0}(\E)$ be the
Grothendieck group of $\E$ and $\mathsf{K}_{0}(\E,\oplus)$ be the
Grothendieck group of $\E$ endowed with the split exact structure;
e.g. if $\X \subseteq \smod\Lambda$ is resolving, then  $\X$,
$\X^{\bot}$ and $\X \cap \X^{\bot}$ are exact subcategories of
$\smod\Lambda$ since they are closed under extensions. Note that
$\mathsf{K}_{0}(\X\cap \X^{\bot}) = \mathsf{K}_{0}(\X\cap
\X^{\bot},\oplus)$ since $\X\cap \X^{\bot}$ is $\Ext$-orthogonal.
Working as in \cite[Theorem 8.9, Corollary 8.10]{B:cm}, we have the
following.

\begin{prop} Let $\X$ be a contravariantly resolving subcategory of $\smod\Lambda$.
\begin{enumerate}
\item There is a decomposition $\mathsf{K}_{0}(\X) \oplus
\mathsf{K}_{0}(\X^{\bot}) \cong \mathsf{K}_{0}(\X\cap \X^{\bot})
\oplus \mathsf{K}_{0}(\smod\Lambda)$.
\item The following are equivalent.
\begin{enumerate}
\item $\X$ is of finite representation type.
\item The set $\big\{[X_{1}]-[X_{2}]+[X_{3}]\big\} \cup
\big\{[X_{\mathfrak{r}P}]-[P] \big\}$ is a free basis of
$\mathsf{K}_{0}(\X,\oplus)$, where $0 \lxr X_{1} \lxr X_{2} \lxr
X_{3} \lxr 0$ is an Auslander-Reiten sequence in $\X$ and
$X_{\mathfrak{r}P}$ is the minimal right $\X$-approximation of
$\mathfrak{r}P$, for any indecomposable projective module $P$.
\end{enumerate}
\end{enumerate}
\end{prop}

\subsection{Weyl Modules and Quasi-Hereditary Algebras} We consider modules having good filtrations
which arise in the study of highest weight categories. Let $\Lambda$
be a quasi-hereditary algebra. We denote by $\Delta$ the Weyl
modules and by $\nabla$ the induced modules, see \cite{Ringel:good}
for more details. Let $\filt \Add \Delta$, resp. $\filt \Delta$, be
the full subcategory of $\Mod\Lambda$, resp. $\smod\Lambda$,
consisting of the modules with Weyl filtration. Dually let $\filt
\Add \nabla$, resp. $\filt \nabla$, the full subcategory of
$\Mod\Lambda$, resp. $\smod\Lambda$ consisting of the modules with
good filtration. Recall from \cite{Ringel:quasi} that $\Lambda$ is
called $\Delta$-{\em finite}, resp. $\nabla$-{\em finite}, if the
full subcategory of finitely generated Weyl, resp. good, modules is
of finite representation type.  By results of Ringel
\cite{Ringel:good} it follows that $\filt\nabla$ is a
contravariantly finite resolving subcategory of $\smod\Lambda$ and
by results of Krause-Solberg \cite{KS:appl} we have $\filt\Add
\Delta = \colim \filt\Delta$ and $\filt\Add \nabla = \colim
\filt\nabla$. Hence by Theorem $3.1$ and Corollary $3.3$ we have the
following consequence.

\begin{cor} Let $\Lambda$ be a quasi-hereditary algebra. Then
$\Lambda$ is $\Delta$-finite, resp. $\nabla$-finite, if and only if
any Weyl, resp. good, module is a direct sum of finitely generated
Weyl, resp. good, modules.
\end{cor}

\subsection{Submodule Categories} For an abelian category $\A$, we denote by $\proj^{\leq 1}\A$ the full subcategory of $\A$ consisting of all objects  $A$ with $\pd A\leq 1$, and by $\Sub(\A)$ the {\em
subobject category} of $\A$; the objects of $\Sub(\A)$ are
monomorphisms $\mu : A \monic B$ in $\A$ and a morphism in
$\Sub(\A)$ between $\mu : A \monic B$ and $\mu^{\prime} : A^{\prime}
\monic B^{\prime}$ is given by a pair of morphisms $g : A \lxr
A^{\prime}$ and $f : B \lxr B^{\prime}$ such that $\mu \circ f = g
\circ \mu^{\prime}$.

If  $\Lambda$ is an Artin algebra of finite representation type, then we denote by
$\mathsf{A}(\Lambda)$ the {\em Auslander algebra} of $\Lambda$, i.e.
$\mathsf{A}(\Lambda) = \End_{\Lambda}(T)^{\op}$, where $T$ is the
direct sum of the isoclasses of indecomposable finitely generated
modules.

\begin{prop} If $\Lambda$ is an Artin algebra, then there is an equivalence
\begin{equation}
\Sub(\smod\Lambda)\big/\smod\Lambda \ \ \stackrel{\approx}{\lxr} \ \  {\proj}^{\leq 1}\smod \, (\smod\Lambda)
\end{equation}
where we identify $\smod\Lambda$ with the full subcategory $\U$ of $\Sub(\smod\Lambda)$ consisting of the isomorphisms. Moreover
the following statements  are equivalent.
\begin{enumerate}
\item $\Sub(\smod\Lambda)$ is of finite representation type.
\item Any monomorphism in $\Mod\Lambda$ is a direct sum of
monomorphisms in $\smod\Lambda$.
\item Any indecomposable object of $\Sub(\Mod\Lambda)$ lies in
$\Sub(\smod\Lambda)$.
\item  ${\proj}^{\leq 1}\smod \, (\smod\Lambda)$ is of finite representation type.
\item $\Lambda$ and $\proj^{\leq 1}\mathsf{A}(\Lambda)$ are of finite representation type.
\end{enumerate}
\begin{proof} Let $\mathsf{Y} \colon \smod\Lambda \lxr \smod \, (\smod\Lambda)$, $A \mapsto \mathsf{Y}(A) = \Hom_{\Lambda}(-,A)$, be the Yoneda embedding which induces an equivalence between $\smod\Lambda$ and $\proj\smod \, (\smod\Lambda)$. Define a functor 
\[
\mathsf{H} \colon \Sub(\smod\Lambda) \, \lxr \, \smod \, (\smod\Lambda), \ \ \  \text{by} \ \ \  \mathsf{H}(A^{\prime} \monic A) = \Coker\bigl(\mathsf{Y}(A^{\prime}) \xr{} \mathsf{Y}(A)\bigr)
\] Clearly the essential image of $\mathsf{H}$ coincides with the full subcategory $\proj^{\leq 1}\smod \, (\smod\Lambda)$.  It is easy to
see that $\mathsf{H}$ induces a full functor $\mathsf{H} :
\Sub(\smod\Lambda) \lxr \proj^{\leq 1}\smod \, (\smod\Lambda)$ which is
surjective on objects and moreover a map in $\Sub(\smod\Lambda)$ is
killed by $\mathsf{H}$ iff it factorizes through an object of the
form $(X \xr{\cong} X)$, i.e. $\Ker\mathsf{H} = \U$. Since the categories $\smod\Lambda$ and $\U$ are equivalent via the functor $A \longmapsto (A \xr{=} A)$, it follows that $\mathsf{H}$ induces the desired equivalence $(3.6)$.

As in \cite{RS1} we identify the morphism category of
$\Mod\Lambda$ with the module category over the triangular matrix
algebra $T_{2}(\Lambda) := \bigl(\begin{smallmatrix} \Lambda & \Lambda\\
0 & \Lambda
\end{smallmatrix}\bigr)$. Under this identification
$\Sub(\smod\Lambda)$ is a full subcategory of $\smod T_{2}(\Lambda)$
and $\Sub(\Mod\Lambda)$ is a full subcategory of $\Mod
T_{2}(\Lambda)$, and both contain the projectives and are closed
under extensions and submodules. Since filtered colimits are exact,
it follows that $\colim\Sub(\smod\Lambda) \subseteq
\Sub(\Mod\Lambda)$. On the other hand  it is well known that any
monomorphism $f : A \monic B$ in $\Mod\Lambda$ can be written as a
filtered colimit of monomorphisms between finitely generated
$\Lambda$-modules, see \cite{Krause:lfp}. Therefore
\[
\Sub(\Mod\Lambda) = \colim\Sub(\smod\Lambda)
\]
Since
$\Sub(\smod\Lambda)$ is contravariantly finite in $\smod
T_{2}(\Lambda)$, see \cite{RS1}, the first three conditions are
equivalent by Theorem $3.1$.

(i) $\Leftrightarrow$ (iv) If $\Sub(\smod\Lambda)$ is of finite representation type, then the equivalence $(3.6)$ shows that
so is the full subcategory ${\proj}^{\leq 1}\smod \, (\smod\Lambda)$. Conversely if this holds, then the full subcategory $\U \approx \smod\Lambda \approx {\proj}\smod \, (\smod\Lambda) \subseteq {\proj}^{\leq 1}\smod \, (\smod\Lambda)$ is of finite representation type, i.e. $\Lambda$ is of finite representation type, and the equivalence $(3.6)$ shows that  $\Sub(\smod\Lambda)$ is of finite representation type. Finally (iv) $\Leftrightarrow$ (v)  follows from the fact that $\smod \, (\smod\Lambda) \approx \smod\mathsf{A}(\Lambda)$, provided that $\Lambda$ is of finite representation type.
\end{proof}
\end{prop}

\begin{exam} For $n \geq 1$, let $\Lambda_{n} = k[t]/(t^{n})$. By
results of Ringel and Schmidmeier \cite{RS2}, $\Sub(\smod\Lambda_{n})$
is of finite representation type iff $n\leq 5$. It is well known
that the Auslander algebra $\mathsf{A}(\Lambda_{n})$ of
$\Lambda_{n}$ is quasi-hereditary and the Weyl modules coincide with
the modules with projective dimension $\leq 1$. So by the above
results, $\Sub(\smod\Lambda_{n})$ is of finite representation type
iff $\mathsf{A}(\Lambda_{n})$ is $\Delta$-finite iff any
monomorphism in $\Mod\Lambda_{n}$ is a direct sum of monomorphisms
in $\smod\Lambda_{n}$ iff $n \leq 5$.
\end{exam}

\subsection{Torsion Pairs} Let $\Lambda$ be an Artin algebra.  We consider full
subcategories $\F$ of $\smod\Lambda$ which are closed under
extensions and subobjects.  Clearly any such subcategory forms the
torsion-free part $\F$ of a torsion pair $(\T,\F)$ in
$\smod\Lambda$, where $\T = \{A \in \smod\Lambda\, | \,
\Hom_{\Lambda}(A,\F) = 0\}$. For instance, since for an algebra $\Lambda$ of finite representation type we have $\gd\mathsf{A}(\Lambda) \leq 2$, we can take $\F$ to be the full subcategory $\proj^{\leq 1}\mathsf{A}(\Lambda)$ considered above.

\begin{prop} Let $\F$ be a full subcategory of $\smod\Lambda$
which contains $\Lambda$ and is closed under extensions and
submodules. Then the following are equivalent.
\begin{enumerate}
\item $\F$ is of finite representation type.
\item $\F$ is contravariantly finite and any module in $\colim\F$
is a direct sum of finitely generated modules.
\item $\F$ is contravariantly finite and any indecomposable module in $\colim\F$
is finitely generated.
\item $\F = {^{\bot}}T$, where $T$ is a finitely generated
cotilting module with $\id T \leq 1$, and any, resp. indecomposable,
module in ${^{\pmb{\bot}}}T$ is a direct sum of finitely generated
modules, resp. finitely generated.
\end{enumerate}
\begin{proof}  Clearly $\F$ is resolving so the first three conditions are equivalent by
Theorem $3.1$, Corollary $3.5$.

  (i) $\Rightarrow$ (iv) Since $\F$ is contravariantly finite, there exists a
  cotorsion pair $(\F,\R)$ in $\smod\Lambda$, see
\cite{AR:applications}, and clearly $\rd_{\F}\smod\Lambda \leq 1$.
Hence by \cite{AR:applications} there exists a cotilting module $T$
in $\smod\Lambda$ with  $\id T \leq 1$ such that ${^{\bot}}T = \F$.
It is not difficult to see that $\colim\F = {^{\pmb{\bot}}}T$, and
then (iv) follows form (ii) and (iii).

(iv) $\Rightarrow$ (i) The hypothesis implies that we have a
cotorsion pair $(\colim\F,\Y)$ in $\Mod\Lambda$ and
$(\colim\F)^{\fin} = \F = {^{\bot}}T$ is contravariantly finite in
$\smod\Lambda$. Then $\F$ is of finite representation type by
Corollary $3.4$.
\end{proof}
\end{prop}

Subcategories of the form ${^{\bbot}}T$, where $T$ is a finitely generated
cotilting module with $\id_{\Lambda}T \leq 1$, are closed under submodules, extensions,
 filtered colimits and contain $\Lambda$. Thus the
following gives a converse to Proposition $3.12$.

\begin{cor} Let $\A$ be a full subcategory of $\Mod\Lambda$
which contains $\Lambda$ and is closed under extensions, submodules
and filtered colimits. Then the following are equivalent:
\begin{enumerate}
\item $\A^{\fin}$ is of finite representation type.
\item Any module in $\A$ is a direct sum of indecomposable modules.
\item Any indecomposable module in $\A$ is finitely generated.
\end{enumerate}
If $\mathrm{(i)}$ holds, then $\A = {^{\pmb{\bot}}}T$, where $T$ is
a finitely generated $1$-cotilting module and $\A$ is definable.
\begin{proof} By \cite{BuK} it follows that $\A  = \colim\A^{\fin}$
and clearly $\A^{\fin}$ is resolving. Then the assertions  follow
from Theorem $3.1$ and Proposition $3.12$.
\end{proof}
\end{cor}

\section{Artin Algebras and Local Rings of Finite CM-Type}
In this section we apply the results of the previous section to the
study of Gorenstein-projective modules over an Artin algebra. In
particular we give several characterizations of virtually Gorenstein
Artin algebras of finite Cohen-Macaulay type in terms of
decomposition properties of the category of Gorenstein-projective
modules. We also prove analogous  results in the commutative
setting.

\subsection{Gorenstein-Projective Modules and Virtually Gorenstein
Algebras}Let $\A$ be an abelian category.  An acyclic complex of
projective, resp. injective, objects of $\A$
$$X^{\bullet}=\;\;\;\;\cdots\lxr X^{i-1}\lxr X^{i}\lxr X^{i+1}\lxr\cdots$$ is called {\bf totally acyclic},
if $\Hom_{\A}(X^{\bullet},P)$, resp. $\Hom_{\A}(I,X^{\bullet})$, is
acyclic, $\forall P \in \Proj\A$, resp. $\forall I \in \Inj\A$.

The following classes of objects have been introduced by Auslander
and Bridger \cite{ABr} in the context of finitely generated modules
over a Noetherian ring, and by Enochs and Jenda \cite{EJ} in the
context of arbitrary modules over any ring.

\begin{defn} Let $\A$ be an abelian category. An object $A \in \A$
is  called:
\begin{enumerate}
\item {\bf Gorenstein-projective} if it is of
the form $A=\Coker (X^{-1}\to X^0)$ for some totally acyclic complex
$X^{\bullet}$ of projective objects.
\item {\bf Gorenstein-injective} if it is of the from $A=\Ker
(X^{0}\to X^1)$ for some totally acyclic complex $X^{\bullet}$ of
injective objects.
\end{enumerate}
\end{defn}

\begin{rem} Depending on the setting, finitely generated Gorenstein-projective modules
appear in the literature under various names, for instance: totally
reflexive modules \cite{JS}, modules of G-dimension zero \cite{ABr,
Yoshino} or (maximal) Cohen-Macaulay modules \cite{AR:cm, BR}.
\end{rem}

From now on we fix an Artin algebra $\Lambda$. We denote by
$\GProj\La$ the full subcategory of $\Mod\La$ which is formed by all
Gorenstein-projective $\La$-modules, and $\GInj\La$ denotes the full
subcategory which is formed by all Gorenstein-injective
$\La$-modules. In \cite{B:cm}, it is shown that $\GProj\Lambda$ is a
definable subcategory of $\Mod\Lambda$ and it is the largest
resolving subcategory of $\Mod\Lambda$ which admits the full
subcategory $\Proj\Lambda$ of projective modules as an Ext-injective
cogenerator.  Moreover there are cotorsion pairs
\[\big(\GProj\La,(\GProj\La)^\perp\big) \quad\text{and} \quad
\big({^\perp(\GInj\La)},\GInj\La\big)\]
 for $\Mod\La$ satisfying: \,
$\GProj\La\bigcap(\GProj\La)^\perp=\Proj\La \quad\text{and} \quad
{^\perp(\GInj\La)}\bigcap\GInj\La=\Inj\La.$
 Throughout we shall
use the category $\Gproj\La$ of finitely generated
Gorenstein-projective modules and the category $\Ginj\La$ of
finitely Gorenstein-injective modules defined as follows:
\[
\Gproj\La:=\GProj\La\bigcap\smod\La\quad\text{and}\quad
\Ginj\La:=\GInj\La\bigcap\smod\La.
\]

\begin{rem}
Note that the $\Lambda$-dual functor $\mathsf{d} :=
\Hom_{\Lambda}(-,\Lambda) : \smod\Lambda \to \smod\Lambda^{\op}$
induces a duality
\[
\mathsf{d} \, : \, (\Gproj\Lambda)^{\op} \,\
\stackrel{\approx}{\lxr} \,\ \Gproj\Lambda^{\op}
\]
Moreover the adjoint pair $\bigl(\mathsf{N}^{+}:=
-\otimes_{\Lambda}\mathsf{D}\Lambda, \mathsf{N}^{-} :=
\Hom_{\Lambda}(\mathsf{D}\Lambda,-)\bigr) \colon \Mod\Lambda
\leftrightarrows \Mod\Lambda$ induces equivalences
\[
(\mathsf{N}^{+},\mathsf{N}^{-}) \, : \, \GProj\Lambda\,\
\stackrel{\approx}{\lxr} \,\ \GInj\Lambda \,\,\,\,\,\,\,\,\
\text{and} \,\,\,\,\,\,\,\,\ (\mathsf{N}^{+},\mathsf{N}^{-}) \, : \,
\Gproj\Lambda\,\ \stackrel{\approx}{\lxr} \,\ \Ginj\Lambda
\]
The stable categories  $\uGProj\Lambda$ and $\oGInj\Lambda$ are
compactly generated triangulated categories, the stable categories
$\uGproj\Lambda$ and $\oGinj\Lambda$ are triangulated subcategories
of $\uGProj\Lambda$ and $\oGInj\Lambda$ respectively  consisting of
compact objects, and  the adjoint pair of Nakayama functors induce
triangulated equivalences
\[
(\mathsf{N}^{+},\mathsf{N}^{-}) \, : \, \uGProj\Lambda\,\
\stackrel{\approx}{\lxr} \,\ \oGInj\Lambda \,\,\,\,\, \text{and}
\,\,\,\,\,\,\, (\mathsf{N}^{+},\mathsf{N}^{-}) \, : \,
\uGproj\Lambda\,\ \stackrel{\approx}{\lxr} \,\ \oGinj\Lambda.
\]
\end{rem}

It is not difficult to see that $(\Gproj\Lambda)^{\bot} =
{^{\bot}}(\Ginj\Lambda)$ \cite{B:cm}. However there are Artin
algebras for which $(\GProj\La)^{\bot} \neq {^{\bot}}(\GInj\La)$
\cite{BK}. To remedy this pathology, the following class of algebras
was introduced in \cite{BR} as a common generalization of Gorenstein
algebras and algebras of finite representation type.

\begin{defn} An Artin algebra $\La$ is called {\bf virtually Gorenstein}
if:
\[(\GProj\La)^\perp \, = \, {^\perp(\GInj\La)}.\]
\end{defn}

\begin{exam} The following list shows that the class of virtually
Gorenstein algebras is rather large, see \cite{B:cm} for more
details and \cite{BK} for an example of an Artin algebra which is
not virtually Gorenstein.
\begin{enumerate}
\item[{\bf 1.}] Gorenstein algebras are virtually Gorenstein.

Recall that an Artin algebra $\Lambda$ is called {\bf Gorenstein} if
$\id{_{\Lambda}}\Lambda < \infty$ and $\id\Lambda_{\Lambda} <
\infty$; in this case it is known that $\id{_{\Lambda}}\Lambda =
\id\Lambda_{\Lambda}$. Important classes of Gorenstein algebras are
the algebras of finite global dimension and the self-injective
algebras. In the first case we have $\GProj\Lambda = \Proj\Lambda$
and in the second case we have $\GProj\Lambda = \Mod\Lambda$.
\item[{\bf 2.}] Any Artin algebra of finite representation type is virtually
Gorenstein.
\item[{\bf 3.}] Let $\Lambda$ be an Artin algebra and assume that any
  Gorenstein-projective module is projective,
i.e. $\GProj\Lambda = \Proj\Lambda$.  Then $\Lambda$ is virtually
Gorenstein.
\item[{\bf 4.}] Let $\Lambda$ be an Artin algebra which is derived equivalent
or stably equivalent to a virtually Gorenstein algebra, e.g. to an
algebra of finite representation type. Then $\Lambda$ is virtually
Gorenstein.
\end{enumerate}
\end{exam}

We refer to \cite{B:cm} for an extensive discussion of virtually
Gorenstein algebras. In particular we need in the sequel the
following properties enjoyed by the class of virtually Gorenstein
algebras, see \cite{B:cm, BK}:

\begin{rem} Let $\Lambda$ be an Artin algebra.
\begin{enumerate}
\item[{\bf 1}.] $\Lambda$ is virtually Gorenstein if and only if
$\GProj\Lambda = \colim\Gproj\Lambda$. In this case $\GProj\Lambda$
is a locally finitely presented additive category with products and
$(\GProj\Lambda)^{\bot} = \colim(\Gproj\Lambda)^{\bot}$.
\item[{\bf 2.}] $\Lambda$ is virtually Gorenstein if and only if the cotorsion pair $\big(\GProj\Lambda,
(\GProj\Lambda)^{\bot}\big)$ is smashing.
\item[{\bf 3.}]  $\Lambda$ is virtually Gorenstein if and only if  $(\GProj\Lambda)^{\bot}\cap \smod\Lambda$ is
contravariantly, or equivalently covariantly, finite, in
$\smod\Lambda$. In this case $(\GProj\Lambda)^{\bot}\cap
\smod\Lambda = (\Gproj\Lambda)^{\bot}$, see {\bf 4.} below.
\item[{\bf 4.}] $\Lambda$ is virtually Gorenstein if and only if
$(\Gproj\Lambda)^{\pmb{\bot}} \subseteq (\GProj\Lambda)^{\bot}$.
Indeed if the last inclusion holds, then since we always have
$(\GProj\Lambda)^{\bot} \subseteq (\Gproj\Lambda)^{\pmb{\bot}}$, ir
follows that $(\GProj\Lambda)^{\bot} =
(\Gproj\Lambda)^{\pmb{\bot}}$. This implies that
$(\GProj\Lambda)^{\bot}$ is closed under coproducts, i.e.  the
cotorsion pair $\big(\GProj\Lambda, (\GProj\Lambda)^{\bot}\big)$ is
smashing. Therefore by {\bf 2.}, $\Lambda$ is virtually Gorenstein.
If this holds, then by \cite{B:cm} we have $(\GProj\Lambda)^{\bot} =
X^{\pmb{\bot}}_{\Lambda/\mathfrak{r}}$ and the minimal right
$\GProj\Lambda$-approximation $X_{\Lambda/\mathfrak{r}}$ of
$\Lambda/\mathfrak{r}$ is finitely generated. This clearly implies
that $(\Gproj\Lambda)^{\pmb{\bot}} \subseteq
X^{\pmb{\bot}}_{\Lambda/\mathfrak{r}} = (\GProj\Lambda)^{\bot}$.
\item[{\bf 5.}] If $\Lambda$ is virtually Gorenstein, then
$\Gproj\Lambda$ is functorially finite in $\smod\Lambda$; in
particular $\Gproj\Lambda$ has Auslander-Reiten sequences, and
$\Gproj\Lambda = \filt(X_{\Lambda/\mathfrak{r}})$ and $\GProj\Lambda
= \filt(\Add X_{\Lambda/\mathfrak{r}})$.
\item[{\bf 6.}] $\Lambda$ is virtually
Gorenstein if and only if $\Lambda^{\op}$ is virtually Gorenstein.
\item[{\bf 7.}] $\uGProj\Lambda$ admits a single compact generator $\unX_{\Lambda/\mathfrak{r}}$.
Moreover
$\Lambda$ is virtually Gorenstein iff $(\uGProj\Lambda)^{\cpt} =
\uGproj\Lambda$, i.e.  the compact objects of $\uGProj\Lambda$ are
coming from the finitely generated Gorenstein-projective modules. In
this case $\uGproj\Lambda$ has Auslander-Reiten triangles and admits
a Serre functor.
\end{enumerate}
\end{rem}

We have the following characterization of virtually Gorenstein
algebras which will be useful later.

\begin{prop} For an Artin algebra $\Lambda$, the following
statements are equivalent.
\begin{enumerate}
\item $\Lambda$ is virtually Gorenstein.
\item
\begin{enumerate}
\item $\Gproj\Lambda$ is contravariantly finite.
\item Any module in $\GProj\Lambda \cap
\colim(\Gproj\Lambda)^{\bot}$ is a direct sum of finitely generated
modules.
\end{enumerate}
\end{enumerate}
\begin{proof} If $\Lambda$ is virtually Gorenstein, then
$\Gproj\Lambda$ is contravariantly finite by Remark $4.6.5$, and
(ii)(b) holds since $\GProj\Lambda \cap \colim(\Gproj\Lambda)^{\bot}
= \GProj\Lambda \cap (\GProj\Lambda)^{\bot} = \Proj\Lambda$ by
Remark $4.6.1$. Conversely if (ii) holds, then by \cite[Corollary
9.6]{B:cm}, (a) implies  the existence of a cotorsion pair
$\bigl(\colim\Gproj\Lambda, \colim(\Gproj\Lambda)^{\bot}\bigr)$ in
$\Mod\Lambda$ and $\colim\Gproj\Lambda
 \subseteq \GProj\Lambda$ since the latter is closed under filtered
 colimits. Let $A$ be in $\GProj\Lambda \cap
\colim(\Gproj\Lambda)^{\bot}$. Then by (b), $A = \oplus_{i\in
I}A_{i}$ where each $A_{i}$ is finitely generated and clearly lies
in $\GProj\Lambda \cap \colim(\Gproj\Lambda)^{\bot} \cap
\smod\Lambda = \Gproj\Lambda \cap (\Gproj\Lambda)^{\bot}$ since both
$\GProj\Lambda$ and $\colim(\Gproj\Lambda)^{\bot}$ are closed under
direct summands. However it is easy to see that $\Gproj\Lambda \cap
(\Gproj\Lambda)^{\bot} = \proj\Lambda$ and then this implies that
$\GProj\Lambda \cap \colim(\Gproj\Lambda)^{\bot} = \Add\proj\Lambda
= \Proj\Lambda$. Now let $A \in \GProj\Lambda$ and let $0 \lxr A
\lxr Y^{A} \lxr X^{A} \lxr 0$ be an exact sequence, where $Y^{A}$ is
a left $\colim(\Gproj\Lambda)^{\bot}$-approximation of $A$ and
$X^{A} \in \colim\Gproj\Lambda$. Clearly $Y^{A}$ lies in
$\GProj\Lambda \cap \colim(\Gproj\Lambda)^{\bot} = \Proj\Lambda$ and
therefore $A$ lies in $\colim\Gproj\Lambda$, since the latter is
resolving. We infer that $\GProj\Lambda = \colim\Gproj\Lambda$ and
therefore $\Lambda$ is virtually Gorenstein by Remark $4.6.1$.
\end{proof}
\end{prop}

\begin{cor} If $\Gproj\Lambda = \proj\Lambda$, then the following
are equivalent.
\begin{enumerate}
\item $\Lambda$ is virtually Gorenstein.
\item Any module in $\GProj\Lambda$ is a direct sum of finitely
generated modules.
\end{enumerate}
\end{cor}

\subsection{Artin algebras of finite CM-type} Our aim in this section is to
characterize Artin algebras enjoying the property that any
Gorenstein-projective module is a direct sum of finitely generated
modules.  This class of algebras, clearly contained in the class of
virtually Gorenstein algebras, is intimately related to the class of
Artin algebras of finite Cohen-Macaulay type, in the sense of the
following definition.

\begin{defn} An Artin algebra $\Lambda$ is said to be of {\bf finite Cohen-Macaulay
type}, {\bf finite CM-type} for short, if the full subcategory
$\Gproj\Lambda$ of finitely generated Gorenstein-projective
$\Lambda$-modules is of finite representation type.
\end{defn}

Clearly an Artin algebra $\Lambda$ is virtually Gorenstein of
finite CM-type if $\Lambda$ is of finite representation type or if
$\GProj\Lambda = \Proj\Lambda$, e.g. if $\gd\Lambda < \infty$. On
the other hand if $\Lambda$ is of infinite CM-type, then obviously
$\gd\Lambda = \infty$. Note that a self-injective algebra is of finite
CM-type iff it is of finite representation type. Finally by Remark $4.3$ it follows that $\Lambda$ is of finite
CM-type if and only if so is $\Lambda^{\op}$.

The following result shows that the class of Artin algebras enjoying
the property that any Gorenstein-projective module is a direct sum
of finitely generated modules coincides with the class of virtually
Gorenstein algebras of finite CM-type, and condition (vii)
characterizes them in terms of an interesting factorization
property.
 To avoid trivialities, we assume that
$\GProj\Lambda \neq \Proj\Lambda$.

\begin{thm} For an Artin algebra $\Lambda$ the following are
equivalent:
\begin{enumerate}
\item Any Gorenstein-projective $\Lambda$-module is a direct sum of finitely
generated modules.
\item $\Lambda$ is virtually Gorenstein and any Gorenstein-projective
$\Lambda$-module is a direct sum of indecomposable modules.
\item $\Lambda$ is virtually Gorenstein and any indecomposable Gorenstein-projective module is finitely generated.
\item $\Lambda$ is virtually Gorenstein of finite CM-type.
\item $\Lambda$ is virtually Gorenstein and the triangulated category
$\uGProj\Lambda$ is phantomless.
\item $\Lambda$ is virtually Gorenstein and the functor category $\Mod\uGproj\Lambda$ is Frobenius.
\item A map  $f \colon A \lxr B$ between Gorenstein-projective
modules factors through a projective module if and only if for any
finitely generated Gorenstein-projective module $X$ and any map
$\alpha : X \lxr A$, the composition $\alpha \circ f \colon X \lxr
B$ factors through a projective module.
\item The above conditions with $\Lambda$ replaced by
$\Lambda^{\op}$.
\end{enumerate}
\begin{proof} The proof that the first four conditions are equivalent is a direct consequence of Theorem $3.1$ and the
fact that $\Lambda$ is virtually Gorenstein if and only if
$\GProj\Lambda = \colim\Gproj\Lambda$, see $4.6.1$. The equivalence
(iv) $\Leftrightarrow$ (viii) follows from the fact that $\Lambda$
is virtually Gorenstein if and only if so is $\Lambda^{\op}$, see
$4.6.4$, and from the existence of an  equivalence
$(\Gproj\Lambda)^{\op} \approx \Gproj\Lambda^{\op}$, see Remark
$4.3$. If (i) holds, then $\Lambda$ is virtually Gorenstein and
therefore $(\uGProj\Lambda)^{\cpt} = \uGproj\Lambda$. It follows
from $4.6.7$ that any object in $\uGProj\Lambda$ is a direct sum of compact
objects and therefore by  \cite[Theorem 9.3]{B:3cats},
$\uGProj\Lambda$ is phantomless, so (v) holds. Using  $4.6.6$, the
proof that (v) $\Leftrightarrow$ (vi) follows from \cite[Proposition
9.2, Theorem 9.3]{B:3cats}. 

We now show that (vi) implies (i). Since
$\uGProj\Lambda$ is phantomless, by \cite[Theorem 9.3]{B:3cats}, for
any Gorenstein-projective module $A$ there is a decomposition $\unA
= \oplus_{i\in I}\unX_{i}$, where the $X_{i}$ are finitely
generated. This implies that there are projective modules $P$ and
$Q$ such that $A \oplus P \cong \oplus_{i\in I}X_{i} \oplus Q$.
Since $\Lambda$ is virtually Gorenstein, by $4.6.1$, the category
$\GProj\Lambda$ is locally finitely presented and $\GProj\Lambda =
\colim\Gproj\Lambda$. Consider the representation functor
$\mathsf{H} : \GProj\Lambda \lxr \Mod\Gproj\Lambda$. Then
$\mathsf{H}(A) \oplus \mathsf{H}(P) \cong \mathsf{H}(\oplus_{i\in
I}X_{i}) \oplus \mathsf{H}(Q) \cong \oplus_{i\in I}\mathsf{H}(X_{i})
\oplus \mathsf{H}(Q)$. Since the objects $\mathsf{H}(X_{i})$ and
$\mathsf{H}(Q)$ are projective, we infer that  $\mathsf{H}(A)$ is a
projective object of $\Mod\Gproj\Lambda$. Since any flat functor in
$\Mod\Gproj\Lambda$ is of the form $\mathsf{H}(A)$, where $A\in
\GProj\Lambda$, we infer that any flat functor is projective, i.e.
the representation category $\Mod\Gproj\Lambda$ is perfect. It is
well-known
 that in a perfect functor category any projective functor is a
 direct sum of finitely generated projective functors, see
 \cite{JL}. It follows that $\mathsf{H}(A) = \oplus_{j\in
 J}\mathsf{H}(X^{\prime}_{j})$, where the $X^{\prime}_{j}$ lie in $\Gproj\Lambda$. Since $\mathsf{H}$ is fully faithful
 we conclude that $A$ is a direct sum of finitely generated modules
 and (i) holds. 

Finally we show that (v) $\Leftrightarrow$ (vii). Assuming (vii), let
$\un f \colon \unA \lxr \unB$ be a map in $\uGProj\Lambda$ such that
the induced map $\uHom_{\Lambda}(X,A) \lxr \uHom_{\Lambda}(X,B)$ is
zero for any finitely generated Gorenstein-projective module $X$.
This means that for any map $X \lxr A$, the composition $X \lxr A
\lxr B$ factorizes through a projective module. Then condition (vii)
implies $f$ factorizes through a projective module, so $\unf = 0$ in
$\uGProj\Lambda$. Applying this to $A = B$ and  $f = 1_{A}$, we have
that $\uHom(X,A) = 0$, $\forall X\in \Gproj\Lambda$, implies that
$\unA = 0$, so $\uGproj\Lambda$ (compactly) generates
$\uGProj\Lambda$. Then by \cite[Theorem 9.3]{B:3cats} we infer that
$(\uGProj\Lambda)^{\cpt} = \uGproj\Lambda$, so by $4.6.6$, $\Lambda$
is virtually Gorenstein, and any object of $\uGProj\Lambda$ is a
direct sum of objects from $\uGproj\Lambda$, i.e. $\uGProj\Lambda$
is phantomless and (v) holds.  Conversely if $\uGProj\Lambda$ is
phantomless, then (vii) holds since $\uGproj\Lambda \subseteq
(\uGProj\Lambda)^{\cpt}$.
\end{proof}
\end{thm}

Note that by the above proof, if $\Lambda$ is virtually Gorenstein
of finite CM-type and $T$ is a representation generator of
$\Gproj\Lambda$, then  there are equivalences $\GProj\Lambda \approx
\Proj\Gamma$ and $\Gproj\Lambda \approx \proj\Gamma$, where $\Gamma
= \End_{\Lambda}(T)^{\op}$.    Moreover let $\Delta =
\underline{\mathsf{End}}_{\Lambda}(T)^{\op}$ be the stable
endomorphism algebra of $T$. Since $\uGproj\Lambda =
\add\underline{T}$, there are equivalences $\uGProj\Lambda \approx
\Proj\Delta$ and $\uGproj\Lambda \approx \proj\Delta$. Since
$\uGProj\Lambda$ is triangulated, by \cite{B:3cats, Buchweitz},
$\Delta$ is a self-injective Artin algebra, and, using Theorem $3.1$ and Proposition $3.8$, we have the
following consequence.

\begin{cor} For an Artin algebra $\Lambda$, the following statements
are equivalent.
\begin{enumerate}
\item $\Lambda$ is virtually  Gorenstein of finite CM-type.
\item $\GProj\Lambda$ is equivalent to the category of projective
modules over an Artin algebra.
\item $\uGProj\Lambda$ is equivalent to the category of projective
modules over a self-injective  Artin algebra.
\item $\Lambda$ is virtually Gorenstein and  $\Gproj\Lambda$ is equivalent to the category of finitely generated projective
modules over an Artin algebra.
\item $\Lambda$ is virtually Gorenstein and  $\uGproj\Lambda$ is equivalent to the category of finitely generated projective
modules over a self-injective  Artin algebra.
\item $\Lambda$ is virtually Gorenstein and the set $\big\{[X_{1}]-[X_{2}]+[X_{3}]\big\} \cup
\big\{[X_{\mathfrak{r}P}]-[P] \big\}$ is a free basis of
$\mathsf{K}_{0}(\Gproj\Lambda,\oplus)$, where $0 \lxr X_{1} \lxr
X_{2} \lxr X_{3} \lxr 0$ is an Auslander-Reiten sequence in
$\Gproj\Lambda$ and $X_{\mathfrak{r}P}$ is the minimal right
$\Gproj\Lambda$-approximation of $\mathfrak{r}P$, $\forall P \in
\mathsf{Ind}(\proj\Lambda)$.
\item $\Lambda$ is virtually Gorenstein and the set $\big\{[\unX_{1}]-[\unX_{2}]+[\unX_{3}]\big\}$ is a free basis of
$\mathsf{K}_{0}(\uGproj\Lambda,\oplus)$, where $ \unX_{1} \lxr
\unX_{2} \lxr \unX_{3} \lxr \Omega^{-1}\unX_{1}$ is an
Auslander-Reiten triangle in $\uGproj\Lambda$.
\end{enumerate}
\end{cor}

\begin{rem}  By Remark $4.3$, the equivalent conditions of
Theorem $4.10$ and Corollary $4.11$ which characterize virtually
Gorenstein algebras of finite CM-type, are also equivalent to the
conditions resulting by replacing everywhere Gorenstein-projectives
with Gorenstein-injectives and projectives with injectives.
 \end{rem}

Since Gorenstein algebras are virtually Gorenstein and, by
\cite[Theorem 11.4]{B:cm}, virtually Gorenstein algebras with
one-sided finite self-injective dimension are Gorenstein, we have as
a corollary of Theorem $4.10$, the following result which,
extending the main result of Chen \cite[Main Theorem]{Chen}, gives a
characterization of Gorenstein algebras of finite CM-type.

\begin{cor} For an Artin algebra  $\Lambda$ the following are equivalent.
\begin{enumerate}
\item $\Lambda$ is a Gorenstein algebra of finite CM-type.
\item $\id\Lambda_{\Lambda} < \infty$ and any Gorenstein-projective $\Lambda$-module is a direct sum
of finitely generated modules.
\item $\Lambda$ is virtually Gorenstein with $\id \Lambda_{\Lambda} < \infty$ and any indecomposable
Gorenstein-projective is finitely generated.
\item The above conditions with $\Lambda$ replaced by
$\Lambda^{\op}$.
\end{enumerate}
\end{cor}

By \cite{B:cm} if $\Lambda$ and $\Gamma$ are Artin algebras which
are derived equivalent or stably equivalent of Morita type, then
$\Lambda$ is virtually Gorenstein iff so is $\Gamma$; in this case
the stable triangulated categories $\uGproj\Lambda$ and
$\uGproj\Gamma$ are (triangle) equivalent. Since Artin algebras of
finite representation type are virtually Gorenstein, we have the
following consequence. Note that algebras derived equivalent to
algebras of finite representation type may be of tame infinite, even
wild, representation type.

\begin{cor} Let $\Lambda$ and $\Gamma$ be Artin algebras which are derived equivalent or stably equivalent of Morita type.
Then $\Lambda$ is virtually Gorenstein of finite CM-type iff
$\Gamma$ is virtually Gorenstein of finite CM-type. In particular if
$\Lambda$ is an Artin algebra which is derived equivalent or stably
equivalent of Morita type to an Artin algebra of finite
representation type, then $\Lambda$ is virtually Gorenstein of
finite CM-type.
\end{cor}

\begin{cor} Let $\Lambda$ be an Artin algebra. Then the following
are equivalent.
\begin{enumerate}
\item $(\GProj\Lambda)^{\bot} \cap \smod\Lambda$ is of finite
representation type.
\item ${^{\bot}}(\GInj\Lambda) \cap \smod\Lambda$ is of finite
representation type.
\item $\Lambda$ is virtually Gorenstein and any module in
$(\GProj\Lambda)^{\bot}$ is a direct sum of finitely generated
modules.
\item $\Lambda$ is virtually Gorenstein and any module in
${^{\bot}}(\GInj\Lambda)$ is a direct sum of finitely generated
modules.
\end{enumerate}
\begin{proof} By \cite{B:cm} we have $\X := (\GProj\Lambda)^{\bot} \cap
\smod\Lambda = \smod\Lambda \cap {^{\bot}}(\GInj\Lambda)$, so (i)
$\Leftrightarrow$ (ii). Clearly (iii) $\Leftrightarrow$ (iv). If (i)
holds, then $\X$ is contravariantly finite so, by $4.6.3$, $\Lambda$
is virtually Gorenstein. Since $\X$ is resolving, the assertion in
(iii) follows from $4.6.1$, $4.6.3$, and Theorem $3.1$.  The
implication (iii) $\Rightarrow$ (i) is similar.
\end{proof}
\end{cor}

\begin{rem} Let $\Lambda$ be an Artin algebra and assume that
$\id\Lambda_{\Lambda} < \infty$.
\begin{enumerate}
\item Any module in ${^{\bot}}(\GInj\Lambda)$ is a direct sum of
finitely generated modules if and only if the full subcategory
$\inj^{<\infty}\!\Lambda$ of $\smod\Lambda$, consisting of all
modules with finite injective dimension, is of finite representation
type. This follows from Theorem $3.1$ and the fact that
$\id\Lambda_{\Lambda} < \infty$ implies that
$\inj^{<\infty}\!\Lambda$ is resolving and ${^{\bot}}(\GInj\Lambda)
= \colim\inj^{<\infty}\!\Lambda$, see \cite[Theorem 11.3]{B:cm}.
\item If $\Lambda$ is virtually Gorenstein, then:
any module in $(\GProj\Lambda)^{\bot}$ is a direct sum of finitely
generated modules if and only if the full subcategory
$\proj^{<\infty}\!\Lambda$ of $\smod\Lambda$, consisting of all
modules with finite projective dimension, is of finite
representation type. This follows from Theorem $3.1$ and the fact
that virtual Gorensteinness and $\id\Lambda_{\Lambda} < \infty$
implies that $\Lambda$ is Gorenstein and  $(\GProj\Lambda)^{\bot} =
\colim\proj^{<\infty}\!\Lambda$, see \cite[Theorem 11.4]{B:cm}.
\end{enumerate}
\end{rem}

\begin{exam} Let $\Lambda$ be a self-injective algebra. Then the
following are equivalent.
\begin{enumerate}
\item $\Sub(\smod\Lambda)$ is of finite representation type.
\item $T_{2}(\Lambda) := \bigl(\begin{smallmatrix} \Lambda & \Lambda\\
0 & \Lambda
\end{smallmatrix}\bigr)$ is of finite CM-type.
\end{enumerate}
It is well-known that $T_{2}(\Lambda)$ is Gorenstein with $\id
T_{2}(\Lambda) \leq 1$. This easily implies that $\Gproj
T_{2}(\Lambda) = \Sub(\proj T_{2}(\Lambda))$. On the other hand it
is not difficult to see that $\Sub(\smod\Lambda) = \Sub(\proj
T_{2}(\Lambda))$, hence $\Sub(\smod\Lambda) = \Gproj
T_{2}(\Lambda)$. It follows that the submodule category
$\Sub(\smod\Lambda)$ of $\smod\Lambda$ is of finite representation
type if and only if the Gorenstein algebra $T_{2}(\Lambda)$ is of
finite CM-type. E.g. if $\Lambda_{n} = k[t]/(t^{n})$, then by
Example $3.11$,  $T_{2}(\Lambda_{n})$ is of finite CM-type iff $n
\leq 5$.
\end{exam}

We close this subsection with a characterization of Artin algebras
of finite CM-type. Recall that if $\mathcal S$ is a class of objects
in a triangulated category $\T$, then the localizing subcategory
$\loc(\mathcal S)$ of $\T$ generated by $\mathcal S$, is the
smallest full triangulated subcategory of $\T$ which contains
$\mathcal S$ and is closed under arbitrary direct sums.

\begin{thm} For an Artin algebra $\Lambda$ the following statements
are equivalent.
\begin{enumerate}
\item $\Lambda$ is of finite CM-type.
\item $\Gproj\Lambda$ is contravariantly finite and any
$(\Gproj\Lambda)$-filtered module is a direct sum of finitely
generated modules.
\item $\Gproj\Lambda$ is contravariantly finite and any module in $\colim\Gproj\Lambda$ is a direct sum of
finitely generated modules.
\item $\Gproj\Lambda$ is contravariantly finite and the localizing subcategory of $\uGProj\Lambda$ generated by
$\uGproj\Lambda$ is phantomless.
\item $\Gproj\Lambda$ is contravariantly finite and the module category $\Mod\uGproj\Lambda$ is Frobenius.
\end{enumerate}
\begin{proof} The equivalence (i) $\Leftrightarrow$ (ii) follows
from Corollary $3.5$. By \cite[Theorem 9.4]{B:cm} the localizing
subcategory $\mathcal L = \loc(\uGproj\Lambda)$ of $\uGProj\Lambda$
generated by $\uGproj\Lambda$ is compactly generated,  coincides
with the stable category $(\colim\Gproj\Lambda)/\Proj\Lambda$,  and
we have $\{(\colim\Gproj\Lambda)/\Proj\Lambda\}^{\cpt} =
\uGproj\Lambda$. If (i) holds and $\Gproj\Lambda = \add T$ for some
$T \in \Gproj\Lambda$, then $\Gproj\Lambda$ is contravariantly
finite and the representation category $\Mod\uGproj\Lambda$ of
$\mathcal L$ is equivalent to the module category
$\Mod\uEnd_{\Lambda}(T)^{\op}$, hence it is locally finite. Since $\uEnd_{\Lambda}(T)^{\op}$ is self-injective, by
\cite[Theorem 9.3]{B:3cats} $\mathcal L$ is phantomless. If (iv) holds, then any
object of $\mathcal L$ is a direct sum of objects from
$\uGproj\Lambda$. This implies that any module in
$\colim\Gproj\Lambda$ is a direct sum of finitely generated modules,
so (iii) holds. If (iii) holds, then $\Lambda$ is of finite CM-type by Theorem $3.1$. Finally the equivalence (iv) $\Leftrightarrow$ (v) follows from \cite[Theorem 9.3]{B:3cats}.
\end{proof}
\end{thm}

\subsection{Projective Cotorsion Pairs}  Recall that for any Artin algebra $\Lambda$  we have a projective cotorsion
 pair $\big(\GProj\Lambda,(\GProj\Lambda)^{\bot}\big)$ with heart
 $\GProj\Lambda\bigcap (\GProj\Lambda)^{\bot} = \Proj\Lambda$, and the
 stable  category $\uGProj\Lambda$ is a compactly
 generated triangulated category. If $(\X,\Y)$ is an arbitrary projective
 cotorsion pair in $\Mod\Lambda$, then the stable
 category $\un\X := \X/\Proj\Lambda$ is again triangulated, it
 admits all small  coproducts, and $\X \subseteq \GProj\Lambda$.
 So $\big(\GProj\Lambda,(\GProj\Lambda)^{\bot}\big)$ is the
 largest projective cotorsion pair in $\Mod\Lambda$.

 The following result generalizes Theorem $4.10$ to arbitrary
 projective cotorsion pairs. Note that we may recover Theorem $4.10$ by choosing below $(\X,\Y)$ to be the cotorsion
pair $\big(\GProj\Lambda,(\GProj\Lambda)^{\bot}\big)$ and using
$4.6.2$. The dual version concerning injective cotorsion pairs is
left to the reader.

\begin{thm} If $(\X,\Y)$ is a projective cotorsion pair in
 $\Mod\Lambda$ then the following are equivalent.
 \begin{enumerate}
 \item Any module in $\X$ is a direct sum of finitely generated
 modules.
 \item The cotorsion pair $(\X,\Y)$ is smashing and  $\X^{\fin}$ is of
 finite representation type.
 \item The cotorsion pair $(\X,\Y)$ is smashing and any
  module in $\X$ is a direct sum of indecomposable modules.
 \item The cotorsion pair $(\X,\Y)$ is smashing and any
 indecomposable module in $\X$ is finitely generated.
 \item The cotorsion pair $(\X,\Y)$ is smashing and the triangulated category $\un\X$ is phantomless.
 \item The cotorsion pair $(\X,\Y)$ is smashing and the functor
 category $\Mod\un\X^{\fin}$ is Frobenius.
 \item $\X$ is equivalent to the category of projective
modules over an Artin algebra.
\item $\un\X$ is equivalent to the category of projective
modules over a self-injective  Artin algebra.
\item The cotorsion pair $(\X,\Y)$ is smashing and
$\X^{\fin}$ is equivalent to the category of finitely generated
projective modules over an Artin algebra.
\item The cotorsion pair $(\X,\Y)$ is smashing and  $\un\X^{\fin}$ is equivalent
 to the category of finitely generated projective
modules over a self-injective  Artin algebra.
 \item A map  $f \colon A \lxr B$ between
modules in $\X$ factors through a projective module if and only if
for any module $X$ in $\X^{\fin}$ and any map $\alpha : X \lxr A$,
the composition $\alpha \circ f \colon X \lxr B$ factors through a
projective module.
 \end{enumerate}
\begin{proof} (i) $\Rightarrow$ (ii) Let any module in $\X$
be a direct sum of finitely generated modules, i.e. $\X = \Add
\X^{\,\mathsf{fin}}$. Let $\{Y_{i}\}_{i\in I}$ be a small set of
modules in $\Y$. Then $\oplus_{i\in I}Y_{i}$ lies in $\Y$ if
$\Ext^{n}_{\Lambda}(X,\oplus_{i\in I}Y_{i}) = 0$, $\forall n \geq
1$, $\forall X \in \X$.  However $X$ is a direct summand of a
coproduct $\oplus_{j\in J}X_{j}$, where $X_{j} \in \X^{\fin}$, and
$\Ext^{n}_{\Lambda}(X_{j},\oplus_{i\in I}Y_{i}) = \oplus_{i\in
I}\Ext^{n}_{\Lambda}(X_{j},Y_{i}) = 0$, $\forall n \geq 1$, $\forall
j \in J$, since $X_{j}$ is finitely generated. Hence $\oplus_{i\in
I}Y_{i}$ lies in $\Y$ and therefore $\Y$ is closed under all small
coproducts, i.e. the cotorsion pair $(\X,\Y)$ is smashing. Since
$\un\X = \Add \un\X^{\fin}$, it follows by \cite[Theorem
9.3]{B:3cats} that the stable triangulated category $\un\X$ is
compactly generated. Then by \cite[Theorem 9.10]{B:cm} we have $\X =
\colim\X^{\fin}$ and by Theorem $3.1$, $\X^{\fin}$ is of finite
representation type.

(ii) $\Rightarrow$ (i) Since the cotorsion pair $(\X,\Y)$ is
smashing, as in \cite[Theorem 6.6]{B:cm} we infer that the stable
triangulated category $\un\X$ is compactly generated. Then by
 \cite[Theorem 9.10]{B:cm} we have $\X = \colim\X^{\fin}$. By
Theorem $3.1$ then any module in $\X$ is a direct sum of finitely
generated modules.

The rest of proof is similar to the proof of Theorem $4.10$ and is
left to the reader.
\end{proof}
\end{thm}

\subsection{Local rings of Finite CM-type}
Our aim in this subsection is to give characterizations of
commutative Noetherian complete local rings $(R,\mathfrak{m},k)$ of
finite CM-type in terms of decomposition properties of the category
of Gorenstein-projective modules. Note that if $\GProj R = \Proj R$,
then $\Gproj R = \proj R$ and clearly $R$ is of finite CM-type and
any Gorenstein-projective  $R$-module is a direct sum of finitely
generated modules. So it suffices to treat the case in which $R$
admits a non-free finitely generated Gorenstein-projective module.

First recall that a module $A$ over any ring $S$ is called {\em
Gorenstein-flat} if there exists an acyclic complex $F^{\bullet}$ of
flat $S$-modules such that $A \cong \Ker(F^{0} \lxr F^{1})$ and the
complex $F^{\bullet}\otimes_{S}E$ is exact for any injective
$S^{\op}$-module $E$.

\begin{thm} Let $(R,\mathfrak{m},k)$ be a commutative Noetherian complete local ring. If $\Gproj R \neq \proj
R$, then the following statements are equivalent.
\begin{enumerate}
\item $R$ is of finite CM-type.
\item $R$ is Gorenstein and any Gorenstein-projective $R$-module is a direct sum of finitely
generated modules.
\item $\Gproj R$ is contravariantly finite in $\smod R$ and any
$(\Gproj R)$-filtered module is a direct sum of finitely generated
modules.
\item Any Gorenstein-flat module is a filtered colimit of finitely generated Gorenstein-projective modules and any Gorenstein-projective
module is a direct sum of finitely generated modules.
\end{enumerate}
If $R$ is of finite CM-type, then any indecomposable
Gorenstein-projective module is finitely generated.
\begin{proof} (i) $\Rightarrow$ (ii) By the results of \cite{CPST}, $R$ is Gorenstein. Let $T$ be a finitely generated
Gorenstein-projective module such that $\Gproj R = \add T$.  Then
$\uGproj R = \add \underline{T}$ and let $\Lambda =
\uEnd_{R}(T)^{\op}$ be the stable endomorphism ring of $T$. Consider
the triangulated category $\uGProj R$ and the representation functor
$\mathsf{H} : \uGProj R \lxr \Mod\uGproj R = \Mod \Lambda$, given by
$\mathsf{H}(\unA) = \uHom_{R}(T,A)$. Since $R$ is of finite CM-type,
by \cite[Lemma 4.1]{CPST}, the $R$-module $\Ext^{1}_{R}(M,N)$ is of
finite length as an $R$-module, $\forall M,N \in \Gproj R$. In
particular $\Lambda = \uHom_{R}(T,T) \cong
\Ext^{1}_{R}(\Omega^{-1}T,T)$ has finite length as an $R$-module. It
follows that the $R$-algebra $\Lambda$ is Noetherian and therefore,
by \cite{Buchweitz}, $\Lambda$ is self-injective Artinian, so the
Frobenius category $\Mod\Lambda$ is a locally finite. By results of
Krause \cite{Krause:stable}, the stable category $\oGInj R$ is
compactly generated and its compact objects are given, up to direct
factors, by the Verdier quotient ${\bf D}^{b}(\smod
R)/\mathsf{perf}R$, where $\mathsf{perf}R$ is the full subcategory
of ${\bf D}^{b}(\smod R)$ formed by the perfect complexes. However,
since $R$ is Gorenstein, it is well-known that ${\bf D}^{b}(\smod
R)/\mathsf{perf}R$ is triangle equivalent to  $\uGproj R$, see
\cite{Buchweitz}. Also, by \cite[Theorem VI.3.2]{BR}, Gorensteinness
implies the existence of a triangle equivalence $\uGProj R \approx
\oGInj R$. We infer that $\uGProj R$ is compactly generated and
$(\uGProj R)^{\cpt}$ is triangle equivalent to $\uGproj R$. Since
$\Mod\uGproj R$ is locally finite, by \cite[Theorem 9.3]{B:3cats},
any object $\unA \in \uGProj R$ admits a direct sum decomposition
$\unA = \oplus_{i\in I}\unX_{i}$, where $\{X_{i}\}_{i\in I}$ is a
small subset of $\Gproj R$. Then in $\Mod R$ we have an isomorphism
$A \oplus P \cong \oplus_{i\in I}X_{i} \oplus Q$, where $P, Q \in
\Proj R$, and therefore  $A$ is pure-projective as a direct summand
of a direct sum of finitely presented modules.  Since $R$ is
complete, by a classical result of Warfield \cite[Corollary
4]{Warfield} it follows that $A$ is a direct sum of finitely
generated modules. We conclude that any Gorenstein-projective module
is a direct sum of indecomposable finitely generated
(Gorenstein-projective) modules.

(ii) $\Rightarrow$ (i) We have an equality $\GProj R = \Add \Gproj
R$, and therefore an equality $\uGProj R = \Add \uGproj R$. This
implies that the triangulated category $\uGProj R$ is phantomless
 and compactly generated with $(\uGProj R)^{\cpt}$ $=
\uGproj R$. By \cite[Theorem 9.3]{B:3cats}, the module category
$\Mod\uGproj R$ is Frobenius and then by \cite[Lemma 10.1]{B:art},
the category $\uGproj R$ has right Auslander-triangles. Since the
functor $\Hom_{R}(-,R)$ gives a self-duality of $\uGproj R$, it
follows that $\uGproj R$ has Auslander-Reiten triangles and
therefore by \cite[Theorems 10.2, 10.3]{B:art}, the Frobenius
category $\Mod\uGproj R$ is locally finite and
$\lvert\Supp\uHom_{R}(-,X)\rvert <\infty$, for any indecomposable $X
\in \Gproj R$. In particular
$\lvert\Supp\uHom_{R}(-,X_{R/\mathfrak{m}})\rvert <\infty$, where
$X_{R/\mathfrak{m}}$ is the (minimal) $\Gproj R$-approximation of
$R/\mathfrak{m}$, which exists since $R$ is Gorenstein. Let $X$ be a
non-zero indecomposable module in $\uGproj R$ and assume that $X
\notin \Supp\uHom_{R}(-,X_{R/\mathfrak{m}})$, i.e.
$\uHom_{R}(X,X_{R/\mathfrak{m}}) = 0$. Since $R$ is Gorenstein,
there exists an exact sequence $0 \lxr Y_{R/\mathfrak{m}} \lxr
X_{R/\mathfrak{m}} \lxr R/\mathfrak{m} \lxr 0$, where
$Y_{R/\mathfrak{m}}$ has finite projective dimension. Applying
$\Hom_{R}(X,-)$  we see easily that $\uHom_{R}(X,R/\mathfrak{m}) =
0$, i.e. any map $X \lxr R/\mathfrak{m}$ factors through a finite
free $R$-module. This easily implies that $X$ admits $R$ as a direct
summand and this is impossible since $X$ is non-free indecomposable.
We infer that $\mathsf{Ind}\uGproj R =
\Supp\uHom_{R}(-,X_{R/\mathfrak{m}})$. Since $\Mod\uGproj R$ is
locally finite, the functor $\uHom_{R}(-,X_{R/\mathfrak{m}})$ has
finite length. Then as in the proof of Theorem $3.1$ we infer that
$\Supp\uHom_{R}(-,X_{R/\mathfrak{m}}) =
\bigcup^{n}_{i=0}\Supp\mathbb S_{i}$, where the $\mathbb S_{i}$ are
the simple functors $(\uGproj R)^{\op} \lxr \ab$ appearing in the
composition series of $\uHom_{R}(-,X_{R/\mathfrak{m}})$. Since
$\vert \Supp \mathbb S_{i}\rvert < \infty$, $\forall i$, we infer
that $\vert\mathsf{Ind}\uGproj R\rvert < \infty$. This  implies that
$\Gproj R$ is of finite representation type, i.e. $R$ is of finite
CM-type.

(i) $\Leftrightarrow$ (iii) The proof follows as in Theorem $4.18$
 using Corollary $3.5$.

(iv) $\Leftrightarrow$ (ii) By recent results of J{\o}rgensen and Holm,
see \cite[Theorem 2.7]{HJ}, the first part of condition (iv) implies
that $R$ is Gorenstein, so (iv) $\Rightarrow$ (ii). If (ii) holds,
then by \cite[Theorem 10.3.8]{EJ} any Gorenstein-flat module is a
filtered colimit of finitely generated Gorenstein-projective
modules, so (iv) holds.
\end{proof}
\end{thm}

\begin{cor} Let $R$ be a commutative Noetherian complete local Gorenstein
ring.
\begin{enumerate}
\item If $\Gproj R = \proj R$, then $R$ is non-singular, i.e. regular, of finite
CM-type.
\item If  $\Gproj R \neq \proj R$ and any
Gorenstein-projective module is a direct sum of finitely generated
modules,  then $R$ is a simple (hypersurface) singularity of finite CM-type.
\end{enumerate}
\end{cor}

\section{Stable and Torsion-Free Modules} Let $\Lambda$ be a
Noetherian ring. Recall that a finitely generated right
$\Lambda$-module $X$ is called {\bf stable} if
$\Ext^{n}_{\Lambda}(X,\Lambda) = 0$, $\forall n \geq 1$, i.e. $X \in
{^{\bot}}\Lambda$. And $X$ is called {\bf torsion-free} if the
transpose $\mathsf{Tr}X$ of $X$ is stable:
$\Ext^{n}_{\Lambda^{\op}}(\mathsf{Tr} X,\Lambda^{\op}) = 0$,
$\forall n \geq 1$, i.e. $\mathsf{Tr}X \in {^{\bot}}\Lambda^{\op}$.
Here $\mathsf{Tr} \colon \umod\Lambda \lxr \umod\Lambda^{\op}$ is
the Auslander-Bridger transpose duality functor given by
$\mathsf{Tr}(M) = \Coker\big(\Hom_{\Lambda}(P^{0},\Lambda) \lxr
\Hom_{\Lambda}(P^{1},\Lambda)\big)$, where $P^{1} \lxr P^{0} \lxr M
\lxr 0$ is a finite projective presentation of $M$, see \cite{ABr}.
In this terminology a finitely generated module $X$ is
Gorenstein-projective if and only if $X$ is stable and torsion-free.

It was an open question if stable implies torsion-free, i.e. if
stableness is sufficient for a finitely generated module  to be
Gorenstein-projective. Recently Jorgensen and Sega answered this
question in the negative. In fact they constructed a commutative
Noetherian local (graded Koszul) ring and a finitely generated
stable module which is not torsion-free, see \cite{JS} for details.
On the other hand Yoshino proved that the answer to the above
question is positive for Henselian local Noetherian  rings
$(R,\mathfrak{m},k)$ such that the full subcategory of stable
modules is of finite representation type, see \cite{Yoshino}.
Recently R. Takahashi proved a far reaching generalization: if $\X$
is a contravariantly finite resolving subcategory of $\smod R$ and
if any module $X$ in $\X$ satisfies $\Ext^{n}_{R}(X,R) = 0$, for $n
\gg 0$, then $\X = \proj R$ or $R$ is Gorenstein and either $\X =
\Gproj R$, or $\X = \smod R$, see \cite{Takahashi}. In the
non-commutative setting, where in the Artin algebra case Takahashi's
result is  false (e.g. consider $\X = \proj^{<\infty}\Lambda$ for a
Gorenstein algebra $\Lambda$ of infinite global dimension), it is
known that all stable modules are Gorenstein-projective if the Artin
algebra $\Lambda$ is Gorenstein or more generally if
$\id{_{\Lambda}}\Lambda < \infty$, see \cite[Proposition 4.4]{B:cm}.
Also any stable module with finite Gorenstein-projective dimension
is Gorenstein-projective, see \cite[Proposition 3.9]{B:cm}, and all
stable right modules are Gorenstein-projective if $\Lambda$ is left
CoGorenstein in the sense that arbitrary syzygy left modules are
stable, see \cite[Proposition 4.4]{B:cm}.

\subsection{Finitely generated stable modules} One reason of why the condition ${^{\bot}}\Lambda =
\Gproj\Lambda$ is important is that it  allows the computation of
the {\em Gorenstein dimension} $\Gor\bdim M := \rd_{\Gproj\Lambda}M$
of a module $M \in \smod\Lambda$, without using resolutions. Recall
that the {\em Gorenstein dimension} of $\smod\Lambda$ is defined by
\[
\Gor\bdim\Lambda = {\rd}_{\Gproj\Lambda}\smod\Lambda =
\bmax\big\{\id{_{\Lambda}}\Lambda,\id\Lambda_{\Lambda}\big\}
\]
Thus $\Lambda$ is Gorenstein if and only if $\Gor\bdim\Lambda <
\infty$.

\begin{lem} Let $\Lambda$ be a Noetherian  ring. Then the following
are equivalent.
\begin{enumerate}
\item ${^{\bot}}\Lambda = \Gproj\Lambda$.
\item $\forall M \in \smod\Lambda$: \ $\Gor\bdim M = \sup\{t \in \mathbb Z \, | \, \Ext^{t}_{\Lambda}(M,\Lambda) \neq
0\}$.
\end{enumerate}
\begin{proof} Assume that ${^{\bot}}\Lambda =
\Gproj\Lambda$ and let $M$ be in $\smod\Lambda$. If $\Gor\bdim M
<\infty$, then the assertion (ii) is easy, see for instance
\cite[Proposition 3.9]{B:cm}. Assume that $\Gor\bdim M = \infty$. If
$\sup\{t \in \mathbb Z\, | \, \Ext^{t}_{\Lambda}(M,\Lambda) \neq 0\}
= d < \infty$, then $\Ext^{n}_{\Lambda}(M,\Lambda) = 0$, $\forall n
> d$. Hence $\Ext^{k}_{\Lambda}(\Omega^{d}M,\Lambda) = 0$, $\forall
k \geq 1$, so $\Omega^{d}M \in {^{\bot}}\Lambda = \Gproj\Lambda$.
However this implies that $\Gor\bdim M \leq d$ and this is not true.
Hence $d = \infty$. Conversely if $\Gor\bdim M = \sup\{t \in \mathbb
Z\, | \, \Ext^{t}_{\Lambda}(M,\Lambda) \neq 0\}$, then $\Gor\bdim M
= 0$, for any module $M \in {^{\bot}}\Lambda$. Hence
${^{\bot}}\Lambda = \Gproj\Lambda$.
\end{proof}
\end{lem}

 Our aim in this section is
to prove the following result which, in the context of virtually
Gorenstein  Artin algebras, gives a complete answer to the question
of when a stable module is torsion-free, hence
Gorenstein-projective. In particular we extend properly Yoshino's
result in the non-commutative setting to the case where the category
of stable modules is not necessarily of finite representation type.

\begin{thm} Let $\Lambda$ be an Artin algebra.
If $\Gproj\Lambda$ is contravariantly finite, for instance if
$\Lambda$ is virtually Gorenstein or of finite CM-type, then the following are
equivalent.
\begin{enumerate}
  \item ${^{\bot}}\Lambda = \Gproj\Lambda$.
 \item ${^{\bot}}\Lambda \cap (\Gproj\Lambda)^{\bot} =
  \proj\Lambda$.
  \item ${^{\bot}}\Lambda \cap (\Gproj\Lambda)^{\bot}$ is of finite
  representation type.
\end{enumerate}
\end{thm}

Before we proceed to the proof of  Theorem $5.2$ we need some
preparations. Let $\A$ be an abelian category with enough
projectives and let $\X$ be a resolving subcategory of $\A$. We view
the stable category $\un\A$ modulo projectives as a left
triangulated category; then clearly the stable category $\un\X$ of
$\X$ is a left triangulated subcategory of  $\un\A$, see \cite{BM},
so $\un\X$ has weak kernels and therefore the category of coherent
functors $\smod\un\X$ over $\un\X$ is abelian with enough
projectives. Moreover the Yoneda embedding $\un\X \lxr \smod\un\X$,
$\unX \longmapsto \uHom(-,X)$ has image in $\Proj\smod\un\X$ and if
idempotents split in $\un\X$, it induces an equivalence $\un\X
\approx \Proj\smod\un\X$. In the sequel we identify $\smod\un\X$
with the full subcategory of $\smod\X$ consisting of all coherent
functors $F : \X^{\op} \lxr \ab$ which admit a presentation
$\X(-,X_{1}) \lxr \X(-,X_{0}) \lxr F \lxr 0$ such that the map
$X_{1} \lxr X_{0}$ is an epimorphism in $\A$. This clearly coincides
with the full subcategory of $\smod\X$ consisting of all functors
vanishing on the projectives.

\begin{lem} Let $\A$ be an abelian category and $(\X,\Y)$ a
cotorsion pair in $\A$ with heart $\omega = \X \cap \Y$.
\begin{enumerate}
\item The category $\smod\X$ is abelian with enough projectives and the restricted Yoneda
functor $\mathsf{H} : \A \lxr \Mod\X$, $\mathsf{H}(A) =
\A(-,A)|_{\X}$ has image in $\smod\X$.
\item The category $\X/\omega$ is right triangulated and the category
$(\X/\omega\lsmod)^{\op}$ is abelian with enough injectives. If
idempotents split in $\X/\omega$, then we have an equivalence
$\X/\omega \approx \Inj(\X/\omega\lsmod)^{\op}$.
\item If $\A$ has enough projectives, then $\smod\un\X$ has enough
projectives and enough injectives. Moreover if idempotents split in
$\X/\omega$, then the assignment $X \longmapsto
\Ext^{1}_{\A}(-,X)|_{\X}$  gives rise to an equivalence
\[
\X/\omega \,\ \approx \,\ \Inj\smod\underline{\X}
\]
\end{enumerate}
\begin{proof} (i), (ii) Since $\X$ is contravariantly finite in $\A$
 it follows that $\X$ has weak kernels and therefore $\smod\X$ is abelian. For
 any object $A \in \A$, let $X^{0}_{A} \lxr A$ be a right
 $\X$-approximation of $A$, and let $X^{1}_{A} \lxr K$ be a right $\X$-approximation of the kernel $K$
 of $X^{0}_{A} \lxr A$. Applying the functor $\mathsf{H}$ we have
 clearly an exact sequence $\mathsf{H}(X^{1}_{A}) \lxr
 \mathsf{H}(X^{0}_{A}) \lxr \mathsf{H}(A) \lxr 0$ which shows that
 $\mathsf{H}(A)$ is coherent. Hence $\Image\mathsf{H} \subseteq
 \smod\X$. By \cite[Lemma VI.1.1]{BM} the stable category $\X/\omega$ is
right triangulated, in particular $\X/\omega$ has weak cokernels.
Therefore $\X/\omega\lsmod$ is abelian with enough projectives and
part (ii) follows.

(iii)  Consider the Yoneda embedding $\mathsf{H} : \X \lxr \smod\X$,
$\mathsf{H}(X) = \X(-,X)$. Since $\omega$ is an Ext-injective
cogenerator of $\X$, we may choose, $\forall X\in \X$, an
$\omega$-{\em injective copresentation} of $X$, i.e. an exact
sequence $0 \lxr X \lxr T \lxr X^{\prime} \lxr 0$, where $T\in
\omega$ and $X^{\prime}\in \X$. Then we have a projective resolution
in $\smod\X$:
\[
0 \lxr \mathsf{H}(X) \lxr \mathsf{H}(T) \lxr \mathsf{H}(X^{\prime})
\lxr {\Ext}^{1}_{\A}(-,X)|_{\X} \lxr 0
\]
This shows that $\Ext^{1}_{\A}(-,X)|_{\X}$ is coherent and lies in
$\smod\un\X$, $\forall X \in \X$.  We show that
$\Ext^{1}_{\A}(-,X)|_{\X}$ is injective in $\smod\un\X$. Let $F$ be
in $\smod\un\X$ with a presentation $\X(-,X_{1}) \lxr \X(-,X_{0})
\lxr F \lxr 0$ in $\smod\X$ such that the map $X_{1} \lxr X_{0}$ is
an epimorphism. Since $\X$ is resolving, we have an exact sequence
$0 \lxr X_{2} \lxr X_{1} \lxr X_{0} \lxr 0$ in $\X$ which induces a
projective resolution of $F$ in $\smod\X$:
\begin{equation}
0 \, \lxr \,  \mathsf{H}(X_{2}) \, \lxr \, \mathsf{H}(X_{1}) \, \lxr \,
\mathsf{H}(X_{0}) \, \lxr \, F \, \lxr \, 0
\end{equation}
and a triangle $\Omega \underline{X}_{0} \lxr \underline{X}_{2} \lxr
\underline{X}_{1} \lxr \underline{X}_{0}$ in $\un\X$ which induces a
projective resolution of $F$ in $\smod\un\X$:
 \[
 \cdots \, \lxr \, \un\X(-,\Omega \underline{X}_{0}) \, \lxr \, \un\X(-,\underline{X}_{2}) \, \lxr \,
\un\X(-,\underline{X}_{1}) \, \lxr \, \un\X(-,\underline{X}_{0}) \, \lxr \, F
\, \lxr \,  0\] Then
$\Ext^{1}_{\smod\un\X}\bigl(F,\Ext^{1}_{\A}(-,X)|_{\X}\bigr)$ is the
homology of the complex
\[
{\Ext}^{1}_{\A}(X_{0},X) \ \lxr \
{\Ext}^{1}_{\A}(X_{1},X) \ \lxr \ {\Ext}^{1}_{\A}(X_{2},X)
\]
 which is exact.
Hence $\Ext^{1}_{\smod\un\X}\bigl(F,\Ext^{1}_{\A}(-,X)|_{\X}\bigr) =
0$ and this implies that $\Ext^{1}_{\A}(-,X)|_{\X}$ is injective in
$\smod\un\X$. Let $\alpha : X_{1} \lxr X_{2}$ be a map in $\X$ and
$0 \lxr X_{i} \lxr T_{i} \lxr X^{\prime}_{i} \lxr 0$ be
$\omega$-injective copresentations of the $X_{i}$. Since
$\Ext^{1}_{\A}(X^{\prime}_{1},T_{2}) = 0$, there exists an exact
commutative diagram in $\X$:
\[
\begin{CD}
 0 @>  >> X_{1} @>  >> T_{1} @>   >> X^{\prime}_{1} @> >> 0\\
   &\ & @V{\alpha}VV @V{\beta}VV @V{\gamma}VV & \\
0 @>  >> X_{2} @> >> T_{2} @>  >> X^{\prime}_{2} @> >> 0\\
\end{CD}
\]
which induces an exact commutative diagram in $\smod\X$:
\[
\begin{CD}
 0 @>  >> \mathsf{H}(X_{1}) @>  >> \mathsf{H}(T_{1}) @>   >> \mathsf{H}(X^{\prime}_{1}) @> >>
 {\Ext}^{1}_{\A}(-,X_{1})|_{\X} @> >> 0\\
   &\ & @V{\mathsf{H}(\alpha)}VV @V{\mathsf{H}(\beta)}VV @V{\mathsf{H}(\gamma)}VV @V{\alpha^{*}}VV  \\
0 @>  >> \mathsf{H}(X_{2}) @>  >> \mathsf{H}(T_{2}) @>   >>
\mathsf{H}(X^{\prime}_{2}) @> >> {\Ext}^{1}_{\A}(-,X_{2})|_{\X} @> >> 0\\
\end{CD}
\]
By diagram chasing it is easy to see that the objects
${\Ext}^{1}_{\A}(-,X_{i})|_{\X}$ are independent of the chosen
$\omega$-injective copresentations  of the $X_{i}$, and the map
$\alpha^{*}$ is independent of the liftings $\beta$ and $\gamma$ of
$\alpha$. In this way we obtain a functor $\mathsf{H}^{*} : \X \lxr
\smod\underline{\X}$, $\mathsf{H}^{*}(X) =
\Ext^{1}_{\A}(-,X)|_{\X}$, and $\mathsf{H}^{*}(\alpha)= \alpha^{*}$.
Since $\mathsf{H}^{*}(\omega) = 0$, $\mathsf{H}^{*}$ induces a
functor $\mathsf{H}^{\vee} : \X/\omega \lxr \smod\underline{\X}$,
and by the above analysis we have $\Image\mathsf{H}^{\vee} \subseteq
\Inj\smod\underline{\X}$.  If $\alpha^{*} = 0$, then chasing the
above diagram we see that $\alpha$ factors through $X_{1} \lxr
T_{1}$ and therefore $\alpha$ is zero in $\X/\omega$, i.e.
$\mathsf{H}^{\vee}$ is faithful. In the same way, using that
$\mathsf{H}$ is fully faithful and its image consists of projective
objects of $\smod\X$, we see that $\mathsf{H}^{\vee}$ is full. Let
$F$ be an injective object in $\smod\underline{\X}$ with projective
resolution in $\smod\X$ as in $(5.1)$.  Then the induced long exact
sequence
\[
0 \, \lxr \, \mathsf{H}(X_{2}) \, \lxr \, \mathsf{H}(X_{1}) \, \lxr \,
\mathsf{H}(X_{0}) \, \lxr \, {\Ext}^{1}_{\A}(-,X_{2})|_{\X} \, \lxr \, \cdots
\]
shows that we have a monomorphism $F \monic
\Ext^{1}_{\A}(-,X_{2})|_{\X}$ in $\smod\underline{\X}$ which splits
since $F$ is injective. This produces an idempotent endomorphism $f$
of $\Ext^{1}_{\A}(-,X_{2})|_{\X}$ which clearly comes from an
idempotent endomorphism $\underline{e} : \unX_{2} \lxr \unX_{2}$ in
$\X/\omega$. By our hypothesis on $\X/\omega$, the idempotent
$\underline{e}$ splits: there is an object $X_{3} \in \X$ and
morphisms $\underline{\kappa} : \underline{X}_{2} \lxr
\underline{X}_{3}$ and $\underline{\lambda} : \underline{X}_{3} \lxr
\underline{X}_{2}$ such that $\underline{e} = \underline{\kappa}
\circ \underline{\lambda}$ and $1_{\underline{X}_{3}}  =
\underline{\lambda} \circ \underline{\kappa}$. This clearly implies
that $F$ is isomorphic to $\mathsf{H}^{\vee}(\underline{X}_{3})$ and
therefore $\mathsf{H}^{\vee} : \X/\omega \lxr
\Inj\smod\underline{\X}$ is an equivalence.
\end{proof}
\end{lem}

We have the following consequence, see also \cite[Proposition
4.6]{Buchweitz} for a related result.

\begin{cor} Let $\A$ be an abelian category and $(\X,\Y)$ a
cotorsion pair in $\A$ with heart $\omega = \X \cap \Y$.
\begin{enumerate}
\item If $\A$ has enough projectives, then there is an equivalence:
\[
 \smod\un\X \,\, \stackrel{\approx}{\lxr} \,\,
(\X/\omega\lsmod)^{\op}
\]
\item If $\A$ has enough injectives, then there is an equivalence:
\[
 \overline{\Y}\lsmod \,\,
\stackrel{\approx}{\lxr} \,\, (\smod\Y/\omega)^{\op}
\]
\end{enumerate}
\begin{proof} By Lemma $5.3$, the abelian category $\smod\un\X$ has enough injectives and we have a full embedding $\mathsf{H}^{*} : \X/\omega
\lxr \Inj\smod\un\X$ and $\Inj\smod\un\X =
\add\Image\mathsf{H}^{*}$. Taking coherent functors we have
 $\smod\un\X \approx \bigl((\Inj\smod\un\X)\lsmod\bigr)^{\op} \approx
 \bigl((\add\Image\mathsf{H}^{\vee})\lsmod\bigr)^{\op} \approx
 (\X/\omega\lsmod)^{\op}$. The second part is proved in a dual way.
\end{proof}
\end{cor}

\begin{exam} $(\alpha)$ If $\A$ has enough projectives and enough injectives,
then considering the cotorsion pairs $(\A,\Inj\A)$ and
$(\Proj\A,\A)$ in $\A$, we deduce an equivalence
\[
\smod\un\A \  \approx \
(\overline{\A}\lsmod)^{\op}
\]

$(\beta)$ Let $T$, resp. $S$, be a cotilting, resp. tilting, module
over an Artin algebra $\Lambda$. Since we have cotorsion pairs
$({^{\bot}}T,\Y)$ and $(\X,S^{\bot})$ in $\smod\Lambda$,  we deduce
equivalences:
\[
\smod\underline{{^{\bot}}T} \ \approx \
\{({^{\bot}}T/\add T)\lsmod\}^{\op} \ \ \ \ \  \text{and} \ \ \ \ \  \overline{S^{\bot}}\lsmod
\ \approx \ \{\smod\ \!(S^{\bot}/\add S)\}^{\op}
\]
 If $\id T \leq 1$,
resp. $\pd S \leq 1$, then it is easy to see that $\Y/\add T \approx
\Z := \{A \in \smod\Lambda \, | \, \Hom_{\Lambda}(A,T) = 0\}$, resp.
$\X/\add S \approx \W := \{A \in \smod\Lambda \, | \,\
\Hom_{\Lambda}(S,A) = 0\}$, see \cite[Proposition V.5.2]{BR}. It
follows that in this case we have equivalences:
\[
\overline{\Y}\lsmod
\ \approx \ (\smod\Z)^{\op} \ \ \ \ \  \text{and} \ \ \ \ \  \smod\un\X \ \approx \ (\W\lsmod)^{\op}
\]
\end{exam}

For the proof of the following we refer to \cite[2.43]{ABr},
\cite[Theorem 3.3]{ABM} and \cite[Corollary 4.9]{B:freyd}.

\begin{lem} Let $\X$ be a resolving subcategory of an abelian category
$\A$ with enough projectives. If $\X \subseteq
{^{\bot}}\{\Proj\A\}$, then the loop functor
  $\Omega : \un\X \lxr \un\X$ is fully faithful
and any projective object in  $\smod\un\X$ is injective.
\end{lem}

\begin{thm} Let $\A$ be an abelian category and $(\X,\Y)$  a
cotorsion pair in $\A$ with heart $\omega = \X \cap \Y$.
\begin{enumerate}
\item Assume that  $\A$ has enough projectives and $\X \subseteq
{^{\bot}}\Proj\A$. Then $\X \subseteq \GProj\A$ provided that one of
the following conditions holds:
 \begin{enumerate}
 \item $\rd_{\X}\A < \infty$.
\item $\X$ is of finite representation type and idempotents split in
$\X/\omega$.
 \end{enumerate}
\item Assume that $\A$ has enough injectives and $\Y \subseteq (\Inj\A)^{\bot}$.
Then $\Y \subseteq \GInj\A$ provided that one of the following
conditions holds:
\begin{enumerate}
 \item $\cd_{\Y}\A < \infty$.
\item $\Y$ is of finite representation type and idempotents split in
$\Y/\omega$.
 \end{enumerate}
\end{enumerate}
\begin{proof} We only prove (i) since part (ii) follows by duality.

Since $\X \subseteq {^{\bot}}\Proj\A$, it follows that $\Proj\A
\subseteq \Y$ and therefore $\Proj\A \subseteq \omega$. We first show that in both
cases  the conclusion $\X \subseteq \GProj\A$
follows from the equality $\omega = \Proj\A$. Indeed if this
equality holds, then since $\omega$, as the heart of the cotorsion
pair $(\X,\Y)$, is an injective cogenerator of $\X$, it follows that
for any object $M \in \X$ there exists a short exact sequence $0\lxr
M \lxr P \lxr M^{\prime} \lxr 0$, where $P$ is projective and
$M^{\prime} \in \X$. Since $\X \subseteq {^{\bot}}\Proj\A$, this
clearly implies that $M$ is Gorenstein-projective. Hence $\X
\subseteq \GProj\A$.

(a) If $\rd_{\X}\A = d <\infty$, then by \cite[Proposition
VII.1.1]{BR} there exists a cotorsion pair $(\X,\Proj^{<\infty}\A)$
with heart $\Proj\A$. It follows that $\Y = \Proj^{<\infty}\A$ and then clearly
$\omega = \Proj\A$.

(b) Since $\X$ is a Krull-Schmidt category of finite representation
type, it follows that so are $\Proj\A$ and $\omega$. Assume that
$\omega\setminus\Proj\A \neq \emptyset$. Let $\{T_{1}, T_{2},
\cdots, T_{l}\}$, where $l \geq 1$, be the set of indecomposable non-projective objects
of $\omega$ and let $\{X_{1},X_{2},\cdots, X_{k}\}$ be the set of
indecomposable objects of $\X$ which are not in $\omega$. Note that
if $\X = \omega$, then it is easy to see that $\X = \Proj\A
\subseteq \GProj\A$, so we may exclude this trivial case. Set $T =
\oplus^{l}_{i=1}T_{i}$, $S = \oplus^{k}_{i=1}X_{i}$, and $X =
T\oplus S$.  Then $\X = \add(P\oplus X)$, where $P$ is the direct
sum of the set of non-isomorphic indecomposable projective objects
of $\A$. Consider the stable categories $\un\X$ and $\X/\omega$.
Then clearly $\un\X = \add \underline{X}$ and $\X/\omega = \add
\underline{S}$. Since idempotents split in $\X/\omega$, by Lemma
$5.3$ there is an equivalence $\Inj\smod\un\X \approx \X/\omega$. On
the other hand, since $\X \subseteq {^{\bot}}\Proj\A$, by Lemma
$5.6$ it follows that any projective object of  $\smod\un\X$ is
injective, i.e. $\Proj\smod\un\X \subseteq \Inj\smod\un\X$. Putting
things together we infer the existence of a full embedding $\un\X
\lxr \X/\omega$ which is the composite $\un\X \stackrel{\approx}{\lxr}
\Proj\smod\un\X \lxr \Inj\smod\un\X \stackrel{\approx}{\lxr} \X/\omega$.
Therefore  $k + l = \lvert \un\X \rvert \leq \lvert \X/\omega
\rvert = k$. This contradiction shows that  $l = 0$, i.e. any indecomposable object
of the heart $\omega$ is projective. We infer that $\omega = \Proj\A$.
\end{proof}
\end{thm}

\begin{rem} Let $\A$ be an abelian Krull-Schmidt category with enough
projectives and enough injectives. If $\X$ is a contravariantly
finite resolving subcategory of $\A$, then it is easy to see that
there exists a cotorsion pair $(\X,\Y)$ in $\A$, where $\Y =
\X^{\bot}$. Hence Theorem $5.7$ applies to resolving subcategories
of finite representation type in $\A$ such that idempotents split in
$\X\big/\X\cap \X^{\bot}$ and $\X \subseteq {^{\bot}}\Proj\A$.
\end{rem}

A full subcategory $\X$ of an abelian category $\A$ is called {\bf
projectively thick} if $\X$ is resolving and closed under direct
summands and cokernels of $\Proj\A$-monics. By \cite{B:cm}, the map
$\X \longmapsto \un\X$ gives a bijection between projectively thick
subcategories of $\GProj\A$, and thick subcategories of $\uGProj\A$.
The following gives a non-commutative generalization of the above
mentioned  result of Yoshino, see \cite[Theorem 5.5]{Yoshino}.

\begin{prop} Let $\Lambda$ be an Artin algebra and $\X$ a resolving subcategory of
$\smod\Lambda$. Then the following statements are equivalent.
\begin{enumerate}
  \item $\X$ is of finite representation type and $\X \subseteq {^{\bot}}\Lambda$.
  \item
  \begin{enumerate}
  \item $\X$ is contravariantly finite in $\smod\Lambda$,
  \item $\X$ is a projectively thick
  subcategory of $\Gproj\Lambda$, and
  \item The functor category $\Mod\un\X$ is Frobenius.
  \end{enumerate}
  \item
  \begin{enumerate}
  \item $\X$ is contravariantly finite in $\smod\Lambda$,
  \item $\un\X$ is a triangulated subcategory of $\umod\Lambda$, and
  \item The stable triangulated category $\colim\X/\colim\proj\Lambda$ is
phantomless.
\end{enumerate}
  \end{enumerate}
  In particular if $\X$ is a representation-finite resolving subcategory of ${^{\bot}}\Lambda$,
  then $\X\subseteq \Gproj\Lambda$,  $\smod\un\X$ and $\Mod\un\X$ are Frobenius, and there is an equivalence
    $\colim\{\X/\proj\Lambda\} \approx
  \colim\X\big/\colim\proj\Lambda$.
\begin{proof} (i) $\Rightarrow$ (ii) Clearly $\X$ is contravariantly
finite in $\smod\Lambda$. Since $\smod\Lambda$ is Krull-Schmidt with
enough injectives, by Remark $5.8$, there exists a cotorsion pair
$(\X,\Y)$ in $\smod\Lambda$, where $\Y = \X^{\bot}$, and idempotents
split in the stable category $\X/\X\cap \Y$. Then $\X \subseteq
\Gproj\Lambda$ by Theorem $5.7$. Since $\Mod\un\X \approx
\Mod\Gamma$ where $\Gamma$ is the stable endomorphism algebra of a
representation generator of $\X$, it follows that $\Mod\un\X$ is
Frobenius. Finally as in \cite[Proof of Proposition 3.8(iii)]{B:cm}
we see that $\X$ is projectively thick.

(ii) $\Rightarrow$ (i) Clearly  $\X \subseteq {^{\bot}}\Lambda$ and
the stable category $\un\X$ is a thick triangulated subcategory of
$\uGproj\Lambda$. By \cite[Theorem 9.4]{B:cm} the stable category
$\colim\X\big/\colim\proj\Lambda$
 is triangulated and compactly generated, and  we have
an equivalence $ \un\X \approx
\{\colim\X\big/\colim\proj\Lambda\}^{\mathsf{cpt}}$. Hence the
representation category of $\colim\X\big/\colim\proj\Lambda$ is
equivalent to the module category $\Mod\un\X$ which by hypothesis is
Frobenius; this clearly implies that $\Mod\un\X$ is perfect. On the
other hand since $\X$ is contravariantly finite, as in the proof of
Theorem $3.1$, it follows that $\X$ is a dualizing variety. Then by
\cite{AS:subcategories} so is $\un\X$ and consequently $\un\X$ has
left almost split maps. This implies by \cite[Theorem 10.2]{B:art}
that $\Mod\un\X$ is locally finite. Then a similar argument as in
the proof of Theorem $3.1$ shows that $\un\X$, or equivalently $\X$,
is of finite representation type.

(ii) $\Leftrightarrow$ (iii) This follows from the following facts:
$(\alpha)$ $\un\X$ is a triangulated subcategory of $\umod\Lambda$
if and only if $\X$ is projectively thick and consists of
Gorenstein-projective modules \cite{B:cm}, and  $(\beta)$
$\colim\X\big/\colim\proj\Lambda$ is compactly generated and
phantomless if and only if the module category $\Mod\un\X$ is
Frobenius \cite{B:3cats}.
\end{proof}
\end{prop}

If $\A$ is an abelian category, we let $\Omega^{\infty}(\A)$ be the
full subcategory of $\A$ consisting of arbitrary syzygy objects,
i.e. objects $A$ admitting  an exact sequence $0 \lxr A \lxr P^{0}
\lxr P^{1} \lxr \cdots$, with $P^{n} \in\Proj \A$, $\forall n \geq
0$.

\begin{cor} Let $\Lambda$ be an Artin algebra and let $\X$ be a
resolving subcategory of $\smod\Lambda$. If $\X \subseteq
{^{\bot}}\Lambda$ and $\Omega^{\infty}(\smod\Lambda^{\op})$ is of
finite representation type, then $\X \subseteq \Gproj\Lambda$.
\begin{proof} Since $\X \subseteq {^{\bot}}\Lambda$, it follows
directly that $\mathsf{Tr}(\X) \subseteq
\Omega^{\infty}(\smod\Lambda^{\op})$. Since $\mathsf{Tr} :
\umod\Lambda \lxr \umod\Lambda^{\op}$ is a duality, we have a fully
faithful functor $\mathsf{Tr} : \un\X \lxr
\Omega^{\infty}(\umod\Lambda^{\op})$. It follows that $\un\X$ is of
finite representation type and therefore $\X \subseteq
\Gproj\Lambda$ by Proposition $5.9$.
\end{proof}
\end{cor}

Now we can give the proof of Theorem $5.2$.

\medskip

{\em Proof of Theorem 5.2}.  (i) $\Rightarrow$ (ii) $\Rightarrow$
(iii) If ${^{\bot}}\Lambda = \Gproj\Lambda$, then ${^{\bot}}\Lambda
\cap (\Gproj\Lambda)^{\bot} = \Gproj\Lambda \cap
(\Gproj\Lambda)^{\bot}$. Clearly any finitely generated projective
module lies in $\Gproj\Lambda \cap (\Gproj\Lambda)^{\bot}$. If $X
\in \Gproj\Lambda \cap (\Gproj\Lambda)^{\bot}$, then since $X$ is
Gorenstein-projective, there is a short exact sequence $0 \lxr X
\lxr P \lxr X^{\prime} \lxr 0$, where $P \in \proj\Lambda$  and
$X^{\prime} \in \Gproj\Lambda$. This extension splits since $X \in
(\Gproj\Lambda)^{\bot}$, so $X$ is projective as a direct summand of
$P$. Hence $\Gproj\Lambda \cap (\Gproj\Lambda)^{\bot} =
\proj\Lambda$ which is always of finite representation type.

(iii) $\Rightarrow$ (i) Since by \cite{B:cm}, the subcategory
$(\Gproj\Lambda)^{\bot}$ is resolving, it follows that so is
${^{\bot}}\Lambda \cap (\Gproj\Lambda)^{\bot}$. Hence
${^{\bot}}\Lambda \cap (\Gproj\Lambda)^{\bot}$  is a resolving
subcategory of $\smod\Lambda$ of finite representation type. Then by
Proposition $5.9$ we have ${^{\bot}}\Lambda \cap
(\Gproj\Lambda)^{\bot} \subseteq \Gproj\Lambda$. Let $T$ be in
$\smod\Lambda$ such that $\add T = {^{\bot}}\Lambda \cap
(\Gproj\Lambda)^{\bot}$. Then $T$ is Gorenstein-projective and since
$T$ lies in $(\Gproj\Lambda)^{\bot}$ we infer that $T$ is
projective. On the other hand since any projective module lies in
both ${^{\bot}}\Lambda$ and $(\Gproj\Lambda)^{\bot}$, we infer that
${^{\bot}}\Lambda \cap (\Gproj\Lambda)^{\bot} = \proj\Lambda$. Since
$\Gproj\Lambda$ is contravariantly finite, there is a cotorsion pair
$(\Gproj\Lambda,(\Gproj\Lambda)^{\bot})$ in $\smod\Lambda$. Hence
$\forall A \in {^{\bot}}\Lambda$ there exists a short exact sequence
$0 \lxr A \lxr Y^{A} \lxr X^{A} \lxr 0$, where $Y^{A} \in
(\Gproj\Lambda)^{\bot}$ and $X^{A}\in \Gproj\Lambda$. Clearly then
$Y^{A}$ lies in ${^{\bot}}\Lambda \cap (\Gproj\Lambda)^{\bot}$,
hence $Y^{A}$ is projective and therefore $A \cong \Omega X^{A}$.
Since $X^{A}$ is Gorenstein-projective, this implies that so is $A$.
We infer that ${^{\bot}}\Lambda = \Gproj\Lambda$.  $\Box$

\begin{cor} For an Artin algebra $\Lambda$, the following are equivalent.
\begin{enumerate}
  \item ${^{\bot}}\Lambda$ is of finite representation type.
  \item $\Lambda$ is of finite CM-type and ${^{\bot}}\Lambda =
  \Gproj\Lambda$.
\end{enumerate}
\end{cor}

\subsection{Infinitely generated stable modules} Now we turn our attention to the question of when
infinitely generated stable modules, i.e. modules in
${^{\pmb{\bot}}}\Lambda = \{A \in \Mod\Lambda \, | \,
\Ext^{n}_{\Lambda}(A,\Lambda) = 0, \,\ \forall n \geq 1 \}$, are
Gorenstein-projective.

\begin{prop} If $\Lambda$ is an Artin algebra then the following are equivalent.
\begin{enumerate}
\item $\Lambda$ is virtually Gorenstein and $\GProj\Lambda =
{^{\pmb{\bot}}}\Lambda$.
\item Any module in ${^{\pmb{\bot}}}\Lambda$ $\cap
(\Gproj\Lambda)^{\pmb\bot} $ is a direct sum of finitely generated
modules.
\item
\begin{enumerate}
\item ${^{\bot}}\Lambda \cap (\Gproj\Lambda)^{\bot}$ is of finite
representation type.
\item Any finitely generated module admits a left
$\{ {^{\pmb{\bot}}}\Lambda$ $\cap$
$(\Gproj\Lambda)^{\pmb\bot}\}$-approximation which is finitely
generated.
\end{enumerate}
\item ${^{\pmb{\bot}}}\Lambda$ $\cap$
$(\Gproj\Lambda)^{\pmb{\bot}}$ is a pure-semisimple locally finitely
presented category.
\end{enumerate}
If $\mathrm{(i)}$ holds, then: ${^{\bot}}\Lambda = \Gproj\Lambda$.
\begin{proof}  Set $\mathcal H = {^{\pmb{\bot}}}\Lambda
\cap$ $(\Gproj\Lambda)^{\pmb\bot}$ and $\omega = \mathcal H \cap
\smod\Lambda$. Then clearly $\omega = {^{\bot}}\Lambda \cap
(\Gproj\Lambda)^{\bot}$.

(i) $\Rightarrow$ (ii) This follows from Proposition $4.7$ since in this
case we have $(\Gproj\Lambda)^{\bbot} = (\GProj\Lambda)^{\bot} =
\colim(\Gproj\Lambda)^{\bot}$.

(ii) $\Rightarrow$ (iii) We have $\mathcal H = \colim\omega =
\Add\omega$. Clearly then any finitely generated module admits a
finitely generated left $\mathcal H$-approximation. Since $\omega$
is resolving, by Theorem $3.1$, $\omega$ is of finite representation
type.

(iii) $\Rightarrow$ (iv) Since the subcategories
${^{\pmb{\bot}}}\Lambda$ and $(\Gproj\Lambda)^{\pmb{\bot}}$ are
definable and resolving, it follows that so is $\mathcal H$. Then by
\cite{Krause:memoirs} condition (b) implies that $\mathcal H =
\colim\omega$, so $\mathcal H$ is locally finitely presented. On the
other hand condition (a) together with Theorem $3.1$  implies that
$\mathcal H$ is pure-semisimple.

(iv) $\Rightarrow$ (i) Since $\omega$ is resolving, by Theorem $3.1$
it follows that $\omega = \fp\mathcal H$ is of finite representation
type. Then as in the proof of Theorem $5.2$ we infer that $\omega =
\proj\Lambda$ and therefore $\mathcal H = \colim\proj\Lambda =
\Proj\Lambda$. Consider the cotorsion pair
$(\A,(\Gproj\Lambda)^{\pmb{\bot}})$ in $\Mod\Lambda$ cogenerated by
$\Gproj\Lambda$. Clearly then we have inclusions $\A \subseteq
\colim\Gproj\Lambda \subseteq \GProj\Lambda \subseteq
{^{\pmb{\bot}}}\Lambda$. If $M$ lies in ${^{\pmb{\bot}}}\Lambda$,
let $0 \lxr M \lxr Y^{M} \lxr A^{M} \lxr 0$ be a short exact
sequence, where $Y^{M} \in$ $(\Gproj\Lambda)^{\pmb\bot}$ and $A^{M}
\in \A$. Plainly $Y^{M}$ lies in ${^{\pmb{\bot}}}\Lambda$ $\cap$
$(\Gproj\Lambda)^{\pmb\bot}$ $= \mathcal H$, hence $Y^{M}$ is
projective. Then $M$ lies in $\A$ since the latter is resolving. We
infer that $\A =$ ${^{\pmb{\bot}}}\Lambda$ and therefore $\A =
\colim\Gproj\Lambda = \GProj\Lambda = {^{\pmb{\bot}}}\Lambda$. Then
${^{\bot}}\Lambda = \Gproj\Lambda$ and, by $4.6.1$, the equality
$\GProj\Lambda = \colim\Gproj\Lambda$ shows that $\Lambda$ is
virtually Gorenstein.
\end{proof}
\end{prop}

Recall that the Gorenstein Symmetry Conjecture, \textsf{(GSC)} for
short, see \cite{ARS, BR, B:cm},  asserts that an Artin algebra with
finite one-sided self-injective dimension is Gorenstein. This
conjecture is still open. However virtually Gorenstein algebras, in
particular algebras which are derived or stably equivalent to
algebras of finite representation type satisfy (\textsf{GSC}), see
\cite{B:cm}. The next consequence shows that also algebras for which
stable modules are Gorenstein-projective satisfy the conjecture.

\begin{cor} Let $\Lambda$ be an Artin algebra. If $\id\Lambda_{\Lambda} < \infty$, then the following are equivalent.
\begin{enumerate}
 \item $\Lambda$ is Gorenstein.
  \item  $\Gproj\Lambda$ is contravariantly finite and
  ${^{\bot}}\Lambda \cap (\Gproj\Lambda)^{\bot}$ is of finite representation type.
\item Any module in ${^{\bbot}}\Lambda \cap$ $(\Gproj\Lambda)^{\bbot}$ is a
 direct sum of finitely generated modules.
\end{enumerate}
\begin{proof} If (i) holds then $\Lambda$ is virtually Gorenstein and ${^{\pmb{\bot}}}\Lambda = \GProj\Lambda$.
Then assertion (iii) holds by Proposition $5.12$.   If (iii) holds,
then, by Proposition $5.12$, $\Lambda$ is virtually Gorenstein and
then (i) holds, by Example $4.5.1$.  If $\Lambda$ is Gorenstein,
then by $4.6.5$, $\Gproj\Lambda$ is contravariantly finite and
clearly we have ${^{\bot}}\Lambda \cap (\Gproj\Lambda)^{\bot} =
\proj\Lambda$ which is of finite representation type, so (i) implies
(ii).  If (ii) holds, then by Theorem $5.2$ we have
${^{\bot}}\Lambda = \Gproj\Lambda$. Since $\id\Lambda_{\Lambda} = d
<\infty$, it follows that for any module $A \in \smod\Lambda$, we
have $\Ext^{d+k}_{\Lambda}(A,\Lambda) \cong
\Ext^{k}_{\Lambda}(\Omega^{d}A,\Lambda) = 0$, $\forall k \geq 1$.
Hence $\Omega^{d}A \in {^{\bot}}\Lambda = \Gproj\Lambda$, $\forall A
\in \smod\Lambda$. This implies that
$\rd_{\Gproj\Lambda}\smod\Lambda  = \Gor\bdim \Lambda \leq d$.
Therefore $\Lambda$ is Gorenstein.
\end{proof}
\end{cor}

\subsection{Local Rings} Let throughout $R$ be a commutative Noetherian local ring.
 The following result shows that if $R$ is Gorenstein but not Artinian,
 then the condition ${^{\pmb{\bot}}}R = \GProj R$ does not necessarily holds.

\begin{prop} If $R$ is Gorenstein, then the  following statements are equivalent.
\begin{enumerate}
\item $R$ is complete and  ${^{\pmb{\bot}}}R = \GProj R$.
\item \ $R$ is Artinian.
\end{enumerate}
\begin{proof} Since $R$ is Gorenstein, there exists a
cotorsion pair $(\GProj R,\Proj^{<\infty}R)$ in $\Mod R$, where
$\Proj^{<\infty}R$ denotes the full subcategory of all $R$-modules
with finite projective dimension, see \cite{B:gorenstein, EJ}. If
$A$ lies in ${^{\bot}}\Proj R$, let $(\dag): \ 0 \lxr Y_{A} \lxr
X_{A} \lxr A \lxr 0$ be an exact sequence, where $X_{A}$ is
Gorenstein-projective and $Y_{A}$ has finite projective dimension.
Then clearly $Y_{A}$ lies in ${^{\bot}}\Proj R$, and since $Y_{A}$
has finite projective dimension, we infer that $Y_{A}$ is
projective. Then $(\dag)$ splits and therefore $A$ is
Gorenstein-projective, i.e. $\GProj R = {^{\bot}}\Proj R$.

(i) $\Rightarrow$ (ii) Since $R$ is complete, it
follows that $R$ is pure-injective as an $R$-module, see \cite{JL}.
On the other hand it is well-known that if $E$ is a pure-injective
module over a ring $S$, then for any filtered system $\{A_{i}\, | \,
i \in I\}$ of $S$-modules, we have an isomorphism
$\plim\Ext^{n}_{S}(A_{i},E) \cong \Ext^{n}_{S}(\colim A_{i},E)$,
$\forall n \geq 0$, see \cite[Chap. 1, Proposition 10.1]{Auslander:functors}. We infer that
${^{\pmb{\bot}}}R$, and therefore $\GProj R$, is closed under
filtered colimits. Since $\Gproj R \subseteq \GProj R$, we have
$\colim\Gproj R \subseteq \GProj R$. Since $R$ is Gorenstein, by \cite{EJ} we also have
$\GFlat R \subseteq \colim\Gproj R$. Hence any Gorenstein-flat
$R$-module is Gorenstein-projective. In particular any flat
$R$-module is Gorenstein-projective. However since $R$ is local, it
follows that the big finitistic dimension, which is equal to the
Krull dimension of $R$, is finite. By a result of Jensen
\cite{Jensen} this implies that any flat $R$-module has finite
projective dimension. Since clearly any Gorenstein-projective module
with finite projective dimension is projective, we infer that any
flat $R$-module is projective and therefore $R$ is perfect. This is
equivalent to $R$ being Artinian.

(ii) $\Rightarrow$ (i) Since $R$
is Artinian, it follows that $R$ is complete and moreover  we have ${^{\bbot}}R$ $= {^{\bot}}\Proj R$, since in
this case any projective $R$-module is a direct summand of a product
of copies of $R$, see \cite{KS}. Hence $\GProj R =$ ${^{\bbot}}R$.
\end{proof}
\end{prop}

\begin{cor} The  following statements are equivalent.
\begin{enumerate}
\item $R$ is complete and $\Omega^{\infty}(\Mod R) \subseteq {^{\pmb{\bot}}}R \subseteq \GProj R$.
\item $R$ is Artinian Gorenstein.
\end{enumerate}
If $\mathrm{(ii)}$ holds, then the inclusions in $\mathrm{(i)}$ are equalities.  If $\Gproj R \neq \proj R$, then the above are also equivalent to:
\begin{enumerate}
\item[$\mathrm{(iii)}$] $R$ is complete, $\Gproj R$ is contravariantly finite and $\GProj R = {^{\pmb{\bot}}}R$.
\end{enumerate}
\begin{proof} (i) $\Rightarrow$ (ii) Since $\GProj R \subseteq
{^{\bot}}\Proj R$, the hypothesis implies that $\Omega^{\infty}(\Mod
R) \subseteq {^{\bot}}\Proj R$. Let $P^{\bullet}$ be an acyclic
complex of projective $R$-modules. Then $\forall n \in \mathbb Z$,
the $R$-module $A^{n} := \Ker(P^{n} \xr{} P^{n+1})$ lies in
$\Omega^{\infty}(\Mod R)$. Since the latter is contained in
${^{\bot}}\Proj R$, it follows that the complex
$\Hom_{R}(P^{\bullet},P)$ is acyclic, $\forall P \in \Proj R$. Hence
any acyclic complex of projective $R$-modules is totally acyclic.
Since $R$ is complete, it is well-known that $R$ admits a dualizing
complex. Then by \cite{IK} it follows that $R$ is Gorenstein. Then
$R$ is Artinian by Proposition $5.14$.

(ii) $\Rightarrow$ (i) Clearly $R$ is complete, and
${^{\pmb{\bot}}}R = \GProj R$ by Proposition $5.14$. If $d = \id R$,
then any $R$-module admits a finite resolution of length at most $d$
by Gorenstein-projective modules, equivalently $\Omega^{d}(A)$ is
Gorenstein-projective, $\forall A\in \Mod R$. Let $A$ be in
$\Omega^{\infty}(\Mod R)$ and let $0 \lxr A \lxr P^{0} \lxr P^{1}
\lxr \cdots$ be an exact sequence, where the $P^{n}$s are
projective. Setting $A^{0} = A$ and $A^{n} := \Ker(P^{n}\lxr
P^{n+1})$, $\forall n \geq 1$, it follows that $A = \Omega^{d}A^{d}$
and therefore $A$ is Gorenstein-projective. Hence
$\Omega^{\infty}(\Mod R) \subseteq \GProj R$.
\end{proof}
\end{cor}

\begin{cor} If $R$ is of finite CM-type
and $\Gproj R \neq \proj R$, then the following are equivalent.
\begin{enumerate}
\item $R$ is Artinian (Gorenstein).
\item $R$ is complete and any module in ${^{\pmb{\bot}}}R \cap (\GProj R)^{\bot}$ is a
direct sum of finitely generated modules.
\end{enumerate}
\begin{proof} By Theorem $4.20$,  $R$ is Gorenstein and therefore $(\GProj R)^{\bot} = \Proj^{<\infty}R$.

(i) $\Rightarrow$ (ii)  Clearly $R$ is complete and by Proposition $5.14$ we have ${^{\pmb{\bot}}}R  = \GProj R$. Hence ${^{\pmb{\bot}}}R \cap (\GProj R)^{\bot} = \GProj R \cap \Proj^{<\infty}R = \Proj R$ and (ii) follows.

(ii) $\Rightarrow$ (i) Since $R$ is of finite CM-type, we have $\GProj R = \Add\Gproj R$. This clearly implies that ${^{\bot}}R =
\Gproj R$  and $(\GProj R)^{\bot} = (\Gproj R)^{\pmb{\bot}}$. Since
$R$ is complete, this in turn implies that the resolving subcategory
${^{\pmb{\bot}}}R \cap (\GProj R)^{\bot}$ is definable, and in
particular it is closed under products. Since any
module in ${^{\pmb{\bot}}}R \cap (\GProj R)^{\bot}$ is a direct sum
of modules in ${^{\bot}}R \cap (\Gproj R)^{\bot} = \proj R$, it follows that
${^{\pmb{\bot}}}R \cap (\GProj R)^{\bot} = \Proj R$ and therefore
$\Proj R$ is closed under products. Then by Chase's Theorem, see
\cite{JL}, $R$ is perfect and consequently Artinian.
\end{proof}
\end{cor}

\section{Relative Auslander Algebras} Let $\Lambda$ be an Artin
algebra. If $\Lambda$ is of finite representation type and $\Gamma$
is the endomorphism ring of a representation generator of
$\smod\Lambda$, then  by a classical result of Auslander
\cite{Auslander:queen}, $\Gamma$ is an {\em Auslander algebra}, i.e.
$\gd\Gamma \leq 2 \leq \ddom\Gamma$, where $\ddom\Gamma$ denotes the
dominant dimension of $\Gamma$, that is, the smallest integer $d\geq
0$ such that the first $d$ terms of a minimal injective resolution
of $\Gamma$ are projective. In this section we are concerned with
the global dimension of the endomorphism ring of a representation
generator of a resolving subcategory $\X$ of $\smod\Lambda$. Then we
apply our results to the case $\X = \Gproj\Lambda$ for an Artin
algebra $\Lambda$ of finite CM-type.

\subsection{Auslander Algebras of resolving subcategories of finite representation type}
To proceed further we need the following general observation.

\begin{lem} Let $\A$ be an abelian category with enough projectives
 and $(\X,\Y)$ a cotorsion pair in $\A$ with heart $\omega = \X \cap \Y$.
 If $\mathsf{H} : \A  \lxr \smod\X$, $\mathsf{H}(A) =
\A(-,A)|_{\X}$ is the restricted Yoneda functor, then:
\begin{enumerate}
\item $\mathsf{H} : \A \lxr \smod\X$ is fully faithful and admits an exact left
adjoint $\mathsf{T} : \smod\X \lxr \A$.
\item $\forall A \in
\A$:\, $\pd\mathsf{H}(A) = {\rd}_{\X}A$.
\end{enumerate}
\begin{proof} (i) Define a functor $\mathsf{T} : \smod\X \lxr \A$ as
follows. If $\mathsf{H}(X_{1}) \lxr \mathsf{H}(X_{0}) \lxr F \lxr 0$
is a finite presentation of $F \in \smod\X$, then set $\mathsf{T}(F)
= \Coker(X_{1} \to X_{0})$. It is easy to see that this defines a
functor $\mathsf{T} : \smod\X \lxr \A$ which is left adjoint to
$\mathsf{H}$. For any  object $A\in \A$ let $0 \lxr Y^{0}_{A} \lxr
X^{0}_{A} \lxr A \lxr 0$ and $0 \lxr Y^{1}_{A} \lxr X^{1}_{A} \lxr
Y^{0}_{A} \lxr 0$ be right
 $\X$-approximation sequences of
 $A$ and $Y^{0}_{A}$ respectively, and observe that
 $X^{1}_{A} := \omega^{1}_{A}$ lies in $\X \cap \Y = \omega$ since $\Y$ is closed under extensions.
  Continuing in this way we have an exact resolution
\begin{equation}
\cdots \lxr \omega^{n}_{A} \lxr \omega^{n-1}_{A} \lxr \cdots \lxr
\omega^{1}_{A} \lxr X^{0}_{A} \lxr A \lxr 0
\end{equation}
where the $\omega^{i}_{A}$ lie in $\omega$ and the $Y^{i}_{A} =
\Image(\omega^{i}_{A} \to \omega^{i-1}_{A})$ lie in $\Y$. By
construction this sequence remains exact after applying $\mathsf{H}$
and therefore
\begin{equation}
 \cdots \lxr \mathsf{H}(\omega^{n}_{A}) \lxr \mathsf{H}(\omega^{n-1}_{A}) \lxr
\cdots \lxr \mathsf{H}(\omega^{1}_{A}) \lxr \mathsf{H}(X^{0}_{A})
\lxr \mathsf{H}(A) \lxr 0
\end{equation}
is a projective resolution of $\mathsf{H}(A)$.  Applying
$\mathsf{T}$ and using that $\mathsf{T}\mathsf{H}|_{\X} \cong
\mathsf{Id}_{\X}$, we infer that the counit $\mathsf{T}\mathsf{H}
\lxr \mathsf{Id}_{\A}$ is invertible and the resulted complex is
acyclic. Hence $\mathsf{H}$ is fully
faithful and it is easy to see that $\mathsf{T}$ is exact.

(ii) Set $d := \rd_{\X}A$. If $d = 0$, then $A\in \X$ and therefore
 $\pd\mathsf{H}(A) = 0$, since $\mathsf{H}(A)$ is projective. Assume
 that $d \geq 1$. Since $\X$ is resolving and
 $\rd_{\X}A = d$, the object $Y^{d}_{A} \in
\Y$ in (6.1) lies in $\X$, see \cite[Lemma 3.12]{ABr}. Hence
$Y^{d}_{A} := \omega^{d}_{A}$ lies in $\omega$ and therefore (6.2)
gives a projective resolution
\[
0 \lxr \mathsf{H}(\omega^{d}_{A}) \lxr \mathsf{H}(\omega^{d-1}_{A})
\lxr \cdots \lxr \mathsf{H}(\omega^{1}_{A}) \lxr
\mathsf{H}(X^{0}_{A}) \lxr \mathsf{H}(A) \lxr 0
\]
We infer that $\pd\mathsf{H}(A) \leq d$. If $\pd\mathsf{H}(A) = t <
d$, then $\mathsf{H}(Y^{t}_{A})$ will be projective, or equivalently
$Y^{t}_{A} \in \X \cap \Y = \omega$. Then $\rd_{\X}A \leq t$ and
this is not true. We infer that $\pd\mathsf{H}(A) = d = \rd_{\X}A$.
This argument also shows that $\rd_{\X}A = \infty$ if and only if
$\pd\mathsf{H}(A) = \infty$.
\end{proof}
 \end{lem}

Now let $\X$ be a contravariantly finite
 subcategory of $\smod\Lambda$. As before we denote by $\mathcal L(\X) = \Mod\X$
 the representation category of $\colim\X$. In this subsection we denote by
 $\mathsf{H}$  the functor  $\Mod\Lambda \lxr \mathcal L(\X)$,
  which is given
 by  $\mathsf{H}(A) = \Hom_{\Lambda}(-,A)|_{\X}$.  Contravariant
 finiteness of $\X$  implies that the category $\smod\X$ of coherent functors
 $\X^{\op} \to \ab$  is abelian and then $\mathsf{H}$ induces a functor
 $\mathsf{H}  \colon \smod\Lambda \lxr \smod\X$.

\begin{cor} Let $\X$ be a contravariantly finite
resolving subcategory of $\smod\Lambda$. Then:
\begin{enumerate}
\item $\forall A \in \smod\Lambda$:\, $\pd\mathsf{H}(A) =
{\rd}_{\X}A$.
\item  $\forall A \in \Mod\Lambda$:\, $\fd\mathsf{H}(A) =
{\rd}_{\colim\X}A$.
\item $\forall A \in \Mod\Lambda$:\, $\pd\mathsf{H}(A) =
\rd_{\colim\X}A$ iff $\X$ is of finite representation type.
\end{enumerate}
\begin{proof} By a result of Auslander-Reiten, see
\cite{AR:applications}, there exists a cotorsion pair $(\X,\Y)$ in
$\smod\Lambda$. Then (i) follows from Lemma $6.1$.  By a result of
Krause-Solberg, see \cite{KS:appl}, $(\X,\Y)$ extends to a cotorsion
pair $(\colim\X,\colim\Y)$ in $\Mod\Lambda$. Since $\X$ consists of
finitely presented modules, the functor $\mathsf{H} : \Mod\Lambda
\lxr \Mod\X$ commutes with filtered colimits. This implies that
$\mathsf{H}$ induces an equivalence between $\colim\X$ and
$\Flat(\Mod\X)$ and an equivalence between $\Add\X$ and
$\Proj(\Mod\X)$.  Then working as in Lemma $6.1$, part (ii) follows.
Finally (iii) holds iff any flat $\X$-module is projective iff
$\Mod\X$ is perfect. As in the proof of Theorem $3.1$ this is
equivalent to saying that $\X$ is of finite representation type.
\end{proof}
\end{cor}

\begin{defn} Let $\X$ be a resolving
subcategory of $\smod\Lambda$ which is of finite representation
type. If $T$ is a representation generator of $\X$, i.e. $\X =
\add(T)$,  then the Artin algebra
\[\mathsf{A}(\X) := {\End}_{\Lambda}(T)^{\op}\]
 is called the {\bf
Auslander algebra} of $\X$, or the $\X$-Auslander algebra of
$\Lambda$.
\end{defn}

Clearly the Auslander algebra of $\X$ is unique up to Morita
equivalence, that is $\mathsf{A}(\X)$  is up to Morita equivalence
independent of the choice of the representation generator of $\X$,
and in case $\X = \smod\Lambda$, i.e. $\Lambda$ is of finite
representation type,  $\mathsf{A}(\smod\Lambda)$ is the Auslander
 algebra of $\Lambda$. Auslander proved that
$\gd\mathsf{A}(\smod\Lambda) = 2$ provided that $\Lambda$ is not
semisimple. The following result treats the general case.

\begin{prop} Let $\X$ be a representation-finite resolving subcategory of
$\smod\Lambda$.
\begin{enumerate}
\item $\gd\mathsf{A}(\X) = 0$ if and only if $\Lambda$ is semisimple  (and then $\X = \proj\Lambda$).
\item $\gd\mathsf{A}(\X) = 1$ if and only if $\gd\Lambda = 1$ and $\X = \proj\Lambda$.
\item $\gd\mathsf{A}(\X) = 2$ if and only if  either
\begin{enumerate}
\item $\X = \proj\Lambda$ and $\gd\Lambda = 2$, or
\item $\X \neq \proj\Lambda$ and ${\rd}_{\X}\smod\Lambda \leq 2$.
\end{enumerate}
\item If  ${\rd}_{\X}\smod\Lambda \geq 3$, then: $\gd\mathsf{A}(\X) =
{\rd}_{\X}\smod\Lambda$.
\end{enumerate}
In any case we have the following bounds:
\[
\gd\mathsf{A}(\X) =
\begin{cases} \, \gd\Lambda \,\,\, & \,\,\, \text{if \,\, $\X = \proj\Lambda$,}
\\
\, \bmax\big\{2,\rd_{\X}\smod\Lambda\big\} \,\,\, & \,\,\, \text{if
\,\, $\X\neq \proj\Lambda$.}
\end{cases}
\]
\begin{proof}
(i) If $\Lambda$ is semisimple, then $\X = \proj\Lambda =
\smod\Lambda$ and clearly $\gd\mathsf{A}(\X) = 0$. Conversely if
this holds, then let $A\in \smod\Lambda$ and let $f \colon P \lxr A$
be the projective cover of $A$. Then the projection $\mathsf{H}(P)
\to \Image\mathsf{H}(f)$ splits and therefore $\Image\mathsf{H}(f)$
is projective, hence $\Image\mathsf{H}(f) \cong \mathsf{H}(X)$ for
some direct summand $X$ of $P$. Applying  $\mathsf{T}$ we have $A
\cong X$. Hence $A$ is projective and therefore $\Lambda$ is
semisimple.

(ii) If $\X = \proj\Lambda$ and $\gd\Lambda = 1$, then
$\gd\mathsf{A}(\X) = 1$, since $\mathsf{A}(\X)$ is Morita equivalent
to $\Lambda$. Conversely assume that $\gd\mathsf{A}(\X) = 1$ and
there is a non-projective module $X\in \X$. Let $0 \lxr \Omega X
\lxr  P \lxr X \lxr 0$ be exact where $P$ is projective. Then we
have an exact sequence $0 \lxr \mathsf{H}(\Omega X) \lxr
\mathsf{H}(P) \lxr \mathsf{H}(X) \lxr M \lxr 0$ in
$\smod\mathsf{A}(\X)$ and $M$ is not zero since otherwise $X$ will
be projective as a direct summand of $P$. Let $N :=
\Image(\mathsf{H}(P) \to \mathsf{H}(X))$. Then $N$ is projective,
hence $N = \mathsf{H}(X_{1})$ for some module $X_{1}$ in $\X$ which
necessarily is a direct summand of $P$. Applying $\mathsf{T}$ we
then have $\mathsf{T}(M) = 0$ and $X \cong X_{1}$. This means that
$X$ is projective and this is not the case. We conclude that $\X =
\proj\Lambda$. Hence $\mathsf{A}(\X)$ is Morita equivalent to
$\Lambda$ and then $\gd\mathsf{A}(\X) = \gd\Lambda = 1$.

(iii) If (a) holds, then clearly $\gd\mathsf{A}(\X) = 2$ since
$\mathsf{A}(\X)$ is Morita equivalent to  $\Lambda$. If (b) holds,
and $\rd_{\X} \smod\Lambda = 0$, then $\X = \smod\Lambda$ and
therefore $\Lambda$ is a non-semisimple algebra of finite
representation type. Then  by Auslander's Theorem
\cite{Auslander:queen} we have $\gd\mathsf{A}(\X) = 2$. Now assume
that $1 \leq \rd_{\X}\smod\Lambda \leq 2$. Since there is a
non-projective module $X\in \X$, we have a non-split exact sequence
$0 \lxr \Omega X \lxr P \lxr X \lxr 0$ where $P$ is projective. Then
$\Omega X$ lies in $\X$ since $\X$ is resolving, and  we have a
projective resolution $0 \lxr \mathsf{H}(\Omega X) \lxr
\mathsf{H}(P) \lxr \mathsf{H}(X) \lxr M \lxr 0$ in
$\smod\mathsf{A}(\X)$. As in the proof of (ii) above, $M$ is
non-zero and $N$ is not projective. Hence $\pd M = 2$. On the other
hand any finitely generated $\mathsf{A}(\X)$-module $F$ admits a
presentation $0 \lxr \mathsf{H}(A) \lxr \mathsf{H}(X^{1}) \lxr
\mathsf{H}(X^{0}) \lxr F \lxr 0$ where the $X^{i}$ lie in $\X$ and
$A = \Ker(X^{1} \to X^{0})$.  If $B = \Coker(X^{1} \to X^{0})$, then
we have an exact sequence $0 \lxr A \lxr X^{1} \lxr X^{0} \lxr B
\lxr 0$ in $\smod\Lambda$. If $B = 0$, then $A \in\X$ since $\X$ is
resolving. If $B\neq 0$, then $A\in\X$ since $\rd_{\X}B \leq 2$. In
any case the above exact sequence is a projective resolution of $F$
and therefore $\pd F \leq 2$. Hence $\gd \mathsf{A}(\X) = 2$. If
$\gd\mathsf{A}(\X) = 2$, then for any $\Lambda$-module $A$ we have
$\rd_{\X}A = \pd\mathsf{H}(A) \leq 2$, so $\rd_{\X}\smod\Lambda \leq
2$. If $\X = \proj\Lambda$, then (i), (ii) imply that $\gd\Lambda =
2$.

(iv) Let $d = \rd_{\X}\smod\Lambda \geq 3$.  Let $M$ be a finitely
generated $\mathsf{A}(\X)$-module and let $\mathsf{H}(X^{d-1}) \lxr
\mathsf{H}(X^{d-2}) \lxr \cdots \lxr \mathsf{H}(X^{0}) \lxr M \lxr
0$ be part of the projective resolution of $M$. Then we have an
exact complex $X^{d-1} \lxr X^{d-2} \lxr \cdots  \lxr X^{0} \lxr A
\lxr 0$ in $\smod\Lambda$ where the $X^{i}$ lie in $\X$ and $A =
\Coker(X^{1}\to X^{0})$. Since $\rd_{\X}A \leq d$ and $\X$ is
resolving, it follows that $B := \Ker(X^{d-1}\to X^{d-2})$ lies in
$\X$, see \cite{ABr}. Setting $X^{d} = B$, we have a finite
projective resolution $0 \lxr \mathsf{H}(X^{d}) \lxr
\mathsf{H}(X^{d-1}) \lxr \mathsf{H}(X^{d-2}) \lxr \cdots \lxr
 \mathsf{H}(X^{0}) \lxr M \lxr 0$ of $M$ and
therefore $\gd\mathsf{A}(\X) \leq d$. Then in fact
$\gd\mathsf{A}(\X) = d$ since otherwise $\rd_{\X}\smod\Lambda < d$.
Since clearly  $\rd_{\X}\smod\Lambda = \infty$ if and only if
$\gd\mathsf{A}(\X) = \infty$, we infer that $\rd_{\X}\smod\Lambda =
\gd\mathsf{A}(\X)$.
\end{proof}
\end{prop}

\begin{rem} Let $\Lambda$ be an Artin algebra of finite
representation type.  Also let $T$ be a finitely generated cotilting
module of injective dimension $n \geq 0$, see
\cite{AR:applications}. Then we have a cotorsion pair $(\X,\Y)$ in
$\smod\Lambda$, where $\X = {^{\bot}}T$ and clearly $\X$ is a
representation-finite resolving subcategory of $\smod\Lambda$. It is
well-known that $\rd_{\X}\smod\Lambda = n$. Hence by Proposition
$6.4$ we have $\gd\mathsf{A}(\X) = n$, if $n \geq 3$; if $T$ is not
projective and $n \leq 2$, then $\gd\mathsf{A}(\X) = 2$; finally
$\gd\mathsf{A}(\X) = n$, if $T$ is projective and $n \leq 1$. We
infer that the global dimension of the Auslander algebra of a
(non-trivial) resolving subcategory of finite representation type
can take any value $n \geq 0$. On the other hand let $\X$ be a
representation-finite thick subcategory of $\smod\Lambda$ containing
$\Lambda$, i.e. $\X$ is resolving and closed under cokernels of
monomorphisms, e.g. $\X = (\Gproj\Lambda)^{\bot}$. Then $\gd
\mathsf{A}(\X) = \infty$, except if $\X = \smod\Lambda$ in which
case $\gd\mathsf{A}(\X) \leq 2$.
\end{rem}

Let $\X$ be a resolving contravariantly finite subcategory of
$\smod\Lambda$. We set $\inj\X := \X \cap \inj\Lambda$. We say that
$\X$ has  $\X$-dominant dimension at least two: $\dd_{\inj\X}\X \geq
2$, if for any $X \in \X$, there exists an exact sequence $0 \lxr X
\lxr I^{0} \lxr I^{1}$ where the $I^{i}$ lie in $\inj\X$. If $\X =
\proj\Lambda$, then $\dd_{\proj\Lambda\cap\inj\Lambda}\proj\Lambda
:= \dd\Lambda$ is the dominant dimension of $\Lambda$.

 The following result describes when the $\X$-Auslander algebra of
$\Lambda$ is an Auslander algebra.

\begin{prop} {\bf (1)} Let $\Lambda$ be an Artin algebra and $\X$ a
representation-finite resolving subcategory of $\smod\Lambda$. Then
the following are equivalent.
\begin{enumerate}
\item The $\X$-Auslander algebra of $\Lambda$ is an Auslander
algebra.
\item $\rd_{\X}\smod\Lambda \leq 2 \leq \dd_{\inj\X}\X$.
\end{enumerate}
If $\mathrm{(i)}$ holds, then the Artin algebra  $\Gamma =
\End_{\Lambda}(J)^{\op}$, where $\add J = \inj\X$, is of finite
representation type and there are equivalences $\X \approx
\smod\Gamma$ and $\smod\mathsf{A}(\X) \approx
\mathsf{A}(\smod\Gamma)$. In  particular $\X$ is abelian.

{\bf (2)} Let $\Gamma$ be an Artin algebra of finite representation
type. Then $\X = \proj\mathsf{A}(\smod\Gamma)$ is a resolving
representation-finite subcategory of $\smod\mathsf{A}(\smod\Gamma)$
with $\rd_{\X}\smod\mathsf{A}(\smod\Gamma) \leq 2 \leq
\dd_{\inj\X}\X$.
\begin{proof} {\bf (1)} If (ii) holds, then by Proposition $6.4$, the
$\X$-Auslander algebra $\mathsf{A}(\X)$ of $\Lambda$ has
$\gd\mathsf{A}(\X) \leq 2$. Since $\dd_{\Inj\X}\X \geq 2$, $\forall
X\in\X$, there exists an exact sequence $0 \lxr X \lxr I^{0} \lxr
I^{1}$, where the $I^i$ lie in $\inj\X$. Since $\mathsf{H}$
preserves injectives as a right adjoint of the exact functor
$\mathsf{T}$, we infer that $0 \lxr \mathsf{H}(X) \lxr
\mathsf{H}(I^{0}) \lxr \mathsf{H}(I^{1})$ is exact in
$\smod\mathsf{A}(\X)$, where the $\mathsf{H}(I^{i})$ are
projective-injective. Hence $\dd\mathsf{A}(\X)\geq 2$  and therefore
$\mathsf{A}(\X)$ is an Auslander algebra. It is not difficult to see
that the functor $\mathsf{H} : \smod\Lambda \lxr
\smod\mathsf{A}(\X)$ induces an equivalence $\mathsf{H} : \inj\X
\approx \proj\mathsf{A}(\X) \cap \inj\mathsf{A}(\X)$.  If $J$ is
$\Lambda$-module such that $\add J = \inj\X$, then
$\add\mathsf{H}(J) = \mathsf{H}(\inj\X)$. By Auslander's results
\cite{Auslander:queen}, the ring $\Gamma := \End_{\Lambda}(J)^{\op}
\cong \End_{\mathsf{A}(\X)}\mathsf{H}(J)^{\op}$ is of finite
representation type and $\smod\Gamma$ is equivalent to $\X$. Clearly
then $\smod\mathsf{A}(\X)$ is equivalent to
$\mathsf{A}(\smod\Gamma)$. Conversely if (i) holds, then
$\gd\mathsf{A}(\X) \leq 2$. By Proposition $6.4$ we have
$\rd_{\X}\smod\Lambda \leq 2$. Since $\dd\mathsf{A}(\X) \geq 2$, and
since $\mathsf{H}$ is fully faithful and induces an equivalence
between $\inj\X$ and $\proj\mathsf{A}(\X) \cap \inj\mathsf{A}(\X)$,
it follows that $\dd_{\inj\X}\X \geq 2$.

{\bf (2)}  Let $\Gamma$ be an algebra of finite representation type
and let $\Lambda = \mathsf{A}(\smod\Gamma)$ be its Auslander
algebra. Then $\X = \smod\Gamma \approx \proj\Lambda$ is a resolving
subcategory of finite representation type of $\smod\Lambda$,
$\inj\X$ is the full subcategory of projective-injective modules of
$\smod\Lambda$ and $\rd_{\X}\smod\Lambda = \gd\Lambda \leq 2$ and
$\dd_{\Inj\X}(\X) = \dd\Lambda \geq 2$, since $\Lambda$ is an
Auslander algebra.
\end{proof}
\end{prop}

As a consequence we have the following bound on the dimension
$\bdim{\bf D}^{b}(\smod\Lambda)$ of the bounded derived category of
$\smod\Lambda$ in the sense of Rouquier, see \cite{Rouquier:Kth,
BVdB}. Note that if $\X = \proj\Lambda$, then $\bdim{\bf
D}^{b}(\smod\Lambda) \leq \gd\Lambda = \rd_{\X}\smod\Lambda$, see
\cite{Rouquier:Kth}.

\begin{cor} Let $\X$ representation-finite resolving
subcategory of $\smod\Lambda$.  Then:
\[
\bdim{\bf D}^{b}(\smod\Lambda)\,\ \leq \,\
\bmax\big\{2,{\rd}_{\X}\smod\Lambda\big\}
\]
\begin{proof} The exact left adjoint functor $\mathsf{T} :
\smod\mathsf{A}(\X) \lxr \smod\Lambda$ of $\mathsf{H}$ induces a
triangulated functor ${\bf D}^{b}\mathsf{T} : {\bf
D}^{b}(\smod\mathsf{A}(\X)) \lxr {\bf D}^{b}(\smod\Lambda)$ which
admits as a right adjoint the fully faithful right derived functor of
$\mathsf{H}$. Hence ${\bf D}^{b}(\smod\Lambda)$ is a Verdier
quotient of ${\bf D}^{b}(\smod\mathsf{A}(\X))$ and therefore by
\cite{Rouquier:Kth} we have $\bdim{\bf D}^{b}(\smod\Lambda) \leq
\bdim{\bf D}^{b}(\smod\mathsf{A}(\X))$. Since $\bdim{\bf
D}^{b}(\smod\Gamma) \leq \gd\Gamma$, for any Artin algebra $\Gamma$,
we have  $\bdim{\bf D}^{b}(\smod\mathsf{A}(\X)) \leq
\gd\mathsf{A}(\X)$ and the assertion follows from Proposition $6.4$.
\end{proof}
\end{cor}

\subsection{Cohen-Macaulay Auslander algebras} Let $\Lambda$ be an Artin algebra.
If $\Lambda$ is of finite CM-type, then we call
$\mathsf{A}(\Gproj\Lambda)$, resp. $\mathsf{A}(\uGproj\Lambda)$, the
resp. {\em stable}, {\em Cohen-Macaulay Auslander algebra} of
$\Lambda$.

Now as a consequence of Propositions $6.4$ and $6.6$ we have the
following.

\begin{cor} Let $\Lambda$ be an Artin algebra of finite CM-type.
Then we have the following:
\begin{enumerate}
\item $\gd\mathsf{A}(\Gproj\Lambda) = 0$ if and only if $\Lambda$ is semisimple.
\item $\gd\mathsf{A}(\Gproj\Lambda) = 1$ if and only if $\gd\Lambda = 1$.
\item $\gd\mathsf{A}(\Gproj\Lambda) = 2$ if and only if  either
\begin{enumerate}
\item $\Gproj\Lambda = \proj\Lambda$ and $\gd\Lambda = 2$, or
\item $\Gproj\Lambda \neq \proj\Lambda$ and $\Lambda$ is a Gorenstein algebra with $\id\Lambda \leq 2$.
\end{enumerate}
\item If  $\Gor\bdim\smod\Lambda \geq 3$, then:
\[\gd\mathsf{A}(\Gproj\Lambda) = \Gor\bdim\Lambda = \bmax\big\{\id{_{\Lambda}}\Lambda,
\id\Lambda_{\Lambda} \big\}\]
\item $\Lambda$ is Gorenstein iff its Cohen-Macaulay Auslander algebra
$\mathsf{A}(\Gproj\Lambda)$ has finite global dimension.
\item $\mathsf{A}(\Gproj\Lambda)$ is an Auslander algebra if and
only  if $\Gor\bdim\Lambda \leq 2\leq \dd\Lambda$.
\end{enumerate}
\end{cor}

Note that the Artin algebras $\Lambda$ satisfying  $\Gor\bdim\Lambda
= 2 = \dd\Lambda$ have been described in \cite{AS:gor}.

The following consequence gives an upper bound for the dimension of
the bounded derived category of an Artin algebra of finite CM-type.

\begin{cor} Let $\Lambda$ be a virtually Gorenstein Artin algebra of finite CM-type. Then
\[
\bdim{\bf D}^{b}(\smod\Lambda) \, \leq \, \bmax\big\{2,\id\Lambda\big\}
\]
\begin{proof} If $\id\Lambda = \infty$, there is nothing to prove.
Assume that $\id\Lambda_{\Lambda} = d < \infty$. Since $\Lambda$ is
virtually Gorenstein, it follows that $\Lambda$ is Gorenstein
and $\id{_{\Lambda}}\Lambda = \id_{\Lambda}\Lambda = d$. If
$\Gproj\Lambda = \proj\Lambda$, i.e. $\Lambda$ is Morita equivalent
to $\mathsf{A}(\Gproj\Lambda)$, then $\gd\Lambda = d$. Then the
assertion follows from Corollaries $6.7$ and $6.8$.
\end{proof}
\end{cor}

In the commutative case we have the following consequence. Note that
if $R$ is Cohen-Macaulay, then the second inequality below was first
shown by Leuschke \cite{Leuschke}.

\begin{cor} Let $R$ be a commutative Noetherian complete local ring
of finite CM-type. Then:
\[
\bdim{\bf D}^{b}(\smod R) \, \leq \, \gd\mathsf{A}(\Gproj R) \, \leq
\,  \bmax\big\{2,\bdim R\big\}\,  < \, \infty
\]
and the second inequality is an equality if $\bdim R \geq 2$.
\begin{proof} By Theorem $4.20$, $R$ is Gorenstein. Since the injective dimension $\id R$ of $R$ coincides with the
Krull dimension $\bdim R$, the assertions follow as in Corollaries
$6.8$ and $6.9$.
\end{proof}
\end{cor}

\begin{exam} Let $\Lambda$ be a non-semisimple Artin algebra. If the resolving
subcategory  $\Sub(\smod\Lambda)$ of $T_{2}(\Lambda)$ is of finite
representation type, it follows by Proposition $6.4$ that
$\gd\mathsf{A}(\Sub(\smod\Lambda)) = 2$, since $\Sub(\smod\Lambda)$
is closed under submodules in $\smod T_{2}(\Lambda)$.  If in
addition $\Lambda$ is self-injective, e.g. $\Lambda_{n} =
k[t]/(t^{n})$ with $n \leq 5$, then, by Example $4.17$, we have
$\Gor\bdim T_{2}(\Lambda_{n}) = 1$ and therefore
$\gd\mathsf{A}(\Gproj T_{2}(\Lambda_{n})) = 2$.
\end{exam}

The following gives a convenient characterization  of the full subcategory of Gorenstein-projectives.

\begin{lem} Let $\A$ be an abelian category with enough projectives. If $\X$ is a full subcategory of $\GProj\A$, then the following are equivalent:
\begin{enumerate}
\item $\X = \GProj\A$.
\item $\X$ is resolving and $\rd_{\X}\GProj\A < \infty$.
\end{enumerate}
\begin{proof}  If $\X = \GProj\A$, then clearly $\X$ is resolving and $\rd_{\X}\GProj\A = 0$. Conversely    let $\rd_{\X}\GProj\A = d < \infty$ and let $A \in \GProj\A$. Then $A \cong \Image\bigl(P^{-1} \xr{} P^{0}\bigr)$ for some (totally) acyclic complex  $\cdots \lxr P^{-1} \lxr P^{0} \lxr P^{1} \lxr \cdots$ of projectives. Clearly then $B = \Image\bigl(P^{d-1} \xr{} P^{d}\bigr)$ is Gorenstein-projective and $\Omega^{d}B = A$. Since $\rd_{\X}B \leq d$, by \cite[3.12]{ABr} it follows that $A \in \X$.  We infer that $\GProj\A = \X$.
\end{proof}
\end{lem}

Let $\A$ be an abelian category with enough projectives. Following \cite[VII.1]{BR} we say that $\A$ is {\em Gorenstein} if $\rd_{\GProj\A}\A < \infty$. For instance if $\Lambda$ is a Noetherian ring, then $\Lambda$ is Gorenstein if and only if the category  $\Mod\Lambda$ is Gorenstein in the above sense.

Part (ii) of the following consequence gives a trivial proof to (a generalization of) the main result of  \cite{LiZhang}.

\begin{cor} Let $\A$ be an abelian category with enough projectives.
\begin{enumerate}
\item  $\A$ is  Gorenstein, resp. Gorenstein of finite CM-type, if and only if $\GProj\A$ admits a resolving subcategory $\X$, resp.  of finite representation type, such that $\rd_{\X}\A < \infty$.
\item  If $\A$ is $R$-linear over a commutative ring $R$ and the $R$-module $\A(A,B)$ is finitely generated, $\forall A,B \in \A$, then the following are equivalent:
\begin{enumerate}
\item $\A$ is  Gorenstein of finite CM-type.
\item $\GProj\A$ admits a resolving subcategory $\X$ of finite representation type such that:
    \[
    \gd\mathsf{A}(\X) \, < \, \infty
    \]
\end{enumerate}
\end{enumerate}
\begin{proof} Part (i) follows from Lemma $6.12$. For part (ii) note that for any object $T$ of $\A$, the subcategory $\add T$ is contravariantly finite in $\A$. Therefore if $\X$ is of finite representation type, then the category $\X$ has weak kernels or equivalently the category $\mathsf{A}(\X) = \smod\X \approx \smod\End_{\A}(T)^{\op}$  of coherent functors over $\X = \add T$ is abelian.   The rest follows from Lemma $6.12$ using that $\rd_{\X}\A < \infty$ if and only if $\gd\mathsf{A}(\X) < \infty$, see Proposition $6.4$.
\end{proof}
\end{cor}

We close this section with an application of Corollary $6.8$ and
results of Buchweitz \cite{Buchweitz} to algebras of finite CM-type
having Cohen-Macaulay Auslander algebras of global dimension $2$.

Let $\Lambda$ be an finite-dimensional $k$-algebra of finite CM-type
over a field $k$ and assume that $\Lambda$ has infinite global
dimension. Let $T$ be a representation generator of $\Gproj\Lambda$,
so $\mathsf{A}(\Gproj\Lambda)$ is Morita equivalent to $\Gamma :=
\End_{\Lambda}(T)^{\op}$. Then $\underline{T}$ is a representation
generator of $\uGproj\Lambda$ and the stable Cohen-Macaulay
Auslander algebra $\mathsf{A}(\uGproj\Lambda)$ of $\Lambda$ is
self-injective and  Morita equivalent to $\Delta :=
\uEnd_{\Lambda}(T)^{\op}$. By \cite{Buchweitz}, the canonical ring
epimorphism $\Gamma \lxr \Delta$, is pseudo-flat in the sense that
$\mathsf{Tor}^{\Gamma}_{1}(\Delta,\Delta) = 0$.
 For the notion of a quasi-periodic resolution we
refer to \cite{Buchweitz}.

\begin{cor} Let $\Lambda$ be Gorenstein of dimension $2$ and consider the $\Delta$-bimodule $L =
\mathsf{Tor}^{\Gamma}_{2}(\Delta,\Delta)$.
\begin{enumerate}
\item $L \cong
\uHom_{\Lambda}(T,\Omega T)$ is an invertible $\Delta$-bimodule with
inverse 
\[
L^{-1} = {\Ext}^{2}_{\Gamma}(\Delta,\Gamma) \, \cong \,
{\uHom}_{\Lambda}(\Omega T,T)
\]
\item $\Delta$ admits a quasi-periodic projective resolution over
the $k$-enveloping algebra $\Delta\otimes_{k}\Delta^{\op}$ of $\Delta$ of
period $3$ with periodicity factor $\uHom_{\Lambda}(T,\Omega T)$.
\item If  $\uHom_{\Lambda}(T,\Omega T)$ is of finite order $n$ in the Picard group $\mathsf{Pic}\,\Delta$,
then the Hochschild (co)-homology of $\Delta$ is periodic of period
dividing $3n$.
\end{enumerate}
\begin{proof} Since $\Lambda$ is Gorenstein of dimension 2, by
Corollary $6.8$ we have that $\Gamma$ is of global dimension $2$, so
its Hochschild dimension is $2$. Then the assertions follow by
\cite[Theorem 1.5]{Buchweitz}.
\end{proof}
\end{cor}

\section{Rigid Gorenstein-Projective Modules}
Our aim in this section is to study  rigid or cluster-tilting
finitely generated Gorenstein-projective modules $T$ over an Artin
algebra $\Lambda$ in connection with Cohen-Macaulay finiteness of
$\Lambda$ and representa--tion-theoretic properties of the stable
endomorphism algebra of $T$. In particular using recent results of
Keller-Reiten and Amiot we give in some cases convenient
descriptions of the stable triangulated category of Gorenstein-projectives in terms of the cluster category $\C_{Q}$ associated to
the quiver $Q$ of the stable endomorphism algebra of $T$. Moreover
we show that, in this setting, derived equivalent virtually
Gorenstein algebras share the same cluster category.

Recall that if $\C$ is an abelian, resp. triangulated, category,
then an object $X \in \C$ is called {\bf rigid} if
$\Ext^{1}_{\C}(X,X) = 0$, resp. $\C(X,X[1]) = 0$. If $\C$ is
triangulated, then $X$ is called {\bf cluster-tilting}, see
\cite{KR}, \cite{Iyama}, if $\add X$ is functorially finite and 
\[
\big\{C
\in \C\, | \, \C(C,X[1]) = 0\big\} = \add X =  \big\{C \in \C \,  | \,
\C(X,C[1]) = 0 \big\}
\]

Let $\Lambda$ be an Artin algebra. If $T$ is a finitely generated Gorenstein-projective module, then we denote by $(\Gproj\Lambda)^{\leq 1}_{T}$ the full subcategory of $\Gproj\Lambda$ consisting of all modules $A$  admitting an exact sequence $0 \lxr T_{1} \lxr T_{0}\lxr A \lxr 0$, where the $T_{i}$ lie in $\add T$, i.e.
$\rd_{\add T}A \leq 1$.

Recall that the {\bf
representation dimension} $\rp\Lambda$ of $\Lambda$
in the sense of Auslander \cite{Auslander:queen} is defined by
\[
\rp\Lambda =
\inf\big\{\gd{\End}_{\Lambda}\big(\Lambda\oplus\mathsf{D}(\Lambda)\oplus
X\big) \,\, | \,\, X \in \smod\Lambda\big\}
\]
Note that
$\rp\Lambda \leq 2$ iff $\Lambda$ is of finite representation type, see \cite{Auslander:queen}.  The following result gives a connection between Artin algebras of
finite CM-type and representation-theoretic properties of stable
endomorphism algebras of rigid Gorenstein-projective modules.

\begin{thm} Let  $T$ be a finitely generated Gorenstein-projective
module, and assume that $\Lambda \in \add T$ and the module
$\Omega^{n}T$ is rigid,  for some $n \geq 0$.
\begin{enumerate}
\item There is an equivalence
 \[ (\uGproj\Lambda)^{\leq 1}_{T}\big/\add\underline{T}
 \,\ \stackrel{\approx}{\lxr}\,\
 \smod\underline{\End}_{\Lambda}(T)^{\op}
 \]
and the following statements are equivalent.
 \begin{enumerate}
\item $(\Gproj\Lambda)^{\leq 1}_{T}$ is of finite representation type.
\item The stable
endomorphism algebra $\uEnd_{\Lambda}(T)^{\op}$ is of finite representation
type.
 \end{enumerate}
 In particular if $\Lambda$ is of finite CM-type, then $\uEnd_{\Lambda}(T)^{\op}$ is of finite representation type.
\item The inclusion $(\Gproj\Lambda)^{\leq 1}_{T} \subseteq \Gproj\Lambda$ is an equality if and only $\underline{T}$ is a cluster tilting object in $\uGproj\Lambda$. If this is the case, then:
 \begin{enumerate}
\item $\Lambda$ is of finite CM-type if and only if $\uEnd_{\Lambda}(T)^{\op}$
 is of finite representation type.
\item The following statements are equivalent:
\begin{enumerate}
\item[1.] $\U = \{X \in \Gproj\Lambda \, | \, \text{any map} \,\ \Omega T \to X
\,\ \text{factorizes through a module in} \,\ \add\Omega^{-1}T\}$ is
of finite representation type.
\item[2.] $\V = \{X \in \Gproj\Lambda \, | \, \text{any map} \,\ X \to
\Omega^{-3}T \,\ \text{factorizes through a module in} \,\
\add\Omega^{-1}T\}$ is of finite representation type.
\item[3.] The stable endomorphism algebra $\uEnd_{\Lambda}(T)^{\op}$ is of finite CM-type.
\end{enumerate}
If $\mathrm{3.}$ holds, then $\rp\uEnd_{\Lambda}(T)^{\op} \leq 3$.
\end{enumerate}
\end{enumerate}
\begin{proof} (i) Since $T$ is Gorenstein-projective, we have
isomorphisms:
\[
{\Ext}^{1}_{\Lambda}(\Omega^{n}T,\Omega^{n}T)\,\  \cong \,\
{\uHom}_{\Lambda}(\Omega^{n+1}T,\Omega^{n}T)\,\  \cong \,\
{\uHom}_{\Lambda}(\Omega T,T)\,\ \cong \,\
{\uHom}_{\Lambda}(T,\Omega^{-1}T) \,\ \cong \,\
{\Ext}^{1}_{\Lambda}(T,T)
\]
Hence $T$ is rigid if and only if $\Omega^{n}T$ is rigid, $\forall n
\geq 0$, if and only if $\underline{T}$ is rigid in
$\uGproj\Lambda$.  So we may assume that $T$ is rigid. Consider the stable category $(\uGproj\Lambda)^{\leq 1}_{T}$ of $(\Gproj\Lambda)^{\leq 1}_{T}$ modulo projectives, and the functor
\[
\mathsf{H}\,  : \,  (\uGproj\Lambda)^{\leq 1}_{T} \, \lxr \,
\smod{\uEnd}_{\Lambda}(T)^{\op}, \,\,\,\,\, \mathsf{H}(\unA) =
{\uHom}_{\Lambda}(T,A)
\]
which induces an equivalence between $\add \underline{T}$ and
$\proj\uEnd_{\Lambda}(T)^{\op}$. Since $\Lambda \in \add T$, it follows that $\proj\Lambda \subseteq \add T$. It is easy to see that  $(\uGproj\Lambda)^{\leq 1}_{T}$ coincides with the extension category $\add \underline{T} \star\add \Omega^{-1}\underline{T}$ which by definition consists  of all objects $\unA$ for which there exists a triangle $\underline{T}_{0} \lxr \unA \lxr \Omega^{-1}\underline{T}_{1} \lxr \Omega^{-1}\underline{T}_{0}$ in $\uGproj\Lambda$, where the $T_{i}$ lie in $\add T$.  Since $\underline{T}$ is rigid, by a result of Keller-Reiten \cite{KR}, the functor $\mathsf{H}$ induces  an equivalence between the stable category $\bigl(\add \underline{T} \star\add \Omega^{-1}\underline{T}\bigr)\big/\add\underline{T}$ and $\smod\uEnd_{\Lambda}(T)^{\op}$, hence we have an equivalence $\bigl(\uGproj\Lambda)^{\leq 1}_{T} \big/\add\underline{T}\, \approx \,\ \smod\uEnd_{\Lambda}(T)^{\op}$. It follows directly from this that $(\Gproj\Lambda)^{\leq 1}_{T}$ is of finite representation type if and only if $\uEnd_{\Lambda}(T)^{\op}$ is of finite representation type.

(ii) If $\underline{T}$ is a cluster tilting object in $\uGproj\Lambda$, then we have $\uGproj\Lambda = \add \underline{T} \star \add \Omega^{-1}\underline{T}$, see \cite{KR}, and therefore  $\uGproj\Lambda = (\uGproj\Lambda)^{\leq 1}_{T}$. It follows that  $\Gproj\Lambda = (\Gproj\Lambda)^{\leq 1}_{T}$.  Conversely assume that $\Gproj\Lambda = (\Gproj\Lambda)^{\leq 1}_{T}$, or equivalently $\uGproj\Lambda = (\uGproj\Lambda)^{\leq 1}_{T}$. Let $\unA$ be in $\uGproj\Lambda$ such that $\uHom(T,\Omega^{-1}A)= 0$.  Since $\uGproj\Lambda = \add \underline{T} \star \add \Omega^{-1}\underline{T}$, there exists a triangle $\underline{T}_{0} \lxr \Omega^{-1}\unA \lxr \Omega^{-1}\underline{T}_{1} \lxr \Omega^{-1}\underline{T}_{0}$ in $\uGproj\Lambda$, where the $T_{i}$ lie in $\add T$. Then by hypothesis the first map is zero and therefore $\Omega^{-1}\unA$ is a direct summand of $\Omega^{-1}\underline{T}_{1}$, or equivalently $\unA$ is a direct summand of $\underline{T}_{1}$, i.e. $\unA \in \add \underline{T}$. Hence $\add \underline{T} = \{\unA \in \uGproj\Lambda \,\ | \,\ \uHom(\underline{T},\Omega^{-1}\unA) = 0\}$. Since $\add\underline{T}$ is contravariantly finite, by \cite{KZ} it follows that $\underline{T}$ is a cluster tilting object in $\uGproj\Lambda$.

Assume now that $\underline{T}$ is a cluster tilting object in $\uGproj\Lambda$. Then part (a) follows from (i). By
a result of Keller-Reiten \cite{KR}, the endomorphism algebra
$\uEnd_{\Lambda}(T)^{\op}$ is Gorenstein of Gorenstein dimension
$\leq 1$. Clearly in this case we have that
$\Gproj\uEnd_{\Lambda}(T)^{\op}$ consists of the submodules of the
finitely generated projective $\uEnd_{\Lambda}(T)^{\op}$-modules,
and $\Ginj\uEnd_{\Lambda}(T)^{\op}$ consists of the factors of the
finitely generated injective $\uEnd_{\Lambda}(T)^{\op}$-modules.
Using these descriptions it is not difficult to see that the functor
$\mathsf{H}$ induces an equivalence between
$\U/\add\Omega^{-1}\underline{T}$ and
$\Gproj\uEnd_{\Lambda}(T)^{\op}$ and between
$\V/\add\Omega^{-1}\underline{T}$ and
$\Ginj\uEnd_{\Lambda}(T)^{\op}$. Then part (b) is a direct
consequence of these equivalences and the fact that the stable categories $\uGproj\uEnd_{\Lambda}(T)^{\op}$ and $\oGinj\uEnd_{\Lambda}(T)^{\op}$ are equivalent, see Remark $4.3$.  Finally by using a recent result
of Ringel \cite{Ringel:torsionless}, it follows easily that the
representation dimension of a Gorenstein algebra of finite CM-type
and of Gorenstein dimension at most one, is at most $3$. We infer
that $\rp\uEnd_{\Lambda}(T)^{\op}\leq 3$.
 \end{proof}
\end{thm}

We have  the following consequence.

\begin{cor} Let $\Lambda$ be an Artin algebra and assume that $\uGproj\Lambda$ contains a cluster tilting object $\underline{T}$ such that the algebra $\uEnd_{\Lambda}(T)^{\op}$ is of finite CM-type. Then $\Lambda$ is of infinite CM-type if and only if $\rp\uEnd_{\Lambda}(T)^{\op} = 3$.
\begin{proof} Since $\rp\uEnd_{\Lambda}(T)^{\op}\leq 2$ if and only if
$\uEnd_{\Lambda}(T)^{\op}$ is of finite representation type if and
only if $\Lambda$ is of finite CM-type, the assertion follows from part (ii) of Theorem $7.1$.
\end{proof}
\end{cor}

We close this section by giving convenient descriptions of the
triangulated category $\uGproj\Lambda$, for a virtually Gorenstein
$k$-algebra $\Lambda$ of finite CM-type over an algebraically closed
field $k$, under the presence of a cluster tilted object. These
descriptions follow from recent results of Amiot \cite{Amiot} and
Keller-Reiten \cite{KR:acyclic}. First we recall that if $Q$ is a
quiver without oriented cycles, and $d \geq 2$ or $d = 1$ and $Q$ is
a Dynkin quiver, then the $d$-{\bf cluster category} of $Q$ over a
field $k$, is the orbit category $\C^{(d)}_{Q} := {\bf D}^{b}(\smod
kQ)\big/ (\mathbb S^{-1}[d])^{\mathbb Z}$ of the bounded derived
category  ${\bf D}^{b}(\smod
kQ)$ of finite-dimensional modules over the path algebra $kQ$
under the action of the automorphism group generated by $X \mapsto
\mathbb S^{-1}(X[d])$, where $\mathbb S = -\otimes^{\bf L}_{kQ}\mathsf{D}(kQ)$.  As shown by Keller \cite{Keller:orbit},
$\C^{(d)}_{Q}$ is a triangulated category and the projection ${\bf
D}^{b}(\smod kQ) \lxr \C^{(d)}_{Q}$ is a triangle functor. If $d =
2$, we write $\C^{(d)}_{Q} = \C_{Q}$.

If $\T$ is a $k$-linear triangulated category over a field $k$ with
finite-dimensional $\Hom$-spaces, then a {\bf  Serre functor} for
$\T$, in the sense of Bondal-Kapranov \cite{BondalKapranov}, is a
triangulated equivalence $\mathbb S : \T \lxr \T$ such that:
\[
\mathsf{D}{\Hom}_{\T}(A,B) \,\ \stackrel{\approx}{\lxr} \,\
{\Hom}_{\T}(B,\mathbb S A)
\]
naturally $\forall A, B \in \T$. If $\T$ admits a Serre functor
$\mathbb S$, then $\mathbb S$ is unique up to an isomorphism of
triangulated functors, and then $\T$ is called ({\bf weakly})
$d$-{\bf Calabi-Yau} if there exists a natural isomorphism $\mathbb
S(?) \ \stackrel{\approx}{\lxr} \ (?)[d]$ as (additive) triangulated
functors. E.g. the $d$-cluster category $\C^{(d)}_{Q}$ is
$d$-Calabi-Yau \cite{Keller:orbit}. The following describes when for
an Artin algebra $\Lambda$, the triangulated category
$\uGproj\Lambda$ admits a Serre functor.

\begin{lem}
\begin{enumerate}
\item The triangulated category $\uGproj\Lambda$ admits a Serre functor
if and only if  for any $A \in \Gproj\Lambda$, the (minimal) right
$\GProj\Lambda$-approximation $X_{\mathsf{DTr}(A)}$ of
$\mathsf{DTr}A$ is finitely generated.
\item If $\Lambda$ is virtually
Gorenstein, then $\uGproj\Lambda$ admits a Serre functor which is given by
$\mathbb S(\unA) = \Omega^{-1}X_{\mathsf{DTr}A}$.
\item If $\Lambda$ is of finite CM-type and
$\uGproj\Lambda$ is $d$-Calabi-Yau, then the stable Cohen-Macaulay
Auslander algebra $\mathsf{A}(\uGproj\Lambda)$ of $\Gproj\Lambda$ is
self-injective and the stable category
$\umod\mathsf{A}(\uGproj\Lambda)$ is $(3d-1)$-Calabi-Yau.
\end{enumerate}
\begin{proof} (i) If
$\mathbb S$ is a Serre functor on $\uGproj\Lambda$, then
$\mathsf{D}\uHom_{\Lambda}(A,B) \cong \uHom_{\Lambda}(B,\mathbb S
A)$ naturally $\forall A, B \in \Gproj\Lambda$. Since
$\uGProj\Lambda$ is compactly generated, by Brown representability,
we have natural isomorphisms $\mathsf{D}\uHom_{\Lambda}(A,B) \cong
\uHom_{\Lambda}(B,\Omega^{-1}X_{\mathsf{DTr}A})$, see
\cite[Proposition 8.8]{B:cm}. Hence $\forall A \in \Gproj\Lambda$ we
have a natural isomorphism $\mathbb S A \cong
\Omega^{-1}X_{\mathsf{DTr}A}$. It follows that
$\Omega^{-1}X_{\mathsf{DTr}A}$, hence $X_{\mathsf{DTr}A}$, is
finitely generated, $\forall A \in \Gproj\Lambda$. Conversely if
this holds, then clearly  $A \mapsto \Omega^{-1}X_{\mathsf{DTr}A}$
is a Serre functor on $\uGproj\Lambda$.

(ii), (iii) Since  $\Lambda$ is virtually Gorenstein iff $X_{M}$ is
finitely generated for any $M \in \smod\Lambda$, see \cite[Corollary
8.3]{B:cm}, (ii) follows from (i). Part (iii) follows from our
previous results and a result of Keller, see \cite[Lemma
8.5.2]{Keller:orbit}, which asserts that if $\T$ is a $d$-Calabi-Yau
triangulated category, then the stable triangulated category $\umod\T$ of the
Frobenius category $\smod\T$ of coherent functors over $\T$ is
$(3d-1)$-Calabi-Yau.
\end{proof}
\end{lem}

We have the following consequence.

 \begin{cor} Let $\Lambda$ be a virtually Gorenstein finite-dimensional $k$-algebra over a field $k$ and let
 $ T$ be a cluster tilting object in the triangulated category
 $\uGproj\Lambda$. Assume that $\uGproj\Lambda$ is weakly
 $2$-Calabi-Yau:  $X_{\mathsf{DTr}A} \cong \Omega^{-1}A$, naturally
 in $\uGproj\Lambda$, for any finitely generated Gorenstein-projective
 module $A$.
 \begin{enumerate}
 \item If the quiver $Q$ of $\uEnd_{\Lambda}(T)^{\op}$ has no oriented cycles,
 then there is a triangle equivalence
 \[\uGproj\Lambda \,\ \stackrel{\approx}{\lxr} \,\ \C_{Q}\]
 \item If $\Lambda$ is of finite CM-type and $\uGproj\Lambda$ is standard, then there is a triangle equivalence
 \[ \uGproj\Lambda \,\ \stackrel{\approx}{\lxr} \,\  \C_{Q}/G\] where
  $\C_{Q}/G$ is the orbit category of $\C_{Q}$ under the action of some cyclic group $G$ of
 automorphisms.
 \item If $\uGproj\Lambda$ is $2$-Calabi-Yau and the algebra $\uEnd_{\Lambda}(T)^{\op}$ is of finite CM-type, then the stable Cohen-Macaulay Auslander
algebra $\mathsf{A}(\uGproj\uEnd_{\Lambda}(T)^{\op})$ of $\Gproj\uEnd_{\Lambda}(T)^{\op}$ is
self-injective and the stable triangulated category
$\umod\mathsf{A}(\uGproj\uEnd_{\Lambda}(T)^{\op})$ is $8$-Calabi-Yau.
 \end{enumerate}
\begin{proof} (i) Since $\Gamma := \uEnd_{\Lambda}(T)^{\op}$ is Gorenstein with $\id_{\Gamma}\Gamma \leq 1$,
the assumption implies that $\Gamma$ is
hereditary. Then the assertion follows from \cite{KR:acyclic}. Part
(ii) follows easily from Amiot's results \cite{Amiot}. Finally part
(iii) follows from Lemma $7.3$ and the fact that $\uGproj\Lambda$ being $2$-Calabi-Yau implies that $\uGproj\uEnd_{\Lambda}(T)^{\op}$ is $3$-Calabi-Yau, see \cite{KR}.
\end{proof}
 \end{cor}

The following consequence shows that derived equivalent algebras
satisfying the conditions of the above corollary share the same
(orbit category of the) associated cluster category.

\begin{cor} Let $\Lambda$ and $\Gamma$ be derived equivalent algebras. Assume that $\Lambda$ is virtually Gorenstein
and the stable triangulated category $\uGproj\Lambda$ is (weakly)
$2$-Calabi-Yau and admits  a cluster tilting object $T$.
\begin{enumerate}
\item $\uGproj\Gamma$ is (weakly) $2$-Calabi-Yau and admits a cluster tilting object $S$.
\item If the quiver $Q_{\Lambda}$ of $\uEnd_{\Lambda}(T)^{\op}$ has no oriented cycles,
then so does the quiver $Q_{\Gamma}$ of $\uEnd_{\Gamma}(S)^{\op}$
and the associated cluster categories are triangle equivalent:
$\C_{Q_{\Lambda}} \approx \C_{Q_{\Gamma}}$.
\item If $\Lambda$ is of finite CM-type  and $\uGproj\Lambda$ is standard, then $\Gamma$ is of finite CM-type and the
associated orbit categories of the cluster categories of
$Q_{\Lambda}$ and $Q_{\Gamma}$  are triangle equivalent:
$\C_{Q_{\Lambda}}/G \approx \C_{Q_{\Gamma}}/G^{\prime}$.
\end{enumerate}
\begin{proof} It is shown in \cite{B:cm} that if $\Lambda$ and
$\Gamma$ are derived equivalent Artin algebras, then $\Lambda$ is
virtually Gorenstein if and only if so is $\Gamma$. Moreover in this
case there is induced a triangle equivalence between
$\uGproj\Lambda$ and $\uGproj\Gamma$. Then the assertions follow
directly from Corollary $7.4$.
\end{proof}
\end{cor}

\medskip

{\bf Acknowledgement.} The author thanks the referee for the very helpful comments and remarks.


\medskip

\medskip

\begin{thebibliography}{99}

\bibitem{Amiot}
\textsc{C.~Amiot}, \textit{On the structure of triangulated categories
with finitely many indecomposables},  Bull. Soc. Math. France {\bf
135} (2007),  no. 3, 435--474.


\bibitem{AT}
\textsc{L.~Angeleri-H\"{u}gel} and \textsc{J.~Trlifaj}, \textit{Direct limits of modules of
finite projective dimension}, in: Rings, Modules, Algebras, and
Abelian Groups, in: LNPAM, vol. {\bf 236}, M. Dekker, 2004, 27--44.

\bibitem{AST}
\textsc{L.~Angeleri-H\"{u}gel, J.~Saroch and J.~Trlifaj}, \textit{On
the telescope conjecture for module categories},  J. Pure Appl.
Algebra {\bf 212} (2008), no. 2, 297--310.


\bibitem{ABM}
\textsc{I.~Assem, A.~Beligiannis and  N.~Marmaridis}, \textit{Right
Triangulated Categories with Right Semi-equivalence}, in: Algebras
and Modules II (Geiranger 1996), CMS Conf. Proc. {\bf 24} Amer.
Math. Soc. Providence, RI,  (1998), 17--37.

   \bibitem{Auslander:coherent}
\textsc{M.~Auslander},
 \textit{Coherent Functors}, in: Proceedings
of the Conference on Categorical Algebra, La Jolla (1966), 189--231.


\bibitem{Auslander:ENS}
\textsc{M.~Auslander}, {\em Anneaux de Gorenstein, et torsion en alg\`{e}bre
commutative}, S\'{e}minaire d' Alg\`{e}bre Commutative, dirig\'{e} par
Pierre Samuel, 1966/67. Paris (1967), 69 pp.


\bibitem{Auslander:queen}
\textsc{M.~Auslander}, \textit{Representation dimension of Artin algebras},  Queen
Mary College Mathematics Notes, London (1971).

\bibitem{Auslander:II}
\textsc{M.~Auslander}, \textit{Representation theory of Artin algebras. I, II.},
Comm. Algebra  {\bf 1} (1974), 177--268; ibid. 1 (1974), 269--310.

\bibitem{Auslander:large}
\textsc{M.~Auslander}, \textit{Large modules over Artin algebras}, in: Algebra, topology,
and category theory (a collection of papers in honor of Samuel
Eilenberg),  pp. 1--17. Academic Press, New York, 1976.

\bibitem{Auslander:functors}
\textsc{M.~Auslander}, \textit{Functors and morphisms determined by objects},
 Representation theory of algebras (Proc. Conf., Temple Univ., Philadelphia, Pa., 1976),  pp. 1--244.
 Lecture Notes in Pure Appl. Math., Vol. 37, Dekker, New York, 1978.

\bibitem{ABr}
\textsc{M.~Auslander and M.~Bridger}, \textit{Stable module theory}, Mem.
Amer. Math. Soc.,  No. {\bf 94} (1969), 146 pp.

\bibitem{ABu}
\textsc{M.~Auslander and R.-O.~Buchweitz}, \textit{The homological theory of
maximal Cohen-Macaulay approximations}, Colloque en l'honneur de
Pierre Samuel (Orsay, 1987). M\'em. Soc. Math. France (N.S.)  No.
{\bf 38} (1989), 5--37.

\bibitem{AR:dual}
\textsc{M.~Auslander and I.~Reiten}, \textit{Stable Equivalence of
Dualizing $R$-Varieties}, Advances in Math. {\bf 12} (1974),
306--366.



\bibitem{AR:applications}
\textsc{M.~Auslander and I.~Reiten}, \textit{Applications of Contravariantly Finite
subcategories}, Adv. in Math. {\bf 86} (1991), 111--152.

\bibitem{AR:cm}
\textsc{M.~Auslander and I.~Reiten}, \textit{Cohen-Macaulay and Gorenstein
Algebras}, Progr. Math. {\bf 95} (1991), 221--245.

\bibitem{AR:dtr}
\textsc{M.~Auslander and I.~Reiten}, \textit{$D$Tr-periodic modules and functors},
Representation theory of algebras (Cocoyoc, 1994),  39--50, CMS
Conf. Proc., {\bf 18}, Amer. Math. Soc., Providence, RI, (1996).

\bibitem{ARS}
\textsc{M.~Auslander,  I.~Reiten and S.~Smal{\o}}, \textit{Representation
Theory of Artin Algebras}, Cambridge University Press, (1995).

\bibitem{AS:subcategories}
\textsc{M.~Auslander and S.~Smal{\o}}, \textit{Almost split sequences in
subcategories}, J. Algebra {\bf 69} (2) (1981) 426--454.

\bibitem{AS:gor}
\textsc{M.~Auslander and O.~Solberg}, \textit{Gorenstein algebras and
algebras with dominant dimension at least $2$}, Comm. Algebra {\bf
21} (1993), no. 11, 3897--3934.

\bibitem{BS}
\textsc{S.~Bazzoni and J.~Stovicek}, \textit{All tilting modules are of
finite type},  Proc. Amer. Math. Soc.  {\bf 135}  (2007), no. 12,
3771--3781.


\bibitem{B:gorenstein}
\textsc{A.~Beligiannis}, \textit{The Homological Theory of
Contravariantly Finite Subcategories: Auslander-Buchweitz Contexts,
Gorenstein Categories and (Co-)Stabilizations}, Comm. Algebra {\bf
28}(10), (2000), 4547--4596.

\bibitem{B:3cats}
\textsc{A.~Beligiannis}, \textit{Relative Homological Algebra and Purity in
Triangulated Categories},  J. Algebra {\bf 227} (2000), 268--361.


\bibitem{B:freyd}
\textsc{A.~Beligiannis}, \textit{On the Freyd Categories of an Additive Category},
Homology Homotopy Appl. {\bf 2} (2000), 147--185.

\bibitem{B:mathscand}
\textsc{A.~Beligiannis}, \textit{Homotopy Theory of Modules and Gorenstein Rings},
Math. Scand. {\bf 88} (2001), 5--45.


\bibitem{B:art}
\textsc{A.~Beligiannis}, \textit{Auslander-Reiten Triangles, Ziegler Spectra and
Gorenstein Rings}, K-Theory {\bf 32} (2004), 1--82.

\bibitem{B:cm}
\textsc{A.~Beligiannis}, \textit{Cohen-Macaulay Modules, (Co)Torsion Pairs, and
Virtually Gorenstein Algebras}, J. Algebra Vol. {\bf 288}, No.1,
(2005), 137--211.


\bibitem{BM}
\textsc{A.~Beligiannis and N.~Marmaridis}, \textit{Left Triangulated
Categories Arising from Contravariantly Finite Subcategories}, Comm.
in Algebra {\bf 22} (1994), 5021--5036.

\bibitem{BK}
\textsc{A.~Beligiannis and H.~Krause}, \textit{Thick Subcategories and
Virtually Gorenstein Algebras}, Illinois J. Math. {\bf 52}, No. 2, (2008), 551--562.



\bibitem{BR}
\textsc{A.~Beligiannis and I.~Reiten}, \textit{Homological and homotopical
aspects of torsion theories}, Mem. Amer. Math. Soc. {\bf 188}
(2007), no. 883, viii+207 pp.


\bibitem{BondalKapranov}
\textsc{A.~Bondal and M.~Kapranov}, \textit{Representable functors,
Serre functors, and reconstructions}, Math. USSR-Izv.  {\bf 35}
(1990), no. 3, 519--541.


\bibitem{BVdB}
\textsc{A.~Bondal and M.~Van den Bergh},  \textit{Generators and
representability of functors in commutative and noncommutative
geometry},  Moscow Math. J. {\bf 3} (2003), 1--36.

\bibitem{BuK}
\textsc{A.B.~Buan and H.~Krause}, \textit{Cotilting modules over tame
hereditary algebras},  Pacific J. Math. {\bf 211} (2003), no. 1,
41--59.

\bibitem{Buchweitz:tate}
\textsc{R.-O.~Buchweitz}, \textit{Maximal Cohen--Macaulay modules and
Tate-cohomology over Gorenstein rings}, Unpublished manuscript
(1987), 155 pp.

\bibitem{Buchweitz}
\textsc{R.-O.~Buchweitz}, \textit{Finite representation type and periodic Hochschild
(co-)homology}, Trends in the representation theory of
finite-dimensional algebras (Seattle, 1997),  81--109, Contemp.
Math., 229, Amer. Math. Soc., Providence, RI, (1998).



\bibitem{Chen}
\textsc{X.-W. Chen}, \textit{An Auslander-type result for
Gorenstein-projective modules},  Adv. Math.
{\bf 218}, Issue 6, 2043--2050.

\bibitem{Christensen}
\textsc{L.W.~Christensen}, \textit{Gorenstein dimensions},  Lecture
Notes in Mathematics, {\bf 1747}, Springer-Verlag, Berlin, 2000.
viii+204 pp.

\bibitem{CPST}
\textsc{L.W.~Christensen, G.~Piepmeyer, J.~Striuli, and R.~Takahashi},
\textit{Finite Gorenstein representation type implies simple singularity},
Adv. Math. {\bf 218} (2008), no. 4, 1012--1026.

 \bibitem{WCB}
\textsc{W.~Crawley -~Boevey}, \textit{Locally Finitely Presented
Additive Categories}, Comm. in Algebra, {\bf 22}, (1994),
1644--1674.


\bibitem{ET}
\textsc{P.~Eklof and J.~Trlifaj}, \textit{Covers induced by Ext}, J.
Algebra {\bf 231} (2000), 640--651.

\bibitem{EJ}
\textsc{E.E.~Enochs and O.M.G.~Jenda}, \textit{Gorenstein injective and
projective modules},  Math. Z. {\bf 220} (1995), no. 4, 611--633.

\bibitem{EJ:book}
\textsc{E.E.~Enochs and O.M.G.~Jenda}, \textit{Relative homological algebra}, de Gruyter
Expositions in Mathematics, {\bf 30}, Walter de Gruyter $\&$ Co.,
Berlin, (2000), xii+339 pp.


\bibitem{Happel}
\textsc{D.~Happel}, \textit{On Gorenstein algebras}, Progr. Math. {\bf 95} (1991), 389--404.


\bibitem{HJ}
\textsc{H.~Holm and P.~J{\o}rgensen}, \textit{Rings without a Gorenstein analogue of the Govorov-Lazard Theorem},
Q. J. Math. {\bf 62} (2011), 977--988. 


\bibitem{Iyama}
\textsc{O.~Iyama}, \textit{Higher Auslander-Reiten Theory on maximal
orthogonal subcategories}, Adv. Math. {\bf 210} (2007), 22--50.

\bibitem{Jensen}
\textsc{C.U.~Jensen}, \textit{On the vanishing of $\plim\sp{(i)}$}, J.
Algebra {\bf 15} (1970), 151--166.

\bibitem{JL}
\textsc{C.U.~Jensen and H.~Lenzing}, \textit{Model-Theoretic Algebra with
Particular Emphasis on Fields, Rings, Modules, Algebra, Logic and
Applications}, vol. {\bf 2}, Gordon and Breach Science Publishers,
New York, 1989, xiv+443 pp.

\bibitem{JS}
\textsc{D.A.~Jorgensen and L.M.~Sega},  \textit{Independence of the
total reflexivity conditions for modules}, Algebr. Represent. Theory
{\bf 9} (2006),  no. 2, 217--226.

\bibitem{Keller:orbit}
\textsc{B.~Keller}, \textit{On triangulated orbit categories}, Doc.
Math. {\bf 10}  (2005), 551--581.

\bibitem{KR}
\textsc{B.~Keller and I.~Reiten}, \textit{Cluster-tilted algebras are
Gorenstein and stably Calabi-Yau},  Adv. Math.  {\bf 211}  (2007),
no. 1, 123--151.

\bibitem{KR:acyclic}
\textsc{B.~Keller and I.~Reiten}, \textit{Acyclic Calabi-Yau categories. With an appendix by Michel Van den Bergh.},
Compositio Math.,  {\bf 144}, Issue 05, (2008), 1332--1348.


\bibitem{KZ}
\textsc{S.~Koenig and B.~Zhu}, \textit{From triangulated categories to
abelian categories: cluster tilting in a general framework},  Math.
Z. {\bf 258}  (2008),  no. 1, 143--160.

\bibitem{Krause:lfp}
\textsc{H.~Krause}, \textit{Functors on locally finitely presented
additive categories}, Colloq. Math. {\bf 75} (1998), no. 1,
105--132.


\bibitem{Krause:memoirs}
\textsc{H.~Krause}, \textit{The Spectrum of a Module Category}, Mem. Amer. Math.
Soc. {\bf 149}  (2001),  no. 707, x+125 pp.

\bibitem{Krause:stable}
\textsc{H.~Krause}, \textit{The stable derived category of a Noetherian scheme},
Compos. Math.  {\bf 141}  (2005),  no. 5, 1128--1162.

\bibitem{IK}
\textsc{H.~Krause and I.~Iyengar}, \textit{Acyclicity versus total
acyclicity for complexes over Noetherian rings}, Doc. Math. {\bf 11}
(2006), 207--240.

\bibitem{KS}
\textsc{H.~Krause and M.~Saorin}, \textit{On minimal approximations of
modules}, Contemp. Math., {\bf 229}, Amer. Math. Soc., Providence,
RI, (1998), 227--236.

\bibitem{KS:appl}
\textsc{H.~Krause and O.~Solberg}, \textit{Applications of cotorsion
pairs.}  J. London Math. Soc. (2)  {\bf 68}  (2003),  no. 3,
631--650.

\bibitem{Leuschke}
\textsc{G.J.~Leuschke}, \textit{Endomorphism rings of finite global
dimension}, Canad. J. Math. {\bf 59} (2007), no. 2, 332--342.

\bibitem{LiZhang}
\textsc{Z.-W. Li and Pu Zhang}, \textit{Gorenstein algebras of finite Cohen-Macaulay type},
Adv. Math. {\bf 223} (2010), 728--734.

\bibitem{Orlov}
\textsc{D.O.~Orlov}, \textit{Triangulated categories of singularities
and D-branes in Landau-Ginzburg models},  Proc. Steklov Inst. Math.
2004,  no. {\bf 3} (246), 227--248.

\bibitem{Ringel:good}
\textsc{C.M.~Ringel}, \textit{The category of modules with good
filtrations over a quasi-hereditary algebra has almost split
sequences},   Math. Z. {\bf 208} (1991),  no. 2, 209--223.

\bibitem{Ringel:quasi}
\textsc{C.M.~Ringel}, \textit{The category of good modules over a quasi-hereditary
algebra}, Carleton-Ottawa Math. Lecture Note Ser., {\bf 14},
Carleton Univ., Ottawa, ON, 17 pp., (1992).

\bibitem{Ringel:torsionless}
\textsc{C.M.~Ringel}, \textit{The torsionless modules of an artin algebra},
preprint, (2008).

\bibitem{RT}
\textsc{C.M.~Ringel and H.~Tachikawa}, \textit{QF-3 rings}, J. Reine Angew.
Math.  {\bf 272}  (1974), 49--72.


 \bibitem{RS1}
\textsc{C.M.~Ringel and M.~Schmidmeier}, \textit{The Auslander-Reiten
translation in submodule categories}, Trans. Amer. Math. Soc. {\bf
360} (2008), no. 2, 691--716.

\bibitem{RS2}
\textsc{C.M.~Ringel and M.~Schmidmeier}, \textit{Invariant subspaces of nilpotent linear
operators. I.}, J. Reine Angew. Math. {\bf 614}, (2008), 1--52.


\bibitem{Rouquier:Kth}
\textsc{R.~Rouquier}, \textit{Dimensions of triangulated categories},
Journal of K-theory {\bf 1} (2008), 193--256 and errata, 257--258.

\bibitem{Takahashi}
\textsc{R.~Takahashi}, \textit{Contravariantly finite resolving
subcategories over commutative rings},  Amer. J. Math. {\bf 133} (2011), no. 2, 417--436.

\bibitem{Warfield}
\textsc{R.B. Jr.~Warfield}, \textit{Purity and algebraic compactness
for modules},  Pacific J. Math. {\bf 28} (1969), 699--719.

\bibitem{Yoshino:book}
\textsc{Y.~Yoshino}, \textit{Cohen-Macaulay modules over Cohen-Macaulay
rings},  London Mathematical Society Lecture Note Series, {\bf 146},
Cambridge University Press, Cambridge, (1990), viii+177 pp.

\bibitem{Yoshino}
\textsc{Y.~Yoshino}, \textit{A functorial approach to modules of G-dimension zero},
Illinois J. Math.  {\bf 49}  (2005),  no. 2, 345--367.

\end{thebibliography}
\end{document}